\documentclass[11pt]{article} 
\usepackage{amssymb,amsmath,amsthm}
\usepackage[linktoc=all]{hyperref}
\usepackage{tikz}
\usetikzlibrary{arrows.meta}
\usetikzlibrary{decorations.markings}

\newcommand{\RR}{\mathbb{R}}
\newcommand{\ZZ}{\mathbb{Z}}
\newcommand{\NN}{\mathbb{N}}

\newcommand{\xx}{\mathbf{x}}
\newcommand{\yy}{\mathbf{y}}
\newcommand{\zz}{\mathbf{z}}
\newcommand{\ww}{\mathbf{w}}
\newcommand{\xaa}{\mathbf{a}}
\newcommand{\xbb}{\mathbf{b}}
\newcommand{\xee}{\mathbf{e}}

\newcommand{\xtt}{\mathbf{t}}
\newcommand{\xuu}{\mathbf{u}}
\newcommand{\xvv}{\mathbf{v}}
\newcommand{\xone}{\mathbf{1}}
\newcommand{\LL}{\mathcal{L}}

\renewcommand{\AA}{\mathcal{A}}

\newcommand{\CC}{\mathcal{C}}
\newcommand{\VV}{\mathcal{V}}

\newcommand{\eps}{\varepsilon}
\newcommand{\Ran}{\mathrm{ran}\,}
\newcommand{\Ker}{\mathrm{Ker}\,}
\newcommand{\Tr}{\mathrm{Tr}}
\newcommand{\supp}{\mathrm{supp}}
\newcommand{\const}{\mathrm{const}}
\newcommand{\Pl}{\mathrm{Palm}}
\newcommand{\Patt}{P}

\newcommand{\Prob}{\mathrm{Prob}}  % probability
\newcommand{\prE}{\mathbb{E}} % expectation

\newcommand{\Ga}{\Gamma}

\newcommand{\Knoloops}[1]{K^*_{#1}}%{\overleftrightarrow{K_{#1}}}
\newcommand{\Kwloops}[1]{K^\circ_{#1}}%{\hat K_{#1}}

\newcommand{\parent}[1]{\alpha(#1)} %parent vertex (in trees)

\newcommand{\beq}[1]{\begin{equation}\label{#1}}
\newcommand{\eeq}{\end{equation}}
\newcommand{\bref}[1]{{\bf\ref{#1}}}
\newcommand{\brefA}[1]{{\bf\ref{#1}A}}
\newcommand{\brefB}[1]{{\bf\ref{#1}B}}
\newcommand{\bnref}[1]{{\bf\ref*{#1}}}

\newcommand{\ba}[1]{\begin{array}{#1}}
\newcommand{\ea}{\end{array}}
\newcommand{\mat}[2]{\left[\ba{#1}#2\ea\right]}

\newcommand{\dst}{\displaystyle}
\newcommand{\eqline}[1]{\\[#1ex]\dst}

\newcounter{problem}
\newcounter{subproblem}
\renewcommand{\thesubproblem}{\arabic{problem}\Alph{subproblem}}
\newcommand{\prob}[1]{\refstepcounter{problem}\label{#1}\setcounter{subproblem}{0}  \bigskip\noindent{\bf\arabic{problem}.}}
\newcommand{\subprob}[1]{\refstepcounter{subproblem}\medskip\noindent{\bf\thesubproblem.}\label{#1}}
\newcommand{\trivprob}[1]{\refstepcounter{problem}\label{#1}\setcounter{subproblem}{0}\bigskip\noindent{\bf \arabic{problem}$^\circ$.}}
\newcommand{\trivsubprob}[1]{\refstepcounter{subproblem}\label{#1}\medskip\noindent{\bf \arabic{problem}\Alph{subproblem}$^\circ$.}}
\newcommand{\trivcsubprob}[1]{\refstepcounter{subproblem}\label{#1}\medskip\noindent{\bf \arabic{problem}\Alph{subproblem}$^{\bullet}$.}}
\newcommand{\hardprob}[1]{\refstepcounter{problem}\label{#1}\setcounter{subproblem}{0}
\bigskip\noindent{\bf \arabic{problem}$^*$.}}
\newcommand{\hardsubprob}[1]{\refstepcounter{subproblem}\label{#1}\medskip\noindent{\bf \arabic{problem}\Alph{subproblem}$^*$.}}
\newcommand{\thprob}[1]{\refstepcounter{problem}\label{#1}\setcounter{subproblem}{0} \bigskip\noindent{\bf \arabic{problem}$^\bigstar$.}}
\newcommand{\thsubprob}[1]{\refstepcounter{subproblem}\label{#1}\medskip\noindent{\bf \arabic{problem}\Alph{subproblem}$^\bigstar$.}}

\newcommand{\sol}[1]{\bigskip\noindent s{\bf \ref{#1}. }}
 
\newcommand{\solL}[2]{\bigskip\noindent s{\bf \ref{#1}#2. }}

\renewcommand{\thesubsection}{\arabic{subsection}}

\newcommand{\separline}{\medskip\centerline{\rule{0.26\textwidth}{0.5pt}\hspace{0.5em} \raisebox{-0.65em}{$\tilde{}$\; $\tilde{}$}
\hspace{0.5em}\rule{0.26\textwidth}{0.5pt}}\medskip}

\title{Minimization of the sum under product constraints
\\[\bigskipamount] \large(Problems and theorems)}
\author{Sergey Sadov\footnotemark[1]\footnote{E-mail: serge.sadov@gmail.com}}
\date{}

\begin{document}
\maketitle
\thispagestyle{empty}

\begin{abstract}
We systematically explore a class of constrained optimization problems
 with linear objective function and constraints that are linear combinations of logarithms of the optimization variables. Such problems can be viewed as a generalization of the inequality between the arithmetic and geometric means.

The existence and uniqueness of the minimizer is proved under natural assumptions in the general case.

We study in detail special subclasses where the set of constraints is described in  combinatorial terms (oriented graphs, rooted trees). 
In particular, given a directed, strongly connected graph, we seek to minimize the total of all arc values under cyclic product constraints. 

We obtain some estimates and asymptotics for the minimum in problems with given (large) number of variables. 

Also in this context we revisit an asymptotical result known as J.~Shallit's minimization problem. 
%from 1994 SIAM Reviews. 
%discussed in the last section
%together with some generalizations.

The material is presented in the form of a problem book.
Along with problems that constitute main theoretical threads, there are many exercises, some mini-paradoxes, and a touch of numerical 
methods.


\medskip
MSC primary: 
26D20,  % Other analytical inequalities
secondary:
05C35, % Graph theory, Extremal problems 
%05C38 % Paths and cycles
15A99, %Linear algebra- miscellaneous
%
%
%49-XX CALCULUS OF VARIATIONS AND OPTIMAL CONTROL;
%OPTIMIZATION [See also 34H05, 34K35,
%49Jxx Existence theories
%49J21 Optimal control problems involving relations other than differential
%equations
%49Kxx Optimality conditions
%49K21 Problems involving relations other than differential equations
%
%39-XX DIFFERENCE AND FUNCTIONAL EQUATIONS
39A20, %: Multiplicative and other generalized difference equations, e.g. of Lyness type
%49M05, %Methods of successive approximations->Methods based on necessary conditions
49K99 %: Necessary conditions and sufficient conditions for optimality->None of the above, but in this section
%

%
% 90-XX OPERATIONS RESEARCH, MATHEMATICAL PROGRAMMING
%90C35, %Programming involving graphs or networks
%90C27. %Combinatorial optimization

\end{abstract}

\newpage

\tableofcontents

\newpage

\renewcommand{\thesection}{}
\section{Introduction}
\label{sec:intro}
%\addcontentsline{toc}{section}{Introduction}
% But I don't want the page number in this line
%\addtocontents{toc}{\protect\contentsline{section}{\protect{}Introduction}{}{}}

\addtocounter{subsection}{10}
\newcounter{ssecaux}
\setcounter{ssecaux}{0}
\renewcommand{\thesubsection}{\arabic{ssecaux}}

\refstepcounter{ssecaux}
\subsection{The class of optimization problems}

This work is concerned with optimization of a linear objective function 
\beq{sum1}
 \sum_{i=1}^n w_i x_i \to \min
\eeq
under product constraints
\beq{cons_prod}
\prod_{i=1}^n   x_i^{A_{ji}} 
=t_j, \quad 1\leq j\leq r,
\eeq
or, equivalently, under log-linear constraints
\beq{cons_loglin}
\sum_{i=1}^n  A_{ji} \log x_i 
=\log t_j, \quad 1\leq j\leq r.
\eeq
Here $w_i>0$, $t_j>0$ and $A_{ij}\in\RR$ are given constants; 
%\xx=
$(x_1,\dots,x_n)$ is the unknown vector.  

The familiar inequality between the arithmetic and geometric means (AM-GM) is the simplest representative of problems of this class: here, the minimum of the sum $x_1+\dots+x_n$ is sought assuming the product $x_1\dots x_n$ is known. 
%Unsurprisingly, the AM-GM inequality 
%appears as a particular case of most special subclasses of the problem \eqref{sum1}--\eqref{cons_loglin} described here.
% and in many elements of constructions and proofs.

Putting $y_i=\log x_i$, one can, of course, interpret the problem as minimization of the exponential sum, $\sum w_i\exp y_i$, under linear constraints. We take this point of view occasionally.

\refstepcounter{ssecaux}
\subsection{Quick orientation}

%This is a brief orientation to help the reader to decide 
%where to proceed.
%whether to proceed sequentially or to jump. 

%\smallskip
The material is presented in the format of a problem book. In the first part, the facts (problems) are stated, with necessary definitions, and the second part contains solutions and comments. 
Some of the problems are exercises and some are part of a theoretical thread. We further comment on this
%, as well as about everything else, 
in ``A more detailed guide'' later on. % section.

Inspected in the first approximation, each (sub)section of the first part is devoted to more or less
one little theory and related ruminations. 
%And each one can be associated with one of three rather distinct types of contexts.  
With the exception of Section~2, various instances and disguises of problems of type \eqref{sum1}--\eqref{cons_prod} are the object of those
``little theories''.

Section~1 deals with the general problem, when there is no preset structure for the data. 

Section~2 basically fixes terminology pertaining to graphs, which is subsequently used throughout to define
special subtypes of the general problem. 

In Sections~3--5 and 7, some of such special subtypes are studied. 

Section~6 is concerned with questions of
a somewhat different flavor: instead of treating individual minimization problems, we ask about extremal values of the minimum over certain {\em classes}\ of such problems described in appropriate combinatorial terms.

Somewhat hidden among trivia and digressions is a theoretical thread through Sections~2, 3, 5, and 6 (see ``Origins and motivations''), culminating in the asymtpotic formula of Problem~\brefA{prob:mingraph-local-lb-homo}.
%{prob:minextr-graph-global}
 
It's time to start outlining all of the \textit{what}s, \textit{why}s, and \textit{where}s, but the author sees two particularly impatient readers here, Aaron and Beth, and wants to 
% let them go
serve their interests first.

Aaron is a busy person involved in a serious applied research and not looking for exotic distractions.
% is the last thing he needs.
I venture to suggest one result worth his attention here: the uniqueness theorem, Problem~\bref{prob:uniqsol} in Section~1  ---
a property that makes the class of problems \eqref{sum1}--\eqref{cons_prod} rather distinguished.

\smallskip
All readers note: the uniqueness theorem is the pivotal result; without it little would be left of the whole edifice. 

\smallskip
Beth cares more about tangible facts and curiosities than about the general and mainstream. Being a dedicated problem solver, she expects a collection of problems to be a challenge to have fun with. 
So Beth may decide to leap right away.  
%,  Why not to start with  to get your hands dirty over a concrete calculation? (Reading a few definitions is necessary, alas.)
Destinations to try are many: from the little paradoxes
\bref{prob:notmonotone-wrt-products}--\bref{prob:consineq}\ and olympiad-type inequalities \brefA{prob:pnorm}, \brefA{prob:ineq-circuits-arcs}, to that very theorem \bnref{prob:uniqsol}, which is to test your linear algebra and calculus skills in combination, as is \bref{prob:harmonic-const-ratio}, to 
%more challenging \bref{prob:Shallit} (if you haven't seen it elsewhere) and 
my ``time-stamped'' \bref{prob:extrval-teaser},
not yet antique at the time of this writing. 
If some of them catch your interest but resist an effort, it is never late to take a slower start or to check Part II.

Another idea, particularly suitable if you want to get a feel of a concrete calculation pertaining to the main theme, is to leap to problem \bref{prob:ex-quasibasis}. Some browsing backwards will be needed to understand the notation and terminology, but it's a quick way into the middle of things.  
%author-dreamed 
%extreme --- down to (possibly trivial and boring) exercises at the beginning of Section~1.

%\smallskip
%There is a thematical index of problems at the end of Introduction. 

\smallskip
The author hopes that this work can be useful as a reference by researches who might encounter a problem of the described type or, in accordance with its layout, as a source for seminars and students' projects. %Naturally, 
The author strived to uphold the style and standards set by 
the classics (see e.g.\ Preface to \cite{PolyaSzego_1927}) to the extend to his ability. The level of difficulty is far from uniform: simple exercises are deliberately dispersed among more difficult problems to help the readers to assert trivialities and clear up misconceptions.
%Aaron and Beth may still stay with us.

The ``little theories'' found here may indeed seem exotic, but an effort is made to avoid them looking esoteric. 
The author would be particularly happy if 
meaningful applications or connections of the exposed threads were reported by the readers.

The requiremens to reader's background are modest: elementary calculus (Lagrange multipliers) and linear algebra will suffice. While the terminology of graph theory is used extensively and traces of the probability theory appear occasionally, no special knowledge of those branches of mathematics is assumed.

\refstepcounter{ssecaux}
\subsection{Origins and motivations}
This work resulted from
%came into being as a byproduct of 
observations made by the author while browsing through various analytical inequalities and techniques of proof. The initial point was Harold Shapiro's cyclic inequality \cite[Ch.~XCVI]{Mitrinovic_1973}.
In Section~3.1 of the book \cite{Finch_2003} next to Shapiro's problem, a tangential mention is made of a peculiar inequality due to J.~Shallit.
(The Shapiro and Shallit inequalities are described here: %\bref{prob:Shapiro}, 
\brefA{prob:Nesbitt},
\bref{prob:Shallit} and comments.)
Common to both problems is a linear (in the number of variables) asymptotic behaviour of the lower bound with a numerical coefficient. It is not however this banal  similarity that eventually brought this work into being but a realization, separately and differently for each of the two problems, of a presence in the background of certain sums to be minimized under product constraints.  

Cyclic inequalities similar to Shapiro's are tricky.
So it makes sense to play around and try to find variations that are manageable. A problem of finding the lower bound for the cyclic sum 
$$
 \sum_{i\!\!\!\mod n}\frac{x_i}{\min(x_{i+1},\dots,x_{i+k})}
$$
is simple, as is its analog with function `$\max$' in the denominators. (It is much harder with function `$\Sigma$': Diananda's generalization of Shapiro's sum.)

What if the denominator of the $i$-th fraction of the sum contains minimum (or maximum) of some other set of variables? Does the problem remain trivial or will it become, as a rule, hopelessly messy? Investigating this question, the author was led to the structured sums of Sections~3 and the extremal problem of Section~6.2, which in turn required some constructions of Section~5. Full solution is the subject of a forthcoming paper
\cite{Sadov_Graphic_inequalities}. We come very close and the technical preparation is essentially complete. However,
the motivating problem would be buried here; in author's view it deserves a separate title.
% here it would seem off-topic.
% shielded forever by the title.
%; also out of considerations of size and balance. 

%It is not so easy however to explain what unites them here. 

It is easier to describe the relevance of Shallit's problem. It can be put in the form \eqref{sum1}--\eqref{cons_prod}. In \cite{Shallit_sol_1995}, having stated the critical point system, the authors remark: ``It may be regarded as intractable except for very small $n$''. 
The present author had realized that the nonlinear system for the minimizer is not as hopeless as may seem at first sight; for instance, there is enough structure to easily prove the uniqueness of a positive solution by linear algebra means --- and the proof works for the whole class of problems \eqref{sum1}--\eqref{cons_prod}. 

\iffalse
This work began as a research paper aiming at results related to a ``perversed'' Shapiro's problem as described that now take place in Section~6.
In view of many small lemmas, short proofs, remarks, examples and digressions, it seemed an attractive idea to present the material as a collection of problems.
%diagram of main theor thread
\fi

Having two 
% rather ((1,1),1)
nontrivial albeit %, arguably, 
obscure problems and a simple general theorem seemed sufficient to undertake the effort of building a common roof.

Optimization problems with linear constraints and a quadratic objective function need no advertisement, while
the problem of minimization of an exponential sum under linear constraints, being not-too-distant relative graphically, seems to have been largely ignored in the literature, except, of course, the fundamental particular case, the AM-GM inequality. The author sees a plausible explanation in the apparent lack of applications, but part of it can be the belief in general intractability of the critical point system, which we hope to dispel.  

The theory would be much more useful and would have a multitude of connections with classical inequalities (power means, elementary symmetric functions, Carleman's inequality, etc.), if it allowed at least one linear constraint in addition to the product constraints. 
%However, the author at this time is not able to formulate reasonable general assumptions for such theory. 
However, the uniqueness of an extremizer is lost in such generality as simple considerations show. For instance, in the problem of finding extrema of the function $f(x,y)=x$ in the two-variable problem with constraints $x+y=3, xy=2$ and no free variable, there are two critical points already. (If unconvinced by this degenerate example, see Problem \bref{prob:prod-lin-constraints}.)
%Such a (highly desirable) theory probably would have to include conditions of majorization type.

As mentioned earlier, the problems of Shapiro and Shallit keep in focus some type of asymptotic behaviour of the minimum as a function of problem's ``size''. 
While this is not the feature upon which our unification is built, we are mindful of asymptotical questions.
Here and there (especially in Sections~5,~6) we aim to explore the behavior of the minimum as a function of ``size'' parameter in families of problems 
of the type \eqref{sum1}--\eqref{cons_prod} with a defined structural pattern. (The simplest prototypical case is, again, the AM-GM inequality: assuming that $\prod_{i=1}^n x_i=1$, we have $\min \sum_{i=1}^n x_i=n$, where the ``size'' is $n$.)

We also attempt to view the structural pattern as a probability measure space and to pose a few questions about ``typical'' behavior of the minimum.
The author is the least satisfied with the elaboration of this line among all else; however, one result, Problem~\bref{prob:mean-min-fungraph}, is clear-cut and related to a well-known problem of graph theory. 

Since the solution of a nonlinear system is always on the stage, a few problems touch upon algorithmic and numerical aspects, not pretending to any depth.
Complexity of computing the minimum (or the minimizer) is not discussed at all.

Besides the threads and problems prompted by the outlined considerations, there are some simply pulled out of a hat, such as \bref{prob:countcycles-Kn} and Section~4.3.

Many small problems are spread over: exercises, examples, counterexamples. Quite a few are sourced from this author's errors and fallacies; they are here to keep the reader alert. Early examples are \bref{prob:notmonotone-wrt-products} and \bref{prob:consineq}.

\refstepcounter{ssecaux}
\subsection{A more detailed guide}

References to problems are typed in {\bf boldface} with no parentheses, so no confusion with references to numbered formulas should arise. 

Each problem is either a theoretical fact (theorem, lemma) or an exercise or a problem in the proper sense (a little or not so little challenge). 
Of course, there is no strict boundary between these categories, though for some problems the classification  is obvious. 

The order of the problems within a section is chosen, generally, to facilitate the development of the current branch of theory. Definitions are followed by exercises to help
the reader to get accustomed to the terminology. 
Examples are given either before the corresponding theoretical facts (to warm up for a proof of the latter) or after; counterexamples showing relevance of certain conditions and assumptions are usually placed afterwards.  

Oftentime problems are gathered into clusters. The ``leading'' or  ``master'' problem of the cluster sets the main theme; further problems in the cluster have the same number followed by a letter (like {\bf 23A}). Often they are corollaries of the master problem or subordinate facts (examples, counterexamples), but sometimes a lettered problem states the ultimate result. 

For problems marked with a  circle (like $\mathbf{2^\circ}$), no solution or comment is given in Part II. Such problems can be trivial exercises
or they require routine (possibly not too quick) calculation (like %\bnref{prob:Jacobian-AMGM}
\brefA{prob:Jacobian_for_minimizer}). 
In a few cases, a circled problem is hardly simple on its own, but being a slight variation of nearby problems
should present no difficulty in the context
(like %\bnref{prob:ineq-circuits-nodes}
\brefA{prob:countcycles-Kn}).

Some problems of a similar nature are marked with a filled circle (like $\mathbf{9B^{\bullet}}$), which means there is a comment in Part II. 

%[About half of the problem are marked with a circle, %they are given without solution. Such problems can be 
%either trivial (like \ref{prob:relaxcons1}) or solved by routine calculations (like \ref{prob:Jacobian-AMGM}) or be slight variations of nearby problems

Two problems, \bref{prob:uniqsol} and 
\brefA{prob:mingraph-local-lb-homo}, %{prob:minextr-graph-global}, 
are starred:
to mark what are, in author's opinion, 
the pivotal result of this work and the final result of the main theoretical thread.

A few items are marked with asterisk  ($^*$); it means that the author does not have a solution. There will be a comment on most of them in Part II.  

In a collection of small, numerous and diverse facts like this one, the views about their relative 
importance are necessarily subjective and the author prefers to leave the ultimate judgement with the reader. 

\smallskip
Below is a brief overview of the contents.

\smallskip
In Section~1 we study general properties of problems of type \eqref{sum1}--\eqref{cons_prod}, which include: basic transformations, dependence on parameters, making use of symmetries, a numerical approach, but first and foremost --- existence and uniqueness of a (the) minimizer. The section concludes with two examples of optimization problems of similar but wider classes to compare. 

It is possible that a portion of this content can be transferred to the case where the logarithmic function in the log-additive constraints is replaced by some other concave function. The author did not explore matters in that direction.

\smallskip
Section~2 is ancillary; its purpose is to set the preliminaries for combinatorial description of special types of constraints. It introduces the terminology of graph theory with an accent to fit our context. Problems here are mostly exercises to let the reader get accustomed to the terminology. Even so, some earlier ones are not quite sleeping pills. To taste the accent, compare the version of Euler's formula \brefB{prob:cycle_basis_linalg}
%{prob:Euler_formula} 
with one you know.

\smallskip
Section~3 introduces problems with constraints written as products of values of graph's arcs along some cycles.
Manipulations and reductions are discussed, but there is no deeper theory. Problems are mostly simple exercises; an exception is \bref{prob:mean-min-fungraph}, closely related to a well-known enumeration problem of graph theory.

\smallskip
Section~4 is a digression before more theoretically charged sections~5 and 6. 
A look at quotient sums and cyclic inequalities in Subsection~4.1 will strengthen reader's ability to recognize relevant types of problems, which may not be immediately obvious. Also observe the interplay between similar but not identical classes of problems: whether the object of optimization is a $0$-cochain (node values) or $1$-cochain (arc values).
%An interplay between one way and the other will take place in Sections~5 and~6. 

Sums of quotients introduced here will appear again in Section~5 and Shallit's problem of Section~7 has some visual resembance to Shapiro's; otherwise this material is not referred to in the sequel. 

The fundamental cause of difficulty in problems like Shapiro's \bnref{prob:Shapiro} is that the symmetric extremizer is no longer the global minimizer.
The present technique doesn't help, but the lesson is not entirely negative: it reasserts the value of uniqueness. 

Subsection~4.3 stands on its own;
its contents may be viewed as a peculiar generalization of the AM-GM inequality. 

\smallskip
In Section~5 the underlying combinatorial structure 
is a node-weighted tree (or forest). The corresponding optimization problem can also be viewed as a very special case of the problem of Section~3. On the other hand,
an extremal problem for arc-weighted graphs considered in Section~6.2 is reduced in \bref{prob:mingraph-to-tree} to a problem for node-weghted trees; this is why the ``tree sums'' are introduced. 

A notable particular case occurs when there is at most one node of indegree greater than 1; the optimization problem in this case is conveniently
reformulated in terms of node-weighted trees.

Towards the end of Section~5 we discuss computational matters for the tree problem.

\smallskip
The goal of Section~6 is to solve the extremal problem: to find the minimum value of the sum of arc values 
%$\sum x_i$ 
in a graph over the class of all strongly connected graphs with a given number of arcs (representing the optimization variables) and, possibly, nodes (which corresponds to specifying the number of constraints) under the requirement that all cycle products be equal to 1. 

The solution of an individual optimization problem
of the form \eqref{sum1}--\eqref{cons_prod}
amounts to solving a nonilnear system of algebraic equations, which in general is not possible to
do in a closed form. However the extremal problems
of Section~6 allow for an analytical treatment. 

%\bref{prob:mintree-to-min-special-ge}

\smallskip
Finally, in Section~7 we discuss Shallit's problem and some of its analogs, again using the notions introduced
in Sections~2 and~3. 

\clearpage

\newcounter{secaux}
\renewcommand{\thesection}{\Roman{secaux}}

\addtocounter{secaux}{1}
\section{Problems and theorems}
\label{part:problems}
\setcounter{subsection}{0}
\renewcommand{\thesubsection}{\arabic{subsection}}
%\addcontentsline{toc}{section}{Part I. Problems and theorems}
%\addtocontents{toc}{\protect\contentsline{section}{\protect{}Part I. Problems and theorems}{}{}}

\subsection{General 
case:
%properties of the minimizer: 
existence, uniqueness, dependence on data}
\label{sec:general}

Vectors with components $x_i$ (or $w_i$ etc.) are written as $\xx$ (resp.\ $\ww$ etc.), the dimension is defined by the context. We write $\xx=(x_1,\dots,x_n)$ if $\xx$ is but an ordered list of variables and $\xx=[x_1,\dots,x_n]^T$
if $\xx$ is to participate in matrix calculations as a column vector.

We will use the symbol $\xone$ (boldface one) to denote the vector with all components equal to $1$ and dimension determined by context.

The standard Euclidean inner product in $\RR^n$ is denoted
$\langle \xvv,\xvv'\rangle=\sum v_i v'_i$.

A one-dimensional vector, say, $\xx=(x)$, will be occasionally written simply as $x$.   
The vector with components $\log t_i$ is denoted $\log\xtt$. 
We write $\xx>0$ (resp.\ $\xx\geq 0$) if $x_i>0$ (resp.\ $x_i\geq 0$) for all $i$. If $\xx$ and $\xx'$ are vectors of the same dimension,
we write $\xx\geq\xx'$ if $x_i\geq x'_i$ for all $i$.

The minimization problem \eqref{sum1}, \eqref{cons_loglin} is succintly written as 
\beq{basicproblem}
\langle \xx,\ww\rangle\to\min, 
\qquad 
A\,\log\xx=\log\xtt.
\eeq
It is always assumed that $\ww>0$, the number of constraints $r>0$, $\xtt>0$,  and the constraints are independent 
(that is, the rows of the matrix $A$ are linearly independent). We allow a possibility that $x_i=0$ for some $i$ provided $A_{ji}=0$ for all $j$;
in this case we say that $x_i$ is a \emph{nonessential}\ variable, otherwise $x_i$ is an \emph{essential}\ variable. 
If there is a linear combination of constraints with positive coefficients at all essential variables, then we say that the system of costraints
is \emph{compact}. 

A vector $\xx\geq 0$ is \emph{admissible}\ if it satisfies the given constraints.

We set off cautiously and denote
\beq{inf1}
f(\ww,\xtt,A)=\inf_{\xx}\langle \xx,\ww\rangle, \;\; \mbox{\rm where $\xx\geq 0$ satisfies $A\,\log\xx=\log\xtt$}.
\eeq
A \emph{minimizer}\ %(if it exists) 
is any vector $\xx$ for which $\langle \xx,\ww\rangle=f(\ww,\xtt,A)$ satisfying the constraints. 
In formal considerations it is convenient to allow empty set of constraints, in which case the minimum is obviously zero.
We write $f(\ww,\xtt,-)=0$. By definition, the empty constraint system is compact.

\separline

Exercises {\bf 1}--{\bf 10} %\ref{prob:dynprog}  
are to set a comfortable stage for the sequel.

%---------------------------------------------------------

\prob{prob:consfeasible}
The system of constraints \eqref{cons_loglin} (or \eqref{cons_prod}) is always feasible (that is, it has a solution).  
Consequently, always $f(\ww,\xtt,A)<\infty$. 

\prob{prob:noncompact-example}
Give an example with a non-compact system of constraints.  

\prob{prob:change-weight}
(Change of weights). If the pairs of vectors $\ww$, $\xtt$ and $\ww'$, $\xtt'$ satisfy the relation   
$\log \xtt+A\log\ww=\log\xtt'+A\log\ww'$, then 
$f(\ww,\xtt,A)=f(\ww',\xtt',A)$.  
In particular, $f(\ww,\xtt,A)=f(\xone,\xtt',A)$, where $t_j'=t_j\prod_i w_i^{A_{ji}}$.

\trivprob{prob:dilatation}
(Dilatation). Let $\xaa=A\xone$. (That is, $a_j=\sum_i A_{ji}$.)
If the vectors $\xtt$ and $\xtt'$ are related by   
$\xtt+k\xaa =\xtt'$, then 
$f(\ww,\xtt',A)=e^k f(\ww,\xtt,A)$.  
In particular, a one-constraint problem reduces to a problem with product value $1$. 

\trivprob{prob:relaxcons1} 
(Relaxation of constraints). If $\hat A$ is a submatrix of $A$ obtained by deletion of some rows and the vector $\hat\xtt$
is obtained from $\xtt$ by deletion of the corresponding components, then
$f(\ww,\hat\xtt,\hat A)\leq f(\ww,\xtt,A)$. 

\trivprob{prob:monotonicity} 
(Monotonicity with respect to weights.) If $\ww'\geq\ww$, then $f(\ww',\xtt,A)\geq f(\ww,\xtt,A)$.
(For a quantitative version see \bref{prob:Jacobian_for_minimizer}(c).)

\subprob{prob:concave-weights}\ 
(Concavity with respect to weights.) If $0<\alpha<1$ and $\ww=\alpha\ww'+(1-\alpha)\ww''$, then 
$f(\ww,\xtt,A)\geq \alpha f(\ww',\xtt,A)+(1-\alpha)f(\ww'',\xtt,A)$.

%\subprob{prob:weightconcavity}
%(Concavity with respect to weights.)
%Let $\{p_i\}_{i=1}^k$ be a probability distribution.
%Prove that if $\ww=\sum_1^k p_i\ww^{(i)}$, then
%$$
%f(\ww,\xtt,A)\geq\sum_{i=1}^k p_i f(\ww^{(i)},\xtt,A).
%$$

\prob{prob:notmonotone-wrt-products}
(No monotonicity with respect to the right-hand sides.)
Show that it is possible that  all $A_{ji}\geq 0$, $\xtt'\geq \xtt$ yet $f(\ww,\xtt',A)< f(\ww,\xtt,A)$.

\medskip
(One can get a numerical feel for \bref{prob:relaxcons1} and \bref{prob:notmonotone-wrt-products} in Exercise  \bref{prob:ex-quasibasis}.)

\prob{prob:consineq} 
(Constraints in the form of inequalities.) Define
$$
f_{\geq}(\ww,\xtt,A)=\inf_{\xx}\langle \xx,\ww\rangle, \;\; \mbox{\rm where $\xx\geq 0$ satisfies $A\,\log\xx\geq \log\xtt$}.
$$
and similarly define $f_{\leq}$. Then $f_{\geq}(\ww,\xtt,A)\leq f(\ww,\xtt,A)$, $f_{\leq}(\ww,\xtt,A)\leq f(\ww,\xtt,A)$,
and it is possible in both cases that the inequality is strict even if all $A_{ji}>0$.

\prob{prob:zeroweight} 
(Degeneration of a weight). Suppose the variable $x_n$ is essential. 
Eliminating $\log x_n$ from all but one (say, $r$-th) constraints (in the usual linear algebra sense), we will assume that $A_{jn}=0$ for $j=1,\dots,r-1$. 
Let $\hat A$ be the $(r-1)\times (n-1)$ truncation of the matrix $(A_{ji})$ with shortened ranges of indices: $j=1,\dots,r-1$, $i=1,\dots,n-1$.
Denote by $\hat\ww$ the vector $\ww$ with $n$-th component deleted and define the vector $\hat\xtt$ likewise.
Suppose all data of the problem except for the value of $w_n$ are kept constant.
Then $\lim_{w_n\to 0} f(\ww,\xtt,A)=f(\hat\ww ,\hat\xtt,\hat A)$.

\trivprob{prob:dynprog}\ (A dynamic programming approach.) 
Consider a problem \eqref{basicproblem} with $\xx,\ww\in\RR^n$ and $\xtt\in\RR^r$. 
%and a compact system of constraints. 
Suppose $i\in\{1,\dots,n\}$ is an index such that the system of constraints does not determine the value of $x_i$ uniquely. Let $A^{(i)}$ be the matrix obtained by deletion of $i$-th column (corresponding to the variable $x_i$) from $A$ and let $\ww^{(i)}\in\RR^{n-1}$ be the vector obtained from $\ww$ by deletion of the $i$-th component. Let
$\xtt^{(i)}(s)\in\RR^r$ be the vector with components
$t^{(i)}_j(s)=t_j s^{A_{ji}}$.
Then
$$
 f(\ww,\xtt,A)=\inf_{s\in\RR^+}\left(\frac{w_i}{s}+
f(\ww^{(i)},\xtt^{(i)}(s),A^{(i)})\right).
$$

\separline

Problems {\bf 11}--{\bf 12} state convexity inequalities of Jensen and H\"older type.

\prob{prob:logconvex}\ 
(Log-convexity with respect to constraints.) 
%(Analog of Jensen's inequality.)
%Suppose $0<\alpha<1$ and $\beta=1-\alpha$. Let vectors
%$\xtt$, $\xtt'$ and $\xtt''$ be related by
%$\log\xtt=\alpha\log\xtt'+\beta\log\xtt''$.
%Then
%$$
%f(\ww,\xtt,A)\leq f(\ww,\xtt',A)^{\alpha} 
%f(\ww,\xtt'',A)^{\beta}.
%$$
Suppose $\{p_i\}_{i=1}^k$ is a probability distribution
and the vectors $\xtt^{(1)},\dots,\xtt^{(k)}$ and $\xtt$ are related by $\log\xtt=\sum_{i=1}^k p_i \log\xtt^{(i)}$.
Then
$$
\log f(\ww,\xtt,A)\leq \sum_{i=1}^k p_i \log f(\ww,\xtt^{(i)},A) .
$$

\trivprob{prob:pnorm}\ (Minimization of power sums.)
Let $p>0$. Put $f_p(\ww,\xtt,A)=\inf_{\xx}\langle \xx^p,\ww\rangle$ under the same constraints
as in \eqref{inf1}. Then
$$
 f_p(\ww,\xtt,A)=f(\ww,\xtt^p,A).
$$

\subprob{prob:pnorm-convexity}\ (Log-convexity of minimum $p$-norm.)
For $0<p_0<p_1$ and $0<\theta<1$ let $p_{\theta}$ be defined by $1/p_{\theta}=\theta/p_0+(1-\theta)/p_1$. 
Prove that
%The inequality
\beq{log-convex-pnorm}
 f_{p_\theta}(\ww,\xtt,A)^{1/p_\theta}\leq f_0(\ww,\xtt,A)^{(1-\theta)/p_0} f_1(\ww,\xtt,A)^{\theta/p_1}.
\eeq
%holds.

\trivcsubprob{prob:psum-convexity}\ (Log-convexity of minimum $p$-sum)
Prove that if the function $u\mapsto F(1/u)$ is log-convex, that is, $\log F(1/u_\theta)\leq
(1-\theta)\log F(1/u_0)+\theta \log F(1/u_1)$ for $0<u_0\leq u_1$
and $u_\theta=\theta u_0+(1-\theta)u_1$, then the function
$p\mapsto pF(p)$ is also log-convex. 
Conclude that the function $p\mapsto f_p(\ww,\xtt,A)$
is log-convex.

\separline

The next set of problems, up to \bnref{prob:nonuniqsol-polynomial},
deals with existence and uniqueness of the minimizer.
An important corollary of the uniqueness, \bnref{prob:symmetry}, concerns problems with symmetry.

\prob{prob:critpteqs} (Necessary condition for extremum).
A minimizer $\xx$, if it exists, is a solution of the system of $n+r$ equations
\beq{critpteqs}
\ba{l}\dst
w_i x_i=\sum_{j=1}^r \lambda_j A_{ji}\quad (i=1,\dots,n),
\eqline{2}
\sum_{i=1}^n A_{ji}\log x_i=\log t_j \quad (j=1,\dots,r)
\ea
\eeq 
for the unknowns $x_i\geq 0$ ($i=1,\dots,n)$ and $\lambda_j$ ($j=1,\dots,r$).

\smallskip
{\em Geometric interpretation}: the set of critical points is the intersection of the solution set of the constraints 
with range of the matrix $\tilde A^T$, where $\tilde A_{ji}=A_{ji}/w_i$.

\bigskip
The next problem is standard, but we propose
to employ its corollary \bnref{prob:pnorm-convexity}
%\bref{prob:AMGMsmall}  
to establish the existence of minimizer. 
%in the case of compact system of constraints. 

\prob{prob:AMGM}
(Weighted AM-GM inequality). 
Let $r=1$, $A_{1i}=\rho_i$ ($i=1,\dots,n$), $\xtt=(t)$ (one-dimensional vector).
Put $\rho=\sum_1^n \rho_i$. Then 
$$
f(\ww,t,A)=t^{1/\rho}\prod_{i=1}^n\left(w_i \frac{\rho_i}{\rho}\right)^{\rho_i/\rho}=t^{1/\rho} f(\ww,1,A).
$$
The unique minimizer $\xx$ is given by
$$
 x_i w_i=\frac{\rho_i}{\rho} f(\ww,t,A).
$$

\subprob{prob:AMGMsmall}
Corollary. (Effect of a ``small'' component in the vector $\xx$.) %; to be used in \bref{prob:existminimizer}).
In the same setting, suppose $n>1$, $\rho_i>0$. 
Prove that
$$
\langle\xx,\ww\rangle\geq C_i x_i^{-\rho_i/(\rho-\rho_i)},  
$$   
with $C_i=f(\hat\ww^{(i)},t,\hat A^{(i)})$, where
the vector $\hat\ww^{(i)}$ is obtained from $\ww$ by deletion of the $i$-th component
and $\hat A^{(i)}$ is the matrix obtained from $A$ by crossing out the $i$-th column.
%If $x_1\leq \eps$, then 

\prob{prob:compact-sublevel}
Suppose the set of constraints is compact, all variables are essential, and $c>f(\ww,\xtt,A)$. Then the set $X_c$ of admissible vectors 
$\xx$ for which $\langle\xx,\ww\rangle\leq c$ is compact
and non-empty.

\subprob{prob:almost-compact-means-nonzero}
If there is some nontrivial linear combination of constraints with nonnegative coefficients (for instance, if the set of constraints is compact and nonempty),
then $f(\ww,\xtt,A)>0$.

%\subprob{prob:noncompact-nonzero}
%Exhibit a problem with $f(\ww,\xtt,A)>0$ and non-compact system of constraints.

\prob{prob:existminimizer}\
{\bf Theorem.} A minimizer in \eqref{inf1} 
exists if and only if the system of constraints is compact. 

\thprob{prob:uniqsol}\
{\bf Theorem.} A solution of the system \eqref{critpteqs} with $\xx\in\RR^n_{\geq 0}$ and $\lambda_1,\dots,\lambda_r\in\RR$, if it exists, is unique.
Therefore, in a problem \eqref{inf1} with compact system of constraints there exists a unique minimizer. 

\prob{prob:nonuniqsol-polynomial} (Uniqueness relies on specifics of our problem). 
Suppose all $A_{ji}$ are non-negative integers. Upon elimination of the $x$-variables the system \eqref{critpteqs} 
reduces to a system of $r$ polynomial equation for $\lambda_1,\dots,\lambda_r$ with non-negative coefficients and the right-hand side $\xtt$. %,
%$$
% \prod_{i=1}^n \left(\sum_{k=1}^r A_{ik}\lambda_k\right)^{A_{ji}}=t_j \quad(j=1,\dots,r).
%$$
Give a counterexample that refutes a ``generalization'' of \bref{prob:uniqsol} in this special case: %more general uniqueness assertion:
\emph{if $\mathbf{g}:\RR^r_{\geq 0}\to R^r_{\geq 0}$ is a polynomial map with independent, strictly monotone component
functions, then the system $\mathbf{g}(\xx)=\xtt$ has at most one solution in $\RR^r_{\geq 0}$ for any $\xtt\in\RR^r_{\geq 0}$.}
(Note: for $r=1$ the italicized assertion is true.) 

\prob{prob:symmetry} (Symmetric problem has symmetric minimizer.) If $\sigma$ is a permutation of $n$ symbols
such that $w_i=w_{\sigma(i)}$, $A_{ji}=A_{j,\sigma(i)}$ for $i=1,\dots,n$ and $j=1,\dots,r$,
then the minimizer $\xx$ has the same symmetry: $x_i=x_{\sigma(i)}$.

\trivsubprob{prob:symmetry-reduction} (Reduction by symmetry.) Suppose the given problem is invariant under a permutation group $G$ (that is, for every $\sigma\in G$ the above stated conditions are satisfied). Let $\{\Omega_k,\,k=1,\dots,m\}$ be the $G$-orbits on the index set $\{1,\dots,n\}$. Introduce the $m$-vector $\ww^{G}$
and the $r\times m$ matrix $A^G$:
$$
 w^G_k=\sum_{i\in\Omega_k} w_i,\qquad A^G_{jk}=\sum_{i\in\Omega_k} A_{ji}.
$$
Then $f(\ww,\xtt,A)=f(\ww^G,\xtt,A^G)$.

%\subprob{prob:symmetry-example} Example (after Crux 2938a).
%If $x_i$, $y_i$ are positive real numbers, then
%$$
% \sqrt[n]{(x_1+y_1)\cdots(x_n+y_n)}\geq \sqrt[n]{x_1\cdots x_n}+\sqrt[n]{y_1\cdots y_n}
%$$
%(part b of the problem is a trivial consequence of the AGM)
%
%Solution.
%This is a particular case of the multilinear H\"older inequality. No need to take special advantage of symmetry.

\separline

The rest of this section is more formula-heavy and includes a theoretical {\bf 20} and more practical {\bf 21} methods to approach the minimizer, a quantitative result on the dependence of the minimum on data, and two sample optimization problems from larger classes, to compare with.   

\prob{prob:gradient_flow} (Gradient flow converging to the minimizer).
Consider a problem \eqref{inf1} with compact system of constraints and no non-essential variables. 
Suppose that $\ww=\xone$ (there is no loss of generality, cf.\ \bref{prob:change-weight}). 
Let $P$ be the orthogonal projection onto the row space of the matrix $A$.
% : [if the rows are mutually orthogonal then]
% $$
%  (P\xvv)_i=\sum_{j=1}^r \frac{\langle \xvv, A_{j,\bullet}\rangle}{\langle A_{j,\bullet}, A_{j,\bullet}\rangle} A_{ji}.
% $$ 
Let $\xx_*$ be the minimizer and
$\xx_0$ be any admissible vector. Consider the dynamical system
 $$
  \frac{dx_i}{ds}= (P\xx)_i x_i -x_i^2, \quad i=1,\dots,n,
$$
with initial condition $\xx(0)=\xx_0$. Then $\lim_{s\to +\infty} \xx(s)=\xx_*$.

\prob{prob:newton_method} (Newton's iteration.) 
%
%Construct the tangent plane to the constraint manifold
%at the point $\xx$; find its intersection with the plane $(\Ker A)^\bot=\Ran A^T$.  
%Show that the iterations converge to the minimizer.
%
Suppose that the system of constraints 
in \eqref{inf1} is compact. 
Let $M_y$ be the affine subspace of $\RR^n$ of codimension $r$ defined by the equations
$\sum_{i} A_{ji} y_i=\tau_j$, where $\tau_j=\log t_j$, $j=1,\dots,r$.
Let $N_x$ be the linear subspace of $\RR^n$ of dimension $r$ defined by the parametric presentation:
$x_i=w_i^{-1}\sum_j A_{ji}\lambda_i$, $\forall\lambda_i\in\RR$.

Define a map $\Phi$ from $N_x\times M_y$ to itself as follows.
Let $\xx\in N_x$ and $\yy\in M_y$. 
Put $\eps_i=y_i-\log x_i$. Consider the linear system for $2n+r$ unknowns $\{u_i,v_i,\mu_j\}$
\beq{newton-iter}
\begin{array}{l}
\dst
 w_i u_i=\sum\nolimits_{j} A_{ji}\mu_j,\quad i=1,\dots,n;\\
\dst
 \sum\nolimits_{i} A_{ji} v_i=0,\quad j=1,\dots,r;\\[2ex]
\dst 
 \frac{u_i}{x_i}-v_i=\eps_i,\quad i=1,\dots,n. 
\end{array} 
\eeq

(a) Prove that the system \eqref{newton-iter} is uniquely solvable.

\smallskip
Having the existence and uniqueness of the solution, define $\Phi(\xx,\yy)=(\xx+\xuu,\yy+\xvv)$. Given $\xx^{(0)}\in N_x$ and $\yy^{(0)}\in M_y$,
put $(\xx^{(k)},\yy^{(k)})=\Phi^k(\xx^{(0)},\yy^{(0)})$ and $\eps^{(k)}_i=y^{(k)}_i-\log x^{(k)}_i$.

\smallskip
(b) Prove that there exists a positive constant $C$ (depending on the data $A$, $\xtt$) such that
if $C\|\eps^{(0)}\|<1$, then $C\|\eps^{(k)}\|\leq (C\|\eps^{(k-1)}\|)^2$, $k=1,2,\dots$,
where $\|\eps\|=\max_i |\eps_i|$.
 
%%%%%%%%%%%%%%% 
\prob{prob:Jacobian_for_minimizer}
(Dependence on data: partial derivatives).
(a) Consider the minimization problem \eqref{inf1} with fixed matrix $A$ and varying vectors $\xtt$ and $\ww$. Assuming that the constraints are compact, let $\xx=\xx(\ww,\xtt)$ be the minimizer. Denote by $G$ the $r\times r$ matrix%
\begin{NoHyper}%
\footnote{Note: to compute $G$, as well as $B$ and $Q$ below, %the knowledge of the minimizer $\xx$ is required.}
one needs to know the minimizer $\xx$.} 
\end{NoHyper}
with entries $G_{jk}=\sum_{i=1}^n A_{ji} (w_i x_i)^{-1} A_{ki}$.
Then  
\beq{dxdt}
\frac{\partial x_i}{\partial t_j}=\frac{1}{w_i t_j}\sum_{k=1}^r A_{ki} (G^{-1})_{kj} \quad (i=1,\dots,n;\;\; j=1,\dots,r).
\eeq
(b) Let $Q$ be the orthogonal projector %(with respect to the standard Euclidean scalar product in $\RR^n$) 
onto the subspace $\ker B$ of $\RR^n$,
where $B$ is the $r\times n$ matrix with entries $B_{ji}=A_{ji} (x_i w_i)^{-1/2}$. 
Then
\beq{dxdw}
\frac{\partial x_i}{\partial w_j}=-\sqrt{\frac{x_i x_j}{w_i w_j}}\, Q_{ij} \quad (i,j=1,\dots,n).
\eeq
(c) $\partial f(\ww,\xtt,A)/\partial w_j=x_j$. That is, the minimizer can be found as the gradient of $f(\dots)$ with respect to the weight argument.
%(In particular, we get another proof of uniqueness of the minimizer; note however that the assertion \bref{prob:uniqsol} is logically stronger.)

\trivsubprob{prob:Jacobian-AMGM}\
Check the validity of the above formulas
% Eqs.\ \eqref{dxdt} and \eqref{dxdw} 
in the particular case of Problem \bref{prob:AMGM}.

\prob{prob:noncommutative}\ (A sample of a ``matrix generalization''.) Let us look at the minimization problem \eqref{inf1} in a slightly different way.
Let $B_1,\dots,B_r$ be $n\times n$ diagonal matrices such that the diagonal entries of $B_j$ are
$A_{j1},\dots,A_{jn}$. Let $W$ be the diagonal matrix with diagonal entries $w_1,\dots,w_n$.
Then 
$$
 f(\ww,\xtt,A)=\inf_Y \Tr(W^T e^Y),
$$
where $Y$ runs over the set of diagonal $n\times n$ matrices satisfying the linear constraints $\Tr(B_j^T Y)=\log t_j$, $j=1,\dots,r$. Explore the following situation where the matrices are not all diagonal.

Let $n=2$, $W=\mat{cc}{1 & 0\\0 &1}$, $r=3$, $B_1=\mat{cc}{1 & 0\\0 &-1}$, $B_2=\mat{cc}{0 & 1\\-1 &0}$, 
$B_3=\mat{cc}{1 & 1\\1 &1}$, $t_1=t_2=1$, $t_3=t$.
Find $\inf_Y \Tr e^Y$ as a funciton of $t$ (nondiagonal matrices $Y$ are allowed). Does a minimizer exist? 

\prob{prob:prod-lin-constraints} (A mix of product and linear constraints).
Find the extrema of the objective function
$$
f(u,v,w)=uv+vw+wu,
$$
under the constraints
$$
 uvw=p, \qquad u+v+w=s,
$$
where $p>0$ and $s>0$ are given.
Interpret this problem as a problem of a type similar to
\eqref{sum1}--\eqref{cons_prod} with one additional 
{\em linear}\ constraint.
 
Analyze the outcome in light of the results %\bref{prob:critpteqs}, 
\bref{prob:uniqsol} and \bref{prob:symmetry}.

\subsection{Graphs, cycles, bases}
%and related minimization problems}
\label{sec:graphs-basics}

%\subsubsection{Terminology concerning graphs}

%Digraph, simple digraph, path, circuit, loop, strongly connected, strong component, condensation, final component, tree, forest, poset.

\begin{NoHyper}
By a {\it graph}\ we mean a {\it directed graph},% 
\footnote{In order to define the minimization problem \eqref{mingraph} below, all one needs is the arc-cycle incidence matrix $A$,
in which the orientation of edges (arcs) plays no role. However, in regard to cycles, using directed graphs versus undirected 
amounts in most places to clear-cut formulations versus wordy and/or unnatural.} 
\end{NoHyper}
except in Section~\ref{ssec:harmonic}.
%(2) directed graphs occur naturally when problem \eqref{mingraph-nodes}\ is considered.}
%
A graph $\Ga$ is specified by the set of nodes $\VV(\Ga)$, the set of arcs $\AA(\Ga)$,
and the functions $\alpha$ (beginning of arc) and $\beta$ (end of arc) from $\AA(\Ga)$ to $\VV(\Ga)$. When the graph is fixed, we write simply $\VV$
and $\AA$. 

The graph is \emph{simple}\ if $\alpha(a)=\alpha(a')=v_1$ and $\beta(a)=\beta(a')=v_2$ imply $a=a'$; in a simple graph an arc is uniquely determined
by its ends and we conveniently write $a=(v_1\to v_2)$.
%(The case of non-simple graphs can be reduced to ours by changing weights appropriately.) 
Loops $a=(v\to v)$ are allowed. 

\smallskip
A {\it path}\ %in a graph $\Ga$ 
between nodes $v$ and $v'$ is a sequence of arcs $a_1,\dots,a_k$, $k\geq 1$, where all arcs are distinct and  $\alpha(a_{i+1})=\beta(a_i)$, $i=1,\dots,k-1$, $\alpha(a_1)=v$, $\beta(a_k)=v'$. 

A path with $v=v'$ is called \emph{closed}.

Shifting arcs cyclically along a closed path yields an equivalence relation on the set of closed paths. The equivalence classes are called \emph{circuits}. 

%A \emph{cycle}\ is an equivalence class of closed paths with the same set of arcs that do 
%not pass through the same node twice.
%
\begin{NoHyper}
A \emph{cycle}\ is a circuit that does not pass through the same node twice.% 
\footnote{More precisely, this property is required to hold for {\em some}\ closed path representing the given circuit;
then it holds for {\em all}\ representing paths.}
\end{NoHyper}

Denote by $\CC=\CC(\Ga)$ the set of all circuits in the graph $\Ga$. Let $\CC_0(\Ga)\subset \CC(\Ga)$ be the set of all cycles in $\Ga$.

\smallskip
The cardinality of a finite set $X$ is denoted $|X|$.

The \emph{support}\ of a path, circuit or cycle is the set of arcs in the respective object. 

Any circuit is a concatenation of cycles. 
A cycle is uniquely determined by its support, while a circuit in general is not. 

We will use a slightly abusive notation %like 
$a\in\gamma$ ($a$ an arc, $\gamma$
a circuit) 
and $|\gamma|$ 
instead of the pedantic but clumsy `$a\in\supp(\gamma)$'. 
and `$|\supp(\gamma)|$'.

A circuit $\gamma$ is an {\em Eulerian circuit}\ if $|\gamma|=|\AA(\Ga)|$.

%\smallskip

\separline

%Before we proceed with the optimization problem as such, some more terminology will be introduced
%interspersed with warm-up (or stay-awake) problems \bref{count-circuits-cycles-small-n}--\bref{prob:cycle_basis_linalg}.
%The problems \bref{setcircuits-Eul} and \bref{setcircuits-Ham} below are not connected with our optimization problem but can help the reader to stay awake
%amongst a possibly boring roll of notation and mostly standard terminology. 

We begin with a few warm-up and hopefully entertaining counting exercises that are to facilitate the understanding of cycles and circuits. 
%They are not related to our main line

\trivprob{prob:count-circuits-cycles-small-n} 
Let $\Knoloops{n}$ be the complete directed graph on $n$ nodes without loops
and $\Kwloops{n}$ be the complete directed graph on $n$ nodes with loops.
Thus 
%$\AA(\Knoloops{n})=\{(a\to b), a\neq b\in \VV(\Knoloops{n})\}$, 
$|\AA(\Knoloops{n})|=n(n-1)$,  
and  
%$\AA(\Kwloops{n})=\{(a\to b), a,b\in \VV(\Kwloops{n})\}$, 
$|\AA(\Kwloops{n})|=n^2$. 
Check that $|\CC(\Knoloops{2})|=1$, $|\CC_0(\Kwloops{2})|=3$, $|\CC_0(\Knoloops{3})|=5$, $|\CC(\Kwloops{2})|=6$, and $|\CC(\Knoloops{3})|=9$.

\trivsubprob{prob:count-circuits-vs-supports}\ 
Let $\Ga$ be the graph with one node and 3 loops. Check that $|\CC_0(\Ga)|=3$, $|\CC(\Ga)|=8$, and the number of different supports of circuits in $\Ga$ is $7$.

\prob{prob:ineq-circuits-arcs}
Prove that in any graph $|\CC(\Ga)|< e |\AA(\Ga)|!$.  

\subprob{prob:ineq-circuits-arcs-sharp}
More precisely, if $|\AA|=m$, then $|\CC|< (e+o(1)) (m-1)!$
as $m\to\infty$ and the coefficient $e$ cannot be replaced by a smaller constant.

%---------------------
%\iffalse
\prob{prob:countcycles-Kn}\
%Let $\CC_0(\Ga)\subset \CC(\Ga)$ be the set of all cycles in $\Ga$. % without repeated nodes. 
%
Prove that
$$|\CC_0(\Knoloops{n})|=n!\sum_{k=2}^n \frac{1}{k(n-k)!}
=(e+o(1))(n-1)!
$$ 
and
$|\CC_0(\Kwloops{n})|=|\CC_0(\Knoloops{n})|+n$. 

%\prob{setcircuits-Ham} 
%Let $\CC_0(\Ga)\subset \CC(\Ga)$ be the set of all 
%simple cycles in $\Ga$. % without repeated nodes. 
%Let $\hat K_n$ be the complete (directed) graph on $n$ nodes with loops; thus  and 
%Then
%$$|\CC_0(\hat K_n)|=n!\sum_{k=2}^n \frac{1}{k(n-k)!}
%$$ 

\trivsubprob{prob:ineq-circuits-nodes}\
Corollary: if $\Ga$ is a simple graph with $|\VV(\Ga)|=n$, then $|\CC_0(\Ga)|< (e+o(1))(n-1)!$ and the bound is asymptotically sharp as $n\to\infty$.

\prob{prob:countEul-circuits-Kn}\
Let $k_n$ be the number of Eulerian circuits in $\Knoloops{n}$. %(that is, of length $n(n-1)$ circuits). 
Prove that
$
k_n\geq \prod_{d=1}^{n-1}d^d 
$
and
$
 \log k_n\geq (\frac{1}{2}-o(1))n^2\log n.
$
%[Or replace this problem with anything that gives a lower bound for $|\CC(\Knoloops{n})|$.]

\hardprob{prob:countcircuits-Kn}
Find an asymptotic formula for $|\CC(\Knoloops{n})|$.
%\fi

%Let $\Ga$ be the graph with one node and $m$ loops.
%Then $|\CC(\Ga)|=$, 
%---------------------

%\bigskip

\separline

Let $\LL$ be some set of circuits in a graph $\Ga$. The \emph{arc-circuit incidence matrix}\ $A$ is defined in an obvious way:
for $L\in\LL$ 
we set $A_{L,a}=1$ if $a\in L$ and $A_{L,a}=0$ otherwise.
% and $a\notin L$ we set $A_{L,a}=0$; otherwise $A_{L,a}$ is the multiplicity of the arc $a$ in the cycle $L$; if $L\ni a$ is pre-Eulerian, then $A_{L,a}=1$. 
Circuits $L_1,\dots,L_r$ are {\it independent}\ if the corresponding rows of 
the matrix $A$ are linearly independent. 

A system $\LL$ of independent circuits is \emph{saturated}\ if $\cup_{L\in \LL} L=\AA(\Ga)$.

A maximal set of independent circuits %without repeated arcs
is called a {\it basis of circuits}\ in $\Ga$. 
%It is possible that 
If the set of circuits $\CC(\Ga)$ is empty, 
then
%in which case 
the basis of circuits is empty.
%; the graph $\Ga$ is then called \emph{acyclic}, and $\emptyset$ is the basis of cycles. 
If the circuits concerned are cycles we get the notions of an \emph{arc-cycle incidence matrix}
and of \emph{independent cycles}. A basis of circuits all of which are cycles is called a \emph{basis of cycles}.%
%
%\footnote{There is a subtle apparent ambiguity in the notion of basis of cycles: whether it is a basis of circuits with all circuits being cycles
%or it is a maximal set of indepent cycles (so that no cycle can be added). In fact it doesn't matter, see problem \bref{basis_of_cycles_well_defined}.}
 
%\prob{basis_of_cycles_well_defined} (Bases of cycles are defined unambiguously). 
%Suppose $\LL$ is a basis of circuits that are cycles. Then 

%\trivsubprob{prob:1-cycle-quasibasis-general} 
%Every strongly connected graph possesses a one-element quasi-basis (a cycle possibly with repeated arcs).

\smallskip
In optimization problems that are soon to be introduced, simplifications and reductions will be used throughout. 
From the above ``warm-up'' exercises it should become clear that cycles are more manageable than circuits. Passing from circuit to cycle formulations and decomposing a problem into subproblems corresponding to strong components of the 
defining graph will be a common technique.  
The exercises below are to explore some details and to provide assurance to manipulations of such kind.

\trivprob{prob:1-cycle-quasibasis} 
A circuit $L$ in $\Ga$ is Eulerian if and only if 
%(a circuit containing every arc exactly once), 
$\{L\}$ is a one-element saturated system in $\Ga$.

\prob{prob:bases-circuits-to-cycles} (Reduction of bases of circuits to bases of cycles).
If $L$ is a circuit in the graph $\Ga$, let $[L]$ denote the set of all cycles that are a part of $L$. 
Let $\LL$ be a basis of circuits. Then the set of participating cycles, $\cup_{L\in\LL}[L]$, contains some basis of cycles as a subset.

\bigskip
The graph is {\it connected}\ (resp.\ {\it strongly connected}) if for any pair of nodes $v\neq v'$ there exists a path from $v$ to $v'$ \underline{or} 
(resp., \underline{and}) a path from $v'$ to $v$.
%; in this case the node $v'$ is said to be {\it reachable}\ from $v$.
% if either $v=v'$ or there exists a path from $v$ to $v'$. 

{\it Strong components}\ of $\Ga$ are the maximal strongly connnected subgraphs of $\Ga$. If they are $\Ga_1,\dots,\Ga_r$,
then $\VV(\Ga)$ is the disjoint union of $\VV(\Ga_k)$, $k=1,\dots,r$.
The set of strong components is partially ordered by the relation of reachability (existence of a path from the ``smaller'' to the ``bigger''). 
Maximal elements in that partially ordered set (poset) are called the {\it final components}. The (disjoint) union of the final components is  
denoted $\lim\Ga$. 

\trivprob{prob:relevant_arcs}\
Given a graph $\Ga$, let us call an arc $a\in\AA(\Ga)$ \emph{relevant}\ if $a\in\AA(\Ga_k)$ for some strong component $\Ga_k$ of $\Ga$.
If $\gamma$ is any cycle in $\Ga$ and $a\in \gamma$, then $a$ is relevant. Conversely, if $a$ is a relevant arc, then there exists a cycle $\gamma\ni a$. 

\trivsubprob{prob:onlystrongmatter}\
Corollary: any cycle in $\Ga$ is a cycle in some strong component of $\Ga$. 

\trivsubprob{prob:basis-reduction-to-strong}
Corollary: any basis of cycles in $\Ga$ is a (disjoint) union of bases of cycles in all strong components of $\Ga$. 

%\smallskip
(In all these statements the word {\it cycle}\ can be replaced by {\it circuit}.)

\bigskip
Put $\VV^+(v)=\{a\in\AA:\,\alpha(a)=v\}$ and $\VV^-(v)=\{a\in\AA:\,\beta(a)=v\}$. 
The \emph{outdegree}\ of the node $v$ is $d^+(v)=|\VV^+(v)|$ and the \emph{indegree}\ of $v$ is $d^-(v)=|\VV^-(v)|$. 
A graph $\Ga$ is a \emph{functional graph}\ if $d^+(v)=1$ for each $v\in\VV$. The set of arcs in a functional graph is determined by the 
\emph{adjacency function}\
$\phi:\VV\to\VV$: for each $v\in\VV$, $\VV^+(v)=\{(v\to\phi(v))\}$. 

\trivprob{prob:relevant-fungraph}\
All arcs in a functional subgraph are relevant if and only if the adjacency function is a bijection (a permutation).

\trivsubprob{prob:lim-fungraph}\
The set of relevant arcs in a functional graph $\Ga$ is $\AA(\lim\Ga)$.

\trivcsubprob{prob:basis-fungraph}\
In a functional graph $\Ga$, the set of all cycles $\CC_0(\Ga)$ is the basis of cycles.
Let $\phi_0$ be the restriction of the adjacency function on $\VV_0=\VV(\lim\Ga)$.
%Then $\phi_0$ is a permutation of $\VV_0$ and 
The cycles in $\Ga$ are trajectories of the cycles of the permutation $\phi_0$.  

\prob{prob:cycle_basis_linalg}
Given a strongly connected graph $\Ga$, consider the vector space of \emph{$0$-cochains}\ $\RR^\VV$  (vectors with components indexed by the nodes), the vector space of \emph{$1$-cochains}\ $\RR^\AA$ (vectors with components indexed by the arcs),
and the operator $d:\;\RR^\VV\to \RR^\AA$, $(d\yy)_a=y_{\beta(a)}-y_{\alpha(a)}$.  The space   $\RR^\AA$ is identified with its dual by means of the
Euclidean inner product $\langle\xx,\xx'\rangle=\sum_a x_a x'_a$. Let $U$, the space of \emph{algebraic 1-cycles}, be the cokernel of $d$, that is:
$U=\{\xx\in\RR^\AA:\; \langle\xx,d\yy\rangle=0 \,\forall\yy\in\RR^\VV\}$. Suppose $\LL$ is a basis of circuits in $\Ga$. 
Then the columns of the incidence matrix for $\LL$ form a basis in $U$ (in the usual linear algebra sense).

\subprob{prob:basis_size_inv} 
Corollary: If $\LL$ and $\LL'$ are two bases of circuits in any graph $\Ga$, then $|\LL|=|\LL'|$.

\subprob{prob:Euler_formula} 
(Euler's formula.)
If 
%$\Ga$ is a strongly connected graph and $\LL$ is a basis of cycles in $\Ga$, then 
%$|\VV|-|\AA|+|\LL|=1$. %(In particular, the cardinality of a basis is an invariant of the graph).
%
%More generally, if 
$\Ga$ is a disjoint union of $r$ strongly connected graphs and $\LL$ is a basis of circuits in $\Ga$, then  $|\VV|-|\AA|+|\LL|=r$.

\subprob{prob:basis-of-cycles-exercise}
Find a basis of cycles in the graph $\Ga$ with $\VV=\{1,2,3,4\}$ and $\AA=\{(1\to2),(2\to 3),(3\to 4),(4\to 1),(1\to 3),(2\to 4)\}$.

\subprob{prob:cycle_basis-counterex}
Give a counterexample to the effect that the assertion \bref{prob:cycle_basis_linalg}\ may fail if $\Ga$ is not strongly connected.

\subsection{Arc sums, cyclic constraints}
\label{sec:arcsums}

%\bigskip
We set up a particular case of the minimization problem \eqref{inf1} with reference to a graph. 
The input data consist of (i) a given arc-weighted graph $\Ga$:
to each arc $a\in\AA(\Ga)$ there corresponds a given $w_a>0$, (ii) some (possibly empty) set $\LL$ of independent circuits in $\Ga$, and (iii) prescribed values $t_L>0$
of cyclic products, one for each $L\in\LL$. 

\begin{NoHyper}
Let $\xx$, $\ww$ and $\xtt$ be the vectors with components $x_a$, $w_a$ and $t_L$ respectively.
The \emph{cyclic product}\ corresponding to a circuit $L$ is the function %of the vector $\xx$ of arc values%
\footnote{The definition makes sense for any set $L$ 
of arcs, but will only be used when $L$ is a circuit.}
$$
p_L(\xx)=\prod_{a\in L} x_a.
$$
%$$p_\gamma(\xx)=\prod_{a\in\gamma} x_a^{\rho(\gamma,a)},$$
Introduce the minimization problem: to determine
\end{NoHyper}%
\beq{mingraph}
f_{\Ga,\LL}(\ww,\xtt)=\min_{\xx}\langle \xx,\ww\rangle, \;\; \mbox{\rm where $\xx\geq 0$ %satisfies 
and
$p_{L}(\xx)=t_L$ for all $L\in\LL$}.
\eeq
It is a particular case of problem \eqref{inf1} and, obviously, $f_{\Ga,\LL}(\ww,\xtt)=f(\ww,\xtt,A)$, where $A$ is the arc-circuit incidence matrix
for the pair $(\Ga,\LL)$. 
(The case $\LL=\emptyset$ corresponds to the unconstrained minimum; then $f_{\Ga,\emptyset}(\ww,-)=0$.) 

The undetermined component $x_a$ of a candidate vector $\xx$ will be called the \emph{value}\ of the arc $a$, to distinguish it from the known \emph{weight}\ $w_a$. In formula \eqref{mingraph} we wrote `$\min$' instead of `$\inf$' since  the system of constraints is compact, so there is the unique minimizer by \bref{prob:existminimizer}, \bref{prob:uniqsol}.   

%If the feasibility set $\xx\geq 0:\,A\,\log\xx=\log\xtt$ is empty, we formally set $f_{\Ga,\LL}(\ww,\xtt)=+\infty$.

If $\LL$ is a basis of circuits and $t_L=1$ for any $L\in\LL$, the constraints can be stated 
without reference to a particular basis: %(cf.~\bref{prob:cycle_basis_linalg}): 
{\it for any circuit $L$ in $\Ga$, the cyclic product $p_L(\xx)=1$.} 
(It is of course possible to use the word `cycle' instead of `circuit', cf.\ \bref{prob:bases-circuits-to-cycles}.)
In this special case, the problem \eqref{mingraph} will be called \emph{homogeneous} and a shorter notation
\beq{homomin}
 f_\Ga(\ww)=f_{\Ga,\LL}(\ww,\xone)
\eeq
will be used.  If, in addition, $\ww=\xone$, we use an yet shorter notation
\beq{homo1min}
f_\Ga=f_{\Ga}(\xone)=\min_{\xx\geq 0}\sum_{a\in\AA(\Ga)} x_a, \;\; \mbox{\rm where $p_\gamma(\xx)=1$ for any 
$\gamma\in \CC_0(\Ga)$.}
%cycle $\gamma$ in $\Ga$.}
\eeq

\trivprob{prob:zerominimum} 
%Show that 
$f_{\Ga,\LL}(\ww,\xtt)=0$ if and only if $\LL=\emptyset$.

\iffalse
\trivprob{prob:min_disjoint} 
Suppose $(\Ga',\LL',\ww',\xtt')$ is a data set for the minimization problem \eqref{mingraph} satisfying all the required conditions,
and that $(\Ga'',\LL'',\ww'',\xtt'')$ is a data set for another such minimization problem. Let $\Ga$ be the disjoint union of the graphs $\Ga'$ and $\Ga''$;
denote $\LL=\LL'\cup\LL''$, $\ww=\ww'\oplus \ww''$ (direct sum, or concatenation, of vectors) and $\xtt=\xtt'\oplus\xtt''$.
Then
$
%\beq{min-sumgraph}
f_{\Ga,\LL}(\ww,\xtt)=f_{\Ga',\LL'}(\ww',\xtt')+f_{\Ga'',\LL''}(\ww'',\xtt'').
$
%\eeq
\fi

\trivprob{prob:delete-unused-arcs}\ 
(Deletion of unused arcs.)
Given a graph $\Ga$ and an independent set of circuits $\LL\neq\emptyset$ in $\Ga$, 
define the graph $\Ga'$ as follows: %$\VV(\Ga')=\VV(\Ga)$ and 
$\AA(\Ga')=\cup_{L\in\LL} L$ and $\VV(\Ga')=\cup_{L\in\LL}\alpha(L)\cup\beta(L)$. 
In other words, $\Ga'$ is obtained from $\Ga$
by deletion of 
%arcs that are not relevant (cf.\ \bref{prob:relevant_arcs}).
those arcs that do not belong to any circuit in $\LL$, so that the system $\LL$ is saturated in $\Ga'$.
Then
$f_{\Ga,\LL}(\ww,\xtt)=f_{\Ga',\LL}(\ww',\xtt)$, where $\ww'$ is the vector with components $w_a$, $a\in\AA(\Ga')$.
If $\xx$ is the minimizer for the problem \eqref{mingraph}, then $\AA(\Ga')=\{a:\, x_a\neq 0\}$.

\trivsubprob{prob:reduction-to-strong}
(Reduction to strongly connected graphs).
Let $\Ga_1,\dots,\Ga_r$ be the strong components of the graph $\Ga$ and $\LL=\cup_{k=1}^r \LL_k$ be the coresponding decomposition of 
a given set of independent circuits in $\Ga$ into the sets of circuits in the strong components,
cf.\ 
\brefA{prob:relevant_arcs}.
%\bref{prob:onlystrongmatter}. 
Also, let $\ww_k$ and $\xtt_k$ be the subvectors of $\ww$ and $\xtt$ pertaining
to the $k$-th strong component. Then 
$f_{\Ga,\LL}(\ww,\xtt)=\sum_{k=1}^r f_{\Ga_k,\LL_k}(\ww_k,\xtt_k)$. In particular, $f_\Ga=\sum_{k=1}^r f_{\Ga_k}$.

\trivsubprob{prob:min-for-fungraph}
If $\Ga$ is a functional graph, then $f_\Ga=|\VV(\lim\Ga)|=|\AA(\lim\Ga)|$.
The minimum possible value $f_\Ga=1$ occurs when the adjacency function $\phi$ is the succession function of a linear order on $\VV$ with loop 
at the end. The maximum value $f_\Ga=|\VV(\Ga)|$ occurs when $\phi$ is a bijection. 

\prob{prob:mean-min-fungraph}
Let $|\VV|=n$ and suppose all maps $\VV\to\VV$ are assigned equal probability $n^{-n}$.
Then the expectation $\prE f_\Ga$ %of $f_\Ga$ 
for a random functional graph with node set $\VV$
as a function of $n$ asymptotically behaves as $\sqrt{\pi n/2}$.

\trivprob{prob:isomorph}
(Isomorpisms and anti-isomorphisms).
A graph \emph{isomorphism}\ $\phi:\,\Ga\to\Ga'$ is  a pair $\phi=(\phi_{\VV},\phi_{\AA})$ of bijective maps $\phi_{\VV}:\,\VV(\Ga)\to\VV(\Ga')$ and  $\phi_{\AA}:\,\AA(\Ga)\to\AA(\Ga')$ preserving the arc-node incidences:
$\alpha'(\phi_{\AA}(a))=\phi_{\VV}(\alpha(a))$, $\beta'(\phi_{\AA}(a))=\phi_{\VV}(\beta(a))$ for all $a\in\AA(\Ga)$.
Similarly, $\phi$ is \emph{anti-isomorphism}\ if it has the same properties but reverses the orientation of arcs:
$\alpha'(\phi_{\AA}(a))=\phi_{\VV}(\beta(a))$, $\beta'(\phi_{\AA}(a))=\phi_{\VV}(\alpha(a))$ for all $a\in\AA(\Ga)$.

We call the problems of the form \eqref{mingraph} for $\Ga$ and $\Ga'$ isomorpic by means of $\phi$ which is either an isomorphism or anti-isomorphism
if the data are correspondingly related: $w'_{\phi_{\AA}(a)}=w_{a}$,  the circuits in $\LL'$
are the $\phi_{\AA}$-images of the circuits in $\LL$, and $t'_{\phi_{\AA}(L)}=t_{L}$ for all $L\in\LL$. 

\smallskip
For two problems isomorphic by means of $\phi$, 
 $f_{\Ga,\LL}(\ww,\xtt)=f_{\Ga',\LL'}(\ww',\xtt')$. The minimizers $\xx$ and $\xx'$ are related by $x'_{\phi_{\AA}(a)}=x_{a}$ for all $a\in\AA(\Ga)$.     

\trivsubprob{prob:automorph}\ (Automorphisms).
If $\phi$ is an automorphism or anti-automorphism of the graph $\Ga$ such that the set of cycles $\LL$ and the data vectors $\ww$ and $\xtt$ are $\phi$-invariant, then the minimizer $\xx$ is also $\phi$-invariant.

\trivsubprob{prob:automorphic_minimizer}
Let $\xx$ be the minimizer for the homogeneous problem with unit weights on the graph $\Ga$
 (so that $f_{\Ga}=\langle\xx,\xone\rangle$). If $\phi$ is an automorphism or anti-automorphism of $\Ga$ and $\phi_{\AA}(a)=a'$, then $x_a=x_{a'}$. 

\noindent
\begin{minipage}[c]{0.65\textwidth} 
\subprob{prob:non-anti-automorphic_minimizer}\
(Cheating at anti-automorphisms).  
Consider the graph $\Ga$ with $\VV(\Ga)=\{1,2,3\}$, $\AA(\Ga)=\{a,a',b,c\}$, where $a=(1\to 2)$, $a'=(2\to 1)$, $b=(2\to 3)$,
$c=(3\to 1)$. 
Let $\xx$ be the minimizer for the homogeneous problem on $\Ga$ with unit weights. 
Show that $x_{a}\neq x_{a'}$.
Isn't there a contradiction with \bnref{prob:automorphic_minimizer}\ given that
\end{minipage}
\hspace{2em}
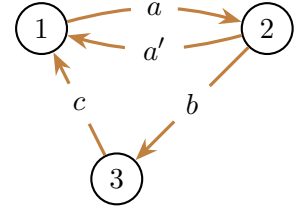
\begin{minipage}[c]{0.3\textwidth} 
%\vspace*{1ex}
\begin{tikzpicture}
\begin{scope}[every node/.style={circle,thick,draw}]
    \node (A) at (0,2) {1};
    \node (B) at (3,2) {2};
    \node (C) at (1,0) {3};
\end{scope}

\begin{scope}[>={Stealth[brown]},
              every node/.style={fill=white,circle},
              every edge/.style={draw=brown,very thick}]
    \path [->] (A) edge[bend left=15] node {$a$} (B);
    \path [->] (B) edge[bend left=15] node {$a'$} (A);
    \path [->] (B) edge node {$b$} (C);
    \path [->] (C) edge node {$c$} (A); 
\end{scope}\end{tikzpicture}
\end{minipage}
\\[0.5ex]
\hspace*{1em} 
(a) The pair (identical map on $\VV$, transposition $a\leftrightarrow a'$ on $\AA$) is an anti-automorphism of $\Ga$; 
\\[0.5ex]
\hspace*{1em} 
(b) The pair (transposition $1\leftrightarrow 2$ on $\VV$,  the induced map on $\VV$) is an anti-automorphism of $\Ga$?

\noindent
\begin{minipage}[c]{0.65\textwidth} 
\prob{prob:ex-quasibasis}
Consider the (non-simple) graph $\Ga$ defined by the data: $\VV(\Ga)=\{1,2,3\}$, $\AA(\Ga)=\{a,b,c,d,e\}$, where $\alpha(a)=\alpha(d)=\beta(c)=1$,
$\beta(a)=\beta(d)=\alpha(b)=\alpha(e)=2$, $\beta(b)=\beta(e)=\alpha(c)=3$.
Let $\gamma_1=\{d,b,c\}$, $\gamma_2=\{a,e,c\}$ and $\gamma_3=\{a,b,c\}$. 
%(Fig.~\ref{fig:ex-quasibasis}).
The system $\LL=\{\gamma_1,\gamma_2\,\gamma_3\}$ is a basis of cycles and the system $\hat \LL=\{\gamma_1,\gamma_2\}$ is a saturated system. 
\end{minipage}
\hspace{2em}
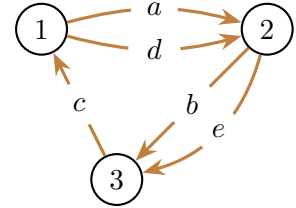
\begin{minipage}[c]{0.3\textwidth} 
%\vspace*{-1ex}
\begin{tikzpicture}
\begin{scope}[every node/.style={circle,thick,draw}]
    \node (A) at (0,2) {1};
    \node (B) at (3,2) {2};
    \node (C) at (1,0) {3};
\end{scope}

\begin{scope}[>={Stealth[brown]},
              every node/.style={fill=white,circle},
              every edge/.style={draw=brown,very thick}]
    \path [->] (A) edge[bend left=15] node {$a$} (B);
    \path [->] (A) edge[bend right=15] node {$d$} (B);
    \path [->] (B) edge node {$b$} (C);
    \path [->] (B) edge[bend left=30] node {$e$} (C);
    \path [->] (C) edge node {$c$} (A); 
\end{scope}\end{tikzpicture}
\end{minipage}

\medskip
\noindent
a) Find $\hat m(\ww)=f_{\Ga,\hat\LL}(\ww,\xone)$.
\\[1ex]
b)
Assuming that $w_a=w_b$ and $w_d=w_e$ (call such a $\ww$ a \emph{symmetric weight}), find $m(\ww,t)=f_{\Ga}(\ww,\xtt)$, where $t_1=t_2=1$, $t_3=t>0$ (a free parameter).
\\[1ex]
c) Check by comparing the explicit formulas in the case of symmetric weight that $\hat m(\ww)\leq m(\ww,t)$ for any $t>0$
(cf.\ \bref{prob:relaxcons1}), and the equality
holds only if $t=p_{\gamma_3}(\xx)$, where $\xx$ is the minimizer for problem (a).  

\subsection{Miscellaneous variants of graphic sums}
\label{sec:graphs-variants}

The adjective {\em graphic}\ in the heading means ``determined by a graph'', ``pertaining to graphs'',
similarly to its use in the %established 
term 
{\em graphic matroid}, say. 

\subsubsection{Sums of quotients and a link to Shapiro's cyclic sums}
\label{ssec:sums-of-quotients}

As before, we assume that given is a graph $\Ga$ with arc weights $w_a$, $a\in\AA$. This time, however, the optimization variables
will be node values $y_{v}$, $v\in\VV$. Denote
\beq{quotient_sum}
S_{\Ga}(\ww|\yy)=\sum_{a\in\AA} w_a \frac{y_{\alpha(a)}}{y_{\beta(a)}}.
\eeq
\emph{Admissible}\ are those vectors $\yy\geq 0$ for which $y_{\beta(a)}>0$ $\forall a\in\AA$.
In particular, if $\Ga$ is strongly connected then $\yy$ is admissible if and only $\yy>0$.
The minimization problem is to find 
\beq{minquotient}
 f^\div_{\Ga}(\ww)=\inf_{\yy} S_{\Ga}(\ww|\yy),
\eeq
where $\yy$ runs over the set of admissible vectors. 
An admissible vector $\yy$ is a \emph{minimizer}\ for problem \eqref{minquotient}\ if $S_{\Ga}(\ww|\yy)=f^\div_{\Ga}(\ww)$.
In the case of $\ww=\xone$ we write $S_{\Ga}(\yy)=S_{\Ga}(\xone|\yy)$ and $f^\div_{\Ga}=f^\div_{\Ga}(\xone)$.

\prob{prob:minquotient-mingraph}
The infimum \eqref{minquotient} coincides with minimum in the homogeneous case of problem \eqref{mingraph},
$$
 f^\div_{\Ga}(\ww)=f_{\Ga}(\ww).
$$
In particular, $ f^\div_{\Ga}=f_{\Ga}$.

\subprob{prob:minquotient-nominimizer}
A minimizer for \eqref{minquotient} does not always exist. The necessary and sufficient condition for the existence of a minimizer is
as follows. In the poset of strong components of the graph $\Ga$ there are no chains of length $\geq 2$ and if there are chains of
length 1 (i.e.\ if $\Ga$ is not strongly connected) then every non-final component is a singleton without loop.

%\subsubsection{A link to Shapiro's cyclic inequality}

\separline
Let us look at a formal, though shallow, connection to classical cyclic inequalities. 
%More interesting is Shallit's problem in Section~7

\prob{prob:Nesbitt} Interpret \emph{Nesbitt's inequality}
$$
 \frac{u_1}{u_2+u_3}+\frac{u_2}{u_3+u_1}+\frac{u_3}{u_1+u_2}\geq \frac{3}{2}
$$ 
as the inequality of the form $S_\Ga^\div(\yy)\geq f^\div_\Ga$
for an appropriate graph $\Ga$. 

\subprob{prob:Shapiro}
It is known that there holds \emph{Shapiro's cyclic inequality} 
$$
\sum_{i=1}^n \frac{u_i}{u_{i+1}+u_{i+2}}\geq C_n,
$$
where $n\geq 3$, $u_{n+1}=u_1$ and $u_{n+2}=u_2$
and $C_n$ is some positive constant (see Comment in Part II).
Show that if $n$ is odd, then Shapiro's inequality can be formally presented in the form  $S_{\Ga_n}^\div(\yy|\ww)\geq f^\div_{\Ga_n}(\ww)$
with an appropriate graph  $\Ga_n$ and some weight vector $\ww=\ww^{(n)}$; however, for $n\geq 5$, the vector $\ww$ has positive as well as negative components.
(Thus we do not really reduce Shapiro's problem with odd $n\geq 5$ to problem \eqref{minquotient}).   

\subsubsection{Exponential sum for a harmonic function}
\label{ssec:harmonic}

Suppose $\Gamma$ is an {\em undirected} 
%\footnote{This section is the only place where undirected graphs appear.} 
graph without loops,
that is, $\VV_+(v)=\VV^-(v)$ (to be denoted $\VV(v)$
in this section) and $v\notin \VV(v)$ for any $v\in\VV$. Here we use the terms {\em vertices}\ and {\em edges}\
instead of `nodes' and `arcs'.
Assume that a subset of vertices $\VV_0\subset\VV$ is marked as ``boundary''.  The complementary subset
$\VV_1=\VV\setminus\VV_0$ is the set of ``inner vertices''.

A function $h$ defined on $\VV$ is {\em harmonic
%\footnote{One could say {\em out-harmonic}\ (and define {\em in-harmonic}\ functions similarly, summing over
%$\VV^-(v)$ instead).
%Usually the Laplace operator and harmonic functions are considered for undirected graphs; that situation corresponds to a symmetric incidence matrix.}
% 
on $\VV_1$}\ if for every $v\in\VV_1$ the Laplace equation 
$$
%\sum_{v'\in\VV^+(v)}(h(v)-h(v'))=0
\sum_{v'\in\VV(v)}(h(v)-h(v'))=0
$$
is satisfied.
The values $h(v)$, $v\in\VV_0$ are called {\em boundary values}.

Consider the optimization problem: to find a function $h$ harmonic on $\VV_1$ that minimizes the functional
$$
 \phi[h]=\sum_{v\in \VV} w_v e^{h(v)}
$$
under the linear constraint on boundary values:
$$
 \sum_{v\in\VV_0} b_v h(v)=\tau.
$$
The set of data consists of the weights $\{w_v\geq 0,\; v\in\VV\}$ and some real numbers $t$ and $\{b_v,\;v\in\VV_0\}$.

\prob{prob:harmonic-graph} 
Interpret this problem as a problem of type 
\eqref{sum1}--\eqref{cons_prod}. 

In the case
$\ww=\xone$
 %$w_v=1$ for all $v\in\VV$ 
and 
$\xbb=\xone$,
denote the value of the minimum by 
%$f_\Gamma^\Delta(t)$.
$F^\Delta(\tau)$.
Prove that always
%$b_v=1$ for all $v\in\VV_0$, always
\beq{harmonic-min}
% \min \phi[h]
F^\Delta(\tau)\leq |\VV|e^{\tau/|\VV_0|}.
\eeq

\prob{prob:harmonic-graph-lb} 
%(An yet another generalization of the AGM inequality.)
 Prove that the inequality 
\eqref{harmonic-min} turns to equality 
%$$
% \min f[\phi] =|\VV| e^{t/|\VV_0|}.
%$$
in the following cases. (Which, neither separately nor in combination, are not necessary.)

%(a) For $t=0$.

(a) In the case where the bipartite graph $(\VV_0,\VV_1,\AA_{01})$ is regular.
Its set of edges consists of those edges from $\AA$ that
connect vertices from $\VV_1$ with vertices from $\VV_0$;
regularity means that there exist positive integer constants $p$ and $q$ such that
(i) for any $v\in\VV_0$, 
$|\VV(v)\cap\VV_1|=p$,
%$|\VV^-(v)\cap\VV_1|=p$,
and (ii) for any $v\in\VV_1$, 
$|\VV(v)\cap\VV_0|=q$.
%$|\VV^+(v)\cap\VV_0|=q$.
(Note that $p/|\VV_1|=q/|\VV_0|$.)

(b) For a linear undirected graph with two-point boundary.
(Precisely: $\VV=\{0,1,\dots,n\}$, $\VV_0=\{0,n\}$, the edges are $(i-1,i)$, $i=1,\dots,n$.)

\prob{prob:harmonic-example}
Consider the (undirected) graph with the set of vertices
%shown on picture:
$\VV=\{1,2,3,4\}$ and the edges $(12)$, $(13)$, $(14)$,
$(23)$. Let the boundary be $\VV_0=\{3,4\}$. 
Compute $F^\Delta(\tau)$ (for one or more values of $t$, or plot a graph). Compare the found values with the upper bound $4e^{\tau/2}$ from \eqref{harmonic-min}.

\prob{prob:harmonic-const-ratio}
Prove that for any graph and any boundary set
$\VV_0\subset\VV$, the ratio
$F^\Delta(\tau)/(|\VV|e^{\tau/|\VV_0|})$ is independent of $\tau$.

%%%%%%%%%%%%%%%%%%%%%%%%%%%%%%%%%%%%%%%%%
\subsection{Trees: node sums, path products}
\label{sec:trees}

%\bigskip
\begin{NoHyper}
A \emph{rooted tree}\ is a %(directed, as everywhere here) 
graph $T$ with a distinguihed node $\rho$, the \emph{root}, such that for any node $v\neq\rho$
there exists a unique path from $\rho$ to $v$. 
The natural partial order on $\VV(T)$ is: $v\le v'$ if and only if  there exists the (necessarily unique) path from $v$ to $v'$. 
Put $[v,v']=\{v''\in\VV(T)\,|\, v\leq v''\leq v'\}$.
The strong components of $T$ are singletons and can be identified with nodes of $T$. 
The nodes $v\in\lim T$
 are the \emph{leaves}.\footnote{A pedant would %prefer to 
write $\{v\}\in\lim T$.}  
We denote the set of leaves $\Lambda(T)$ (or $\Lambda$, if the tree is fixed in the current context). 
Any node $v\neq\rho$ has indegree one and it has the \emph{parent}\ node denoted $\alpha(\rho)$ such that $e(v)=(\alpha(v)\to v)$
is the unique arc ending in $v$.%  
\footnote{Incidentally we have the relation $\alpha(e(v))=\alpha(v)$, where the symbol $\alpha$ on the left denotes the beginning-of-the-arc function, 
and on the right it denotes the parent function.}
The range of the function %$\alpha:\VV(T)\setminus\{\rho\}\to\VV(T)$ is $\VV(T)\setminus\Lambda(T)$.  
$\alpha:\VV\setminus\{\rho\}\to\VV$ is $\VV\setminus\Lambda$.  

\end{NoHyper}

A \emph{rooted forest}\ is a disjoint union of rooted trees; it inherits the partial order and the parent function from its components. 
The minimal elements are the roots of the constituting trees and the maximal elements are their leaves.
Let $F$ be a forest and $\Lambda=\Lambda(F)$ be the set of its leaves. To every leaf $\lambda\in\Lambda$ there corresponds a unique maximal path
in $F$, equivalently, a maximal chain in $F$ considered as a poset; we denote it $\pi(\lambda)$. If $F$ is just a tree, then $\pi(\lambda)=[\rho,\lambda]$.
The \emph{height}\ of a tree or of a forest is $h=\max_{\lambda\in\Lambda}|\pi(\lambda)|-1$.
 
If $\yy$ is a vector with components indexed by the nodes of $F$ and $[v',v'']$ is any interval of nodes, we  
denote
$$
P_{[v',v'']}(\yy)=\prod_{v'\leq v\leq v''} y_v
$$
the product of node values over the path from $v'$ to $v''$. Particularly important in the sequel are maximal chain products $P_{\pi(\lambda)}$,
$\lambda\in\Lambda$.  

%An {\it arc-weighted}\ graph is a pair $(\Ga,\xx)$, where $\xx$ is a vector with components $x_e$ labeled by the arcs of $\Ga$.
%Similarly, a {\it node-weighted}\ graph is a pair $(\Ga,\yy)$, where $\yy$ is a vector with components $y_v$ labeled by $v\in\VV(\Ga)$.
%We will also use vectors with components labeled by other sets, such as $\xtt=(t_j)$ with index $j$ enumerating elements $L_j$ of some basis of circuits $\BB$.      
%A vector $\xone$ has all components equal to $1$; the dimension should be clear from context.

%\smallskip
%The data for the next minimization problem consist of a node-weighted forest $F$ with weight vector $\ww>0$ (its components are $w_v$, $a\in\VV(F)$),
%and the vector $\xtt>0$ of prescribed maximal path products (with components $t_\lambda$, $\lambda\in\Lambda(F)$).  
%The problem is to determine
%\beq{mintreeweighted}
%m_{F}(\ww,\xtt)=\inf_{\yy}\langle \yy,\ww\rangle, \;\; \mbox{\rm where $\yy> 0$ satisfies $p_{\gamma(\lambda)}(\yy)=t_\lambda$ for all $\lambda\in\Lambda$}.
%\eeq

\smallskip
The data for the next minimization problem consist of a forest $F$
and the vector $\xtt>0$ of prescribed maximal path products (with components $t_\lambda$, $\lambda\in\Lambda(F)$).  
The problem is to determine
\beq{mintree}
m_{F}(\xtt)=\inf_{\yy}\sum_{v\in\VV} y_v, \;\; \mbox{\rm where $\yy> 0$ satisfies $P_{\pi(\lambda)}(\yy)=t_\lambda$ for all $\lambda\in\Lambda$}.
\eeq
As before, we use the abbreviation $m_F=m_F(\xone)$.
Obviously, problem \eqref{mintree} is a particular case of problem \eqref{inf1}% 
%\footnote{It is less obvious that \eqref{mintree} can be seen as a particular case of problem \eqref{mingraph}, see \bref{prob:tree-to-general} below.}
\ and the system of constraints is compact; therefore, `$\inf$' in \eqref{mintree} can be replaced by `$\min$'
and there exists the unique minimizer. 

% WRONG! (See Table~1)
%\prob{prob:tree-unitproduct}\
%(Unit products, unrooted trees).
%If $T$ and $T'$ are two rooted trees (not forests) that differ only by the position of the root, 
%then $m_{T}(\xone)=m_{T'}(\xone)$. 
%Consequently, the functional $T\mapsto m_T$ is well-defined on the set of {\em topological}\ (unrooted) trees.

\smallskip
We denote by $T_v$ the subtree of a tree $T$ rooted at the node $v$. Its set of nodes is $\VV(T_v)=\{v'\in\VV(T)\,|\, v'\geq v\}$.
For the given $T$, we write $\Lambda_v=\Lambda(T_v)$.

\trivprob{prob:forest-to-tree} 
Let $F$ be a forest that consists of trees $T_i$, $1\leq i\leq r$ and let $\xtt=\oplus\xtt_i$ be the decomposition of the vector $\xtt$ into subvectors
corresponding to the individual trees. Then
$$
m_F(\xtt)=\sum_{i=1}^r m_{T_i}(\xtt_i).
$$  

\prob{prob:treemin-ineq} (Compare with \bref{prob:consineq}.)
The constraints in problem \eqref{mintree} can be stated in the form of inequalities without changing the result, that is,
\beq{mintree-ineq}
m_{F}(\xtt)=\inf_{\yy}\sum_{v\in\VV} y_v, \;\; \mbox{\rm where $\yy> 0$ satisfies $P_{\pi(\lambda)}(\yy)\geq t_\lambda$ for all $\lambda\in\Lambda$}.
\eeq

\trivsubprob{prob:treemin-monotone-in-t}
Corollary: if $\xtt'\geq \xtt$, $\xtt'\neq\xtt$, then $m_F(\xtt')> m_F(\xtt)$. 
(Compare with \bref{prob:notmonotone-wrt-products}.)
%, \bref{prob:ex-quasibasis}.)

\prob{prob:tree-quotient}
(A ``tree quotient sum'', compare with \eqref{quotient_sum}).
Given a tree $T$ with root $\rho$, 
%and a vector $\xtt$ as above, 
define
\beq{treesum-quotient}
Q_T(\xtt|\yy)=\sum_{v\in\VV(T)\setminus\{\rho\}} \frac{y_{\parent{v}}}{y_v}+\sum_{\lambda\in\Lambda(T)} t_{\lambda} \frac{y_\lambda}{y_{\rho}}
\eeq
and let 
$$
 m_T^{\div}(\xtt)=\inf_{\yy>0} Q_T(\xtt|\yy).
$$
Exhibit a forest $F$ with $|\alpha^{-1}(\rho)|$ component trees and leaves of $F$ being in one-to-one correspondence with leaves of $T$, 
$\Lambda(F)\ni\hat\lambda\leftrightarrow\lambda\in\Lambda(T)$, such that
$$
 m_T^{\div}(\xtt)=m_F(\hat\xtt)
$$
where $\hat t_{\hat\lambda}=t_\lambda$.

\prob{prob:tree-to-general}\ 
(Reduction of the tree problem to a minimization problem under cyclic constraints).
%Reduce problem for $W_T$ to problem for $W_\Ga$.
%
Interpret problem \eqref{mintree} as a particular case of problem \eqref{mingraph}. Specifically,
let $T$ be a rooted tree and $\xtt=(t_\lambda)_{\lambda\in\Lambda}$. Exhibit a graph $\Ga$ with $|\AA(\Ga)|=|\VV(T)|$ and 
a basis $\LL$ of cycles in $\Ga$ being in one-to-one correspondence with leaves of $T$, 
$\Lambda(T)\ni\lambda\leftrightarrow\gamma({\lambda})\in\LL$, such that
$$
m_T(\xtt)=f_\Ga(\xone,\xtt).
$$
(The vector $\xtt$ in the right-hand side determines the values of cycle products in $\Ga$: $p_{\gamma({\lambda})}=t_\lambda$.)

%\subprob{??}
%Reduce problem for $W_\Ga$ to problem for $W_T$.

\prob{prob:tree-minimizer}\
(Critical point equations).  
The minimizer for problem \eqref{mintree} 
is found as the solution of the system of equations
\beq{tree-minimizer}
\begin{array}{l}
\displaystyle
y_v=\sum_{\lambda\in\Lambda_v} y_\lambda,\quad v\in\VV\setminus\Lambda,
\\[2ex]
\displaystyle
P_{\pi(\lambda)}=t_\lambda, \quad \lambda\in\Lambda.
\end{array}
\eeq

\trivsubprob{prob:tree-quotient-minimizer}\ 
The minimizer for the problem of minimization of the quotient sum \eqref{treesum-quotient} satisfying the normalization condition 
$y_\rho=1$ is unique and is found as  the solution of the system of equations identical to \eqref{tree-minimizer} except
that there is no equation with $v=\rho$ in the first group of equations.

\prob{prob:tree-recurrrence}
(Recurrence).
Suppose $T$ is a tree with more than one node.
Let $\xtt_v$ be the subvector of the vector of products $\xtt$ with indices $\lambda\in\Lambda_v$.
The following recurrence relation for the minimum holds:
\beq{rec-trees}
 m_T(\xtt)=\min_{k>0}\left(\frac{1}{k}+\sum_{v\in\alpha^{-1}(\rho)} m_{T_v}(k\xtt_v)\right).
\eeq
Suppose the minimum in \eqref{rec-trees} is attained at $k=k^*$. Let $\yy^{(T)}(\xtt)$ be the minimizer in problem \eqref{mintree} for the tree $T$,
and let $\yy^{(T')}(\xtt)$ be  the minimizer in problem \eqref{mintree} for the forest $\cup_{v\in\alpha^{-1}(\rho)} T_v$ obtained from $T$
by removal of the root.
Then
\beq{rec-trees-minimizer}
y^{(T)}_\rho(\xtt)=\frac{1}{k^*}
\;\quad
\mbox{\rm and}\;\quad
y^{(T)}_v(\xtt)=y^{(T')}_v(k^*\xtt), \;\; v\neq\rho. 
\eeq     

\prob{prob:const-height-tree-radicals}
(Solvability in radicals for constant height trees).
A tree $T$ in which all maximal chains have the same length $\ell$ is called a \emph{constant height tree}; more precisely,
$T$ is a \emph {tree of constant height $h=\ell-1$}. Suppose $T$ is such a tree.  

(a) The function $\xtt\mapsto  m_T(\xtt)$ is homogeneous of order $1/\ell$.

(b) %Show that 
$m_T(\xtt)$ can be expressed in radicals involving rational numbers and the indeterminates $t_\lambda$.

(c) The system of algebraic equations \eqref{tree-minimizer} for the components of
the minimizer is solvable in radicals. 

%\prob{prob:isomorphic-trees}  [WRONG!] 
%(Isomorphic non-rooted trees).
%Let $T_1$ and $T_2$ be two rooted trees that satisfy the conditions:
%\\[0.5ex]
%\hspace*{1em}
%(i) they are isomorphic as non-rooted undirected trees; 
%\\[0.5ex] 
%\hspace*{1em}
%(ii) The roots of $T_1$ and $T_2$ have outdegree $>1$.
%\\[0.5ex]
%Then $m_{T_1}=m_{T_2}$. 

\prob{prob:dilatation-monotone-trees} 
(Motonicity of minimizer in a weak sense).
Let us call a map $M:\,\RR_{>0}^m\to\RR_{>0}^n$ \emph{dilatation-monotone}\ if $M(r\xuu)\geq M(\xuu)$ for any $r\geq 1$.
Given a tree $T$, let $M_T:\xtt\mapsto\yy^{(T)}$ be the function mapping the vector of maximal chain products to the minimizer
for problem \eqref{mintree}. Then $M_T$ is dilatation-monotone. 

\subprob{prob:no-strong-monotone-trees}
(Minimizer fails to be strongly monotone even for constant height trees).
Give an example of a constant height tree $T$ for which the map $M_T$ is not monotone in the strong sense: $\xtt>\xtt'$ 
does not imply $M_T(\xtt)\geq M_T(\xtt')$. (Recall that $\xx>\xx'$ means $x_i>x'_i$ for all components of the vector $x$.)

\subprob{prob:no-dilatation-monotone-graph}
(Weak monotonicity of minimizer fails in problem \eqref{mingraph}). 
Refute the following assertion by a counterexample: {\it ``for any strongly connected graph $\Ga$ and a basis of cycles $\LL$ in $\Ga$, the function 
$\xtt\mapsto \xx^{*}$ whose value is the minimizer for problem \eqref{mingraph} with $\ww=\xone$, is dilatation-monotone.''}

\prob{prob:tree-linearization} (Properties of the linearized critical point system).
Let $\yy^*$ be the minimizer for the problem \eqref{mintree} and $\hat\yy^*$ be the subvector of $\yy$ with component indices $\lambda\in\Lambda$.
According to the critical point equations \eqref{tree-minimizer}, $\hat\yy$ is a fixed point of the map
$$
 M:\,\hat\yy'\mapsto \hat\yy, \qquad \log\hat y_\lambda=
 \log t_\lambda-
 \log\sum_{v\in[\rho,\alpha^{-1}(\lambda)]}\sum_{\lambda'\in\Lambda_v} \hat y_{\lambda'}.
$$ 
For $\lambda\in\Lambda$ and $v\in\VV\setminus\Lambda$ put $Q_{\lambda,v}=1$ if $v\in[\rho,\lambda]$ and $0$ otherwise.
Show that the linearization of the map $M$ at $\hat\yy^*$ has matrix $B$ with entries 
$$
B_{\lambda,\lambda'}=-\sum_{v\in\VV\setminus\Lambda} Q_{\lambda,v} \frac{y^*_\lambda}{y^*_v} Q_{\lambda',v} . 
$$  

\subprob{prob:tree-linearization-spectrum}
Show that the spectrum of the matrix $B$ lies in $[-h,0)$, where $h$ is the height of the tree.

\subprob{prob:tree-numerical-iteration}
Show that the iterations $\hat \yy^{(n+1)}=\tau M\hat\yy^{(n)}+(1-\tau)\hat\yy^{(n)}$ converge to $\hat\yy^*$ provided
$0<\tau<\frac{2}{h+1}$ (where $h$ is the height of the tree) and the initial vector $\hat\yy^{(0)}$ is sufficiently close to $\hat\yy^*$.   

\prob{prob:trees-analytical-numerical}
Find the minima $m_T$ and the minimizers for some small trees
algebraically (in terms of roots of irreducible polynomials) and/or numerically.

\medskip
\noindent
\begin{minipage}{0.78\textwidth}
\prob{prob:snowflake}
Find $m_T$ for the ``snowflake'' tree for three non-isomorphic positions of the root (at the center, in the middle ring and on the perimeter).
\end{minipage}
\hspace{1.6em}
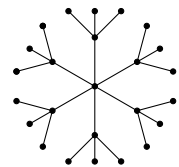
\begin{minipage}{0.2\textwidth}
\begin{tikzpicture}[x=0.5cm,y=0.5cm]
\draw [line width=0.4pt] (2,0)--(2,0.7)--(2,2)--(2,3.3)--(2,4);
\draw [line width=0.4pt] (0.26795,3)--(1.133975,2.5)--(2,2)--(2.866,1.5)--(3.732,1);
\draw [line width=0.4pt] (0.26795,1)--(1.133975,1.5)--(2,2)--(2.866,2.5)--(3.732,3);
\draw [line width=0.4pt] (1.3,4)--(2,3.3)--(2.7,4);
\draw [line width=0.4pt] (1.3,0)--(2,0.7)--(2.7,0);
\draw [line width=0.4pt] (-0.082, 2.3938)--(0.87417,2.65)--(0.61795, 3.60622);
\draw [line width=0.4pt] (4.082, 1.60622)--(3.12583,1.35)--(3.382, 0.3938);
\draw [line width=0.4pt] (-0.082, 1.60622)--(0.87417,1.35)--(0.61795, 0.3938);
\draw [line width=0.4pt] (3.382, 3.60622)--(3.12583,2.65)--(4.082, 2.3938);
\draw [fill=black] (2,0) circle (1pt);
\draw [fill=black] (2,0.7) circle (1pt);
\draw [fill=black] (2,2) circle (1pt);
\draw [fill=black] (2,3.3) circle (1pt);
\draw [fill=black] (2,4) circle (1pt);
\draw [fill=black] (0.26795,3) circle (1pt);
\draw [fill=black] (3.12583,1.35) circle (1pt);
\draw [fill=black] (0.87417,2.65) circle (1pt);
\draw [fill=black] (3.732,1) circle (1pt);
\draw [fill=black] (0.26795,1) circle (1pt);
\draw [fill=black] (3.12583,2.65) circle (1pt);
\draw [fill=black] (0.87417,1.35) circle (1pt);
\draw [fill=black] (3.732,3) circle (1pt);
\draw [fill=black] (1.3,0) circle (1pt);
\draw [fill=black] (2.7,0) circle (1pt);
\draw [fill=black] (1.3,4) circle (1pt);
\draw [fill=black] (2.7,4) circle (1pt);
\draw [fill=black] (-0.082, 2.3938) circle (1pt);
\draw [fill=black] (0.61795, 3.60622) circle (1pt);
\draw [fill=black] (-0.082, 2.3938) circle (1pt);
\draw [fill=black] (0.61795, 3.60622) circle (1pt);
\draw [fill=black] (3.382, 0.3938) circle (1pt);
\draw [fill=black] (4.082, 1.60622) circle (1pt);
\draw [fill=black] (-0.082, 1.60622) circle (1pt);
\draw [fill=black] (0.61795, 0.3938) circle (1pt);
\draw [fill=black] (3.382, 3.60622) circle (1pt);
\draw [fill=black] (4.082, 2.3938) circle (1pt);
\end{tikzpicture}
\end{minipage}

%%%%%%%%%%%%%%%%%%%%%%%%%%%%%%%%%%%%%%%%%%%%%%%%%%%%%%%%%%%%%%%%%%%%%%
\subsection{Extremal problems} % involving the minimum}
\label{sec:extremal_problems}

\subsubsection{Extremal problem for node-weighted rooted trees}
\label{ssec:trees-extremal}

This section is concerned with extremal values of the functional $T\mapsto m_T$ on families of rooted trees with prescribed number of nodes or nodes and leaves.

\smallskip
Let $\mathcal{T}(n,\ell)$ denote the set of (isomorphism classes of) all rooted trees with $|\VV|=n$
nodes and $|\Lambda|=\ell$ leaves. A distinguished member of $\mathcal{T}(n,\ell)$ is the {\em palm tree} $T=\Pl(n,\ell)\in \mathcal{T}(n,\ell)$  defined as follows: 
 $\VV(T)=\{0,\dots,n-1\}$ (node $0$ being the root), $\Lambda(T)=\{n-\ell,\dots,n-1\}$, and $\AA(T)$ is the union of ``trunk'' $\{(i-1\to i), i=1,\dots,n-\ell-1\}$
and ``frond'' $\{(n-\ell-1\to n-j), j=1,\dots,\ell\}$. 
The sets $\mathcal{T}(n,n-1)$ and $\mathcal{T}(n,1)$ are singletons and their unique members are the respective palm trees.
The tree $\Pl(n,n-1)\in \mathcal{T}(n,n-1)$ has empty trunk. 
The tree $\Pl(n,1)\in \mathcal{T}(n,1)$ is a ``linear'' tree; it will be denoted $\mathrm{Lin}(n)$. 

\medskip
\noindent
\begin{minipage}{0.78\textwidth}
\hspace*{1em}
Although $\mathrm{Lin}(n)$ looks like a bare trunk, technically the last arc, $(n-2\to n-1)$ constitutes the singleton frond. The trees $\Pl(5,2)$ and $\mathrm{Lin}(4)$ are shown on the right. The trunk arcs are solid and the frond arcs are dotted.  
\end{minipage}
\hspace{1.6em}
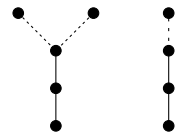
\begin{minipage}{0.2\textwidth}
\begin{tikzpicture}[x=0.5cm,y=0.5cm]
%--Palm(5,2)---------------------------------------
\draw [line width=0.4pt] (1,0)--(1,1)--(1,2);
\draw [line width=0.4pt , dash pattern= on 1pt off 1.3pt] (0,3)--(1,2)--(2,3);
\draw [fill=black] (1,0) circle (2pt);
\draw [fill=black] (1,1) circle (2pt);
\draw [fill=black] (1,2) circle (2pt);
\draw [fill=black] (0,3) circle (2pt);
\draw [fill=black] (2,3) circle (2pt);
%--Lin(4)=Palm(4,1)---------------------------------------
\draw [line width=0.4pt] (4,0)--(4,1)--(4,2);
\draw [line width=0.4pt , dash pattern= on 1.2pt off 2.5pt] (4,2)--(4,3);
\draw [fill=black] (4,0) circle (2pt);
\draw [fill=black] (4,1) circle (2pt);
\draw [fill=black] (4,2) circle (2pt);
\draw [fill=black] (4,3) circle (2pt);
\end{tikzpicture}
\end{minipage}

\medskip\smallskip
Put 
$$
\mathcal{T}^+(n,\ell)=\cup_{k=\ell}^{n-1} \mathcal{T}(n,k)
$$
and
$$
\mathcal{T}(n)=\mathcal{T}^+(n,1).
$$
Thus $\mathcal{T}(n)$ is the set of (isomorphism classes of) all rooted trees with $n$ nodes and $\mathcal{T}^+(n,\ell)$ is its subset whose member trees
have at least $\ell$ leaves. 

%It is convenient to introduce also the set $\mathcal{F}(n)$ of (isomorphism classes of) all rooted forests with $n$ nodes and its subsets
%$\mathcal{F}(n,\ell)$ consisting of forests with $\ell$ leaves.

\trivprob{prob:palm-in-poset-notation}
The tree $\Pl(n,\ell)$ as a poset (with partial order defined in Sec.~\ref{sec:trees}) is isomorphic to
the poset $\mathbf{h}\oplus \ell\xone$ in notation of
\cite[Sec.~3.1,~3.2]{Stanley_1999}, where $h=n-\ell$.
%
%Express $\Pl(n_1,\dots,n_r;\ell_1,\dots\ell_r)$ in terms of direct sums and ordinal sums of elementary posets.

\trivprob{prob:min-palm-tree}
Show that 
\beq{min-palm-tree}
 m_{\Pl(n,\ell)}(\xtt)=(n-\ell+1)\langle\xtt,\xone\rangle^{\frac{1}{n-\ell+1}}.
\eeq
In particular,
\beq{min-palm-tree-homo}
 m_{\Pl(n,\ell)}=(n-\ell+1)\,\ell^{\frac{1}{n-\ell+1}}.
\eeq

\subprob{prob:min-palm-tree-eqneib}
Find all values of $n$ and $\ell$ for which 
$m_{\Pl(n,\ell)}=m_{\Pl(n,\ell+1)}$.

\separline	
In problems {\bf 59}--{\bf 62} we explore facts related to {\em maximization}\ of $m_T$ over classes of trees.

%Problems {\bf 58}--{\bf 60} deal with facts related to {\em maximization}\ of $m_T$ over classes of trees.

\prob{prob:mintree-to-max} (Linear tree maximizes $m_T$ in $\mathcal{T}(n)$.) 
For any tree $T\in\mathcal{T}(n)$ the inequality $m_T\leq n$ holds. The equality occurs if and only if $T=\mathrm{Lin}(n)$. 

\prob{prob:mintree-to-max-ell}\ (An upper estimate for $m_T$ in $\mathcal{T}(n,\ell)$.)
Let $x$ be a positive real number and $\ell\geq x$. Show that for any tree $T\in\mathcal{T}(n,\ell)$ the inequality
$m_T\leq n-\ell-1+x+\ell/x$ holds. 

Some corollaries:

(i) $m_T\leq n-\ell-1+2\sqrt{\ell}$ for 
$T\in\mathcal{T}(n,\ell)$;

(ii)
%In particular, 
$m_T\leq n-\ell/2$ for $T\in\mathcal{T}^+(n,4)$. 

(iii)
%More generally,
$m_T\leq n-(k-1)^2$ for $T\in\mathcal{T}^+(n,k^2)$.

\subprob{prob:mintree-to-max-sqrk}\ (Attainable upper bound.)
Show that $\max\limits_{T\in\mathcal{T}(n,k^2)} m_T=n-(k-1)^2$
if $n\geq k^2+1$.

\subprob{prob:mintree-to-max-not-unique} 
There exist non-isomorphic rooted trees $T_1$, $T_2$ with $|\VV(T_1)|=|\VV(T_2)|=n$, $|\Lambda(T_1)|=|\Lambda(T_2)|=\ell$
and such that $m_{T_1}=m_{T_2}=\max_{T\in\mathcal{T}(n,\ell)} m_T$.
Moreover, the trees $T_1$ and $T_2$ satisfying the stated requirements can be chosen to be non-isomorphic even as unrooted trees. 

\hardsubprob{prob:mintree-to-max-not-unique-odd} 
Does there exist a pair of trees satisfying the requirements of Problem~\bnref{prob:mintree-to-max-not-unique} with  $\ell$ a non-square integer?

\prob{prob:mintree-to-max-lbnd}\ 
(A lower estimate for $\max_{T\in\mathcal{T}(n,\ell)} m_T$.)
Prove that for any $\ell\leq n-1$ there exists a tree $T\in\mathcal{T}(n,\ell)$ such that $m_T\geq 2\sqrt{n-1}=m_{\Pl(n,n-1)}$.

\prob{prob:mintree-to-max-lbnd2}\ 
(A more precise lower estimate.)
If $k^2\leq\ell<\max((k+1)^2,n)$, then there exists a tree
$T\in\mathcal{T}(n,\ell)$ such that 
$
 m_T\geq n-\ell-1+2k
$. The inequality is strict unless $\ell=k^2$.

\subprob{prob:mintree-to-max-doublebnd}\ 
(Corollary: a double-sided bound.)
\beq{maxm-nl}
n-\ell-1+2\lfloor\sqrt{\ell}\rfloor \leq  \max\limits_{T\in\mathcal{T}(n,\ell)}m_T
\leq n-\ell-1+2\sqrt{\ell}.
\eeq

\hardsubprob{prob:max-min-monotone}
The function $\ell\mapsto\max\limits_{T\in\mathcal{T}(n,\ell)} m_T$
decreases for any fixed $n$.

%%%
\separline
Now we turn to {\em minimization}\ of $m_T$.

Let us say that $T$ is a {\em low-branching tree}\ if there is at most one node in $T$ with outdegree $>1$.
We say that a tree $T$ has {\em almost constant height}\ if the lengths of different maximal chains in $T$ (provided there is more than one maximal chain)
differ by at most one. (Constant height trees are, according to this definition, a particular case of trees of almost constant height.)

\prob{prob:mintree-to-min-local} (Low-branching trees as extremizers). 
Let $T\in\mathcal{T}(n,\ell)$ and $t\in\RR$. There exists a low-branching tree $\tilde T\in\mathcal{T}(n,\ell)$
such that
\beq{low-branching-better-than-any}
m_{\tilde T}(t\xone_\ell) \leq m_T(t\xone_\ell). 
\eeq
      
\prob{prob:mintree-to-min-special} (Extremizers in $\mathcal{T}(n,\ell)$). 
Let $n>\ell\ge 1$. For any $t>0$, 
there exists a low branching tree   $T_*\in\mathcal{T}(n,\ell)$ of almost-constant height such that 
 $$
 m_{T_*}(t\xone)=\min_{T\in\mathcal{T}(n,\ell)} m_T(t\xone).
 $$ 
In other words, a minimizer of the functional $T\mapsto m_T(t\xone)$ in  $\mathcal{T}(n,\ell)$ can be found among low branching trees  of almost-constant height. 
        
\prob{prob:mintree-to-min-local-lb}
Let $T\in\mathcal{T}(n,\ell)$. Define
$$
R(x)=\frac{\ln x}{x-1}.
$$
(It is assumed that $R(1)=1$.)

(a) If $t$ is sufficiently small, so that $\ln t\le -n R(\ell)+1+(n-1)/\ell$, then
\beq{mintree-local-lb-small-t}
\ln m_T(t\xone)\ge \frac{\ell \,\ln(t/\ell)}{n+\ell-1}+\ln(n+\ell-1).
\eeq    

(b) If $\ln t>-nR(\ell)$, then
\beq{mintree-local-lb-large-t}
  \ln m_T(t\xone)\ge -R(\ell) +1+\ln(\ln t+nR(\ell) ).
\eeq        

\smallskip
(The conditions on $t$ in (a) and (b) overlap; if both are applicable, then (a) gives a tighter bound.
It is instructive to check the case $\ell=1$.) %, but (b) is a little simpler.)
    
\trivsubprob{prob:mintree-to-min-local-lb-homo}    
As a particular case ($t=1$), for any $T\in\mathcal{T}(n,\ell)$ the inequality
$$
 \ln m_T\ge  \ln n+\ln R(\ell)-R(\ell)+1
$$    
holds. Another form of this estimate is
\beq{min-tree-homo-ml}
 m_T\ge C_\ell n ,
\quad
C_\ell=R(\ell)e^{1-R(\ell)}.
\eeq

\trivsubprob{prob:mintree-to-min-homo-2}
For $\ell=2$ the constant in the inequality \eqref{min-tree-homo-ml} is $C_2=(e\ln 2)/2\approx 0.9421$. Show that it is optimal
(cannot be replaced by a larger number).

Hint: Consider a family of symmetric trees with two leaves.
                        
\prob{prob:mintree-to-min-special-ge} (Palm trees as extremizers in $\mathcal{T}^+(n,\ell)$). 
%
%For any forest $F\in\mathcal{F}(n,\ell)$ and any positive product vector $\xtt=(t_\lambda)$, $\lambda\in\Lambda(F)$, the 
%
Let $T\in\mathcal{T}(n,\ell)$ and $t\in\RR$. There exists $k\ge \ell$ such that
\beq{palm-better-than-any}
m_{\Pl(n,k)}(t\xone_k) \leq m_T(t\xone_\ell). 
\eeq
In other words, $\min_{T\in\mathcal{T}^+(n,\ell)} m_T(t\xone)$ is provided by $T=\Pl(n,k)$ with some $k\ge\ell$. 
Moreover, this minimizing tree is unique in $\mathcal{T}^+(n,\ell)$.

\subprob{prob:mintree-local-counterex} 
Examination of trees with small number of nodes may lead one to conjecture that the tree $\Pl(n,\ell)$ always minimizes $m_T$ in the class $\mathcal{T}(n,\ell)$. Disprove this conjecture
by as simple argument as you can (without reference to \bref{prob:mintree-to-min-local}--\bref{prob:mintree-to-min-local-lb}).

\prob{prob:mintree-to-min-global-homo} (Minimization of $m_T$ in $\mathcal{T}(n)$).
In the homogeneous minimization problem,
$$
\min_{T\in\mathcal{T}(n)} m_T = \min_{1\leq \ell\leq n-1} (n-\ell+1)\,\ell^{\frac{1}{n-\ell+1}}.
$$

\subprob{prob:mintree-to-min-global-explicit} (Convenient explicit lower bound for $m_T$ in $\mathcal{T}(n)$).
For any $T\in\mathcal{T}(n)$
\beq{log-lower-bound-mT-trees}
m_T> e\,\ln(n-\ln n).
\eeq
Along with this lower bound, the following asymptotic formula holds:
$$
\min_{T\in\mathcal{T}(n)} m_T = e\,\ln  n+O\left(\frac{1}{\ln n}\right)
\;\;\mbox{\rm as $n\to\infty$}.
$$

\trivcsubprob{prob:mintree-global-from-local}
Show that inequality \eqref{min-tree-homo-ml} implies the asymptotic lower bound, which is a little weaker than \eqref{log-lower-bound-mT-trees}, but much easier to prove:
$$
\min_{T\in\mathcal{T}(n)} m_T \ge e\ln n-O\left(\frac{(\ln n)^2}{n}\right).
$$

\separline
In line with our systematic exploration of the subject, it is appropriate to ask about statistical properties of the function $m_T$ as a random variable on a suitable probabilistic version of the set $\mathcal{T}(n)$ or $\mathcal{T}(n,\ell)$. It is not easy though to formulate a good concrete question. 

In view of \bref{prob:mintree-to-max-lbnd2}
and \bref{prob:mintree-to-min-global-homo}, the region 
$n-\ell\sim \ln n$ is of special interest, since there $\min m_T=O(\ln n)$, while $\max m_T\sim 2\sqrt{n}$.

In order to define a probabilistic model that would look manageable, we propose to restrict attention to the class $\Pl^*(n,\ell)$ of ``palm bushes''. Trees of this class have been useful in the proofs of our preceding estimates;
in particular, the extremal values of $m_T$ over
$\Pl^*(n,\ell)$ have the same asymptotics as the extremal values over $\mathcal{T}(n,\ell)$.
 
The general element of $\Pl^*(n,\ell)$ is
the tree $T=\Pl(n_1,\ell_1;\dots;n_m,\ell_m)$ --- 
an {\em amalgamated direct sum over the root of the trees $\Pl(n_i+1,\ell_i)$, $1\leq i\leq m$}. Specifically:
$(\ell_1,\dots,\ell_m)$ is a partition of $\ell$,
$(n_1,\dots,n_m)$ is a partition of $n-\ell-1$,
and $T$ is the union of the trees $\Pl(n_i+1,\ell_i)$
sharing the common root.

The probability measure on the set $\Pl^*(n,\ell)$
can be introduced by defining the probability measure on the subset of the Cartesian product of the sets of partitions of $\ell$ and $n-\ell-1$ consisting of the pairs of partitions with equal numbers of parts.

\hardprob{prob:randompalmbush} Consider the set $\Pl^*(n,\ell)$ of ``palm bushes'' as a probability space.
(Define the probability measure precisely.) The functional $T\mapsto m_T$ becomes a random variable on this space. Find the
asymptotics of $\prE m_T$ as $n\to \infty$ under some
assumptions on $\ell$ such as, for instance, $\ell=o(1)$
or $\ell/\ln n\to\const$.

%%%%%%%%%%%%%%%%%%%%%%

\subsubsection{Extremal problems for arc-weighted graphs} 
\label{ssec:graphs-extremal}

Let $\mathfrak{G}_m$ denote the set of of all (isomorphism classes of) strongly connected, not necessarily simple (di)graphs with
$m$ arcs and let $\mathfrak{G}_{m,n}$ ($n\leq m$) denote the subset of $\mathfrak{G}_m$ consisting of graphs with $n$ nodes.
By Euler's formula (see %\bref{prob:Euler_formula}), 
\brefB{prob:cycle_basis_linalg}),
the number of basis circuits for any graph from $\mathfrak{G}_{m,n}$ is the
same, $r = m-n+1$. Denote also $\mathfrak{G}^-_{m,n}=\cup_{k\le n}\mathfrak{G}_{m,k}$ (the set of graphs with $m$ arcs and at most $n$ nodes).
In this section we explore extremal values of the function $f_{\Ga}$ on the above defined sets. 

The problems of minimizing $m_T$ (in Sec.~6.1) and $f_\Gamma$ turn out to be closely related. On the other hand, maximization of $f_\Gamma$ seems to be a rather shallow subject, unlike its counterpart for trees, due to the result {\bf 70}.   

\prob{prob:maxvalue-graph} (Compare \bref{prob:mintree-to-max}.)
Show that 
$$
\max_{\Ga\in\mathfrak{G}_m} f_{\Ga}=m.
$$

\prob{prob:maxgraph} 
The equality $f_\Ga=m$ for $\Ga\in\mathfrak{G}_m$ holds if and only if 
$\Ga$ is a (directed) Eulerian graph.

\bigskip
We pass now to extremal problems concerning minimization of $f_\Ga$.

\medskip
Let $m\ge n\ge 2$. We call a graph $\Ga\in\mathfrak{G}_{m,n}$ {\em special}\ if there exist two nodes $v_*$ (``source'') and $v^*$ (``sink'') such that

(i) $v_*\neq v^*$,

(ii) $v\neq v_*\;\Rightarrow\; d^+(v)=1$,

(iii) $v\neq v^*\;\Rightarrow\; d^-(v)=1$.

\smallskip
In the solutions, whenever there is a need to assign numerical values (components of some vector $\yy\in\RR^{|\VV|}$) to nodes of the graph, we write $y^*$ and $y_*$ instead of
$y_{v^*}$ and $y_{v_*}$ %respectively 
for visual convenience.

\smallskip
The subset of (isomorphism classes of) special graphs in $\mathfrak{G}_{m,n}$ will be denoted $\mathfrak{G}^*_{m,n}$.

\smallskip
A special graph can be thought of as an electric circuit,
where $v^*$ and $v_*$ are battery terminals, several paths outgoing from $v_*$ and arriving at $v^*$ represent 
parallel multi-element loads connected to the battery, and the unique path from 
$v^*$ to $v_*$ represents inner current in serially connected battery cells.

\prob{prob:mingraph-to-special}
Let $m\ge n\ge 2$.
If $\Ga\in\mathfrak{G}_{m,n}$, then there exists $\tilde\Ga\in\mathfrak{G}^*_{m,n}$ with $f_{\tilde\Ga}\le f_{\Ga}$.
In other words, 
$$
\min_{\Ga\in\mathfrak{G}_{m,n}}f_{\Ga}=\min_{\Ga\in\mathfrak{G}^*_{m,n}}f_{\Ga}.
$$

\prob{prob:mingraph-to-tree} The homogeneous extremal minimization problems for graphs and rooted trees are related as follows. 
For any $m\geq n\geq 2$ 
$$
\min_{\Ga\in\mathfrak{G}_{m,n}}f_{\Ga}=\min_{T\in\mathcal{T}(m,\ell)} m_T,
$$
where $\ell=m-n+1$.

\bigskip
Due to \bref{prob:mingraph-to-tree}, the results obtained for trees
 have their counterparts for strongly connected graphs. We state the analogs of
\brefA{prob:mintree-to-min-local-lb}
 %\bref{prob:mintree-to-min-local-lb-homo}
and 
\brefA{prob:mintree-to-min-global-homo}.
%\bref{prob:mintree-to-min-global-explicit}. %The trivial case $n=1$, where there are no corresponding trees, is also included.

%%%
\trivprob{prob:mingraph-local-lb-homo}    
If $G\in\mathfrak{G}_{m,m}$, then $f_\Ga=m$. And if $1\le n\le m-1$, then
for any $G\in\mathfrak{G}_{m,n}$ the inequality
\beq{prob:minextr-graph-mn}
 f_{\Ga}\ge e  m \lambda\,(m-n+1)^{-\lambda}\,\ln(m-n+1),
\qquad
\lambda=1+\frac{1}{m-n}.
\eeq
holds. 

%\subprob{prob:mintree-to-min-global-homo} (Minimization of $f_\Ga$ in $\mathfrak{G}_{m,n}$).
%In the homogeneous minimization problem,
%$$
%\min_{\Ga\in\mathfrak{G}_{m}} m_T = \min_{1\leq \ell\leq n-1} (n-\ell+1)\,\ell^{\frac{1}{n-\ell+1}}.
%$$

\thsubprob{prob:minextr-graph-global} 
(Explicit lower bound for $f_\Ga$ in $\mathfrak{G}_{m}$).
For any $\Ga\in\mathfrak{G}_{m}$
\beq{log-lower-bound-graphs}
f_{\Ga}> e\,\ln(m-\ln m).
\eeq
Along with this lower bound, the following asymptotic formula holds:
$$
\min_{\Ga\in\mathfrak{G}_{m}} f_{\Ga} = e\,\ln  m+O\left(\frac{1}{\ln m}\right)
\;\;\mbox{\rm as $m\to\infty$}.
$$

%%%%
\subsubsection{A ``2021'' problem} %{A teaser problem} 
\label{ssec:teaser-minsum}

The reader who arrived here fresh following our suggestion in the introduction is advised to try parts (a) and (b), which
are more or less recreational problems. Having solved (a), you may want to take a look at Section~\ref{ssec:sums-of-quotients}.
And those who figure out (b) should not be too surprised about estimates like \eqref{log-lower-bound-graphs}. 
 
As regards part (c), it is a peculiar application of the theory developed above, but can also be attempted from scratch. Even if your attempt will not succeed quickly, some constructions employed in solutions of theoretical problems of \S\S~\ref{ssec:trees-extremal}--\ref{ssec:graphs-extremal}
should make more sense.
%and our solution relies on one of the main resluts of that theory. 

\prob{prob:extrval-teaser}
Let $\omega(1),\dots \omega(n)$ be nonempty subsets (not necessarily distinct) of $I=[1:n]$ and $\cup_{i=1}^n \omega_i=I$.
We say that $\omega$ is an {\em assignment}\ (of sets).
If $\xx=(x_1,\dots,x_n)$ is a vector with  positive components, denote
$$
 Y_{\omega}(\xx)=\sum_{i\in I} \frac{x_i}{\min(x_j\,|\,j\in\omega(i))}.
$$
An assignment $\omega(\cdot)$ is {\em irreducible}\ if there is no proper subset $\Sigma\subsetneq I$ such that
$\cup_{i\in\Sigma}\,\omega(i)\subset\Sigma$.

\smallskip
For example, if $n=5$, then the assignment $\omega$ given by $\omega(1)=\{2\}$, $\omega(2)=\{1,3\}$, $\omega(3)=\{4,2\}$,
$\omega(4)=\{5\}$, $\omega(5)=\{1,4,5\}$ is irreducible, but redefining $\omega(5)=\{4,5\}$ leads to a reducible assignment, since the subset $\Sigma=\{4,5\}$ contains both $\omega(4)$ and $\omega(5)$. 

\smallskip
 Now take $n=2021$. 

\smallskip
(a) Find an assignment of sets $i\mapsto \omega(i)$ and $x_1,\dots,x_n>0$ such that
$$
 Y_{\omega}(\xx)<2.021.
$$ 

(b) Find an irreducible assignment $\omega$ and $x_1,\dots,x_n>0$ such that
$$
 Y_{\omega}(\xx)<21.
$$ 

(c)  Prove that if $\omega$ is an irreducible assignment, then for any positive $x_1,\dots,x_n$
$$
 Y_{\omega}(\xx)> 20.
$$ 

%%%%

\subsection{Shallit-type problems}
\label{sec:shallit}

% Publish separately under the title 
% On Shallit's problem concerning minimization of sums of fractions

The following problem was proposed by J.~Shallit \cite{Shallit_1994} in 1994.

%\medskip
\prob{prob:Shallit}
{\it
Let $\xx=(x_1,\dots,x_n)$ be a vector with positive components.
Denote
$$
 S_n(\xx)=\sum_{i=1}^n x_i+\sum_{1\le i\le j\le n} \prod_{k=i}^j\frac{1}{x_k}.
$$
Show that there exists a positive constant $C$ such that
$$
 \min_{\xx>0} S_n(\xx)=3n-C+o(1)
$$
as $n\to\infty$ and find the numerical value of $C$.
}

\medskip
In the Solutions section we give comments and references. 
%on a solution without giving it in full. 
{Problems~{\bf 76}--{\bf 77}}
reveal connection of Shallit's problem to our present subject.

\prob{prob:Shallit-as-particular-case}
Reduce the problem of finding $\min_{\xx}S_n(\xx)$ to the minimization problem of the form \eqref{homo1min}. That is,
exhibit a graph $\Ga_n$ such that 
$\min_{\xx>0} S_n(\xx)=f_{\Ga_n}$.
%\end{prob}

\prob{prob:Shallit-alt-graph} (Simplified graph for Shallit's problem.)
The number of arcs in the graph $\Ga_n$ (and the number of summands in Shallit's sum)
is of order $n^2$. 
Find a graph $\Ga'_n$ with $n+1$ nodes and $O(n)$ arcs such that
$f_{\Ga_n}=f_{\Ga'_n}$.

\hardsubprob{prob:Shallit-alt-generalization}
Provide an explanation to the existence of the transformation $\Gamma_n\to\Gamma'_n$.
Are there other (families of) examples of pairs of graphs
with different number of arcs but the same values of the functional $f_\Gamma$? 

\hardprob{prob:Shallit-no-double-edges}
(A variant of Shallit's problem.)
Define
$$
 \hat S_n(\xx)=\sum_{i=1}^n x_i +\sum_{1\le i< j\le n} \prod_{k=i}^j\frac{1}{x_k}=S_n(\xx)-\sum_{i=1}^n \frac{1}{x_i}.
$$
Show that there exist positive constants $\hat\lambda$ and $\hat C$ such that
$$
 \min_{\xx>0} S_n={\hat\lambda}n-\hat C+o(1)
$$
as $n\to\infty$ and find the values of $\hat\lambda$ and $\hat C$.

\prob{prob:Shallit_general_pattern-graph}
(Generalized Shallit's problem.)
Let $\Patt$ be any nonempty (finite or infinite) subset of $\ZZ_{\geq 0}$ (a ``pattern'').
Denote
$$
 S_n(\xx|\Patt)=\sum_{i=1}^n x_i +\sum_{j-i\in \Patt} \prod_{k=i}^j \frac{1}{x_i}
$$
and
$$
 m_n(\Patt)=\min_{\xx>0} S_n(\xx|\Patt). 
$$
Exhibit a graph $\Ga_n(\Patt)$ with $n$ nodes such that
$$
 m_n(\Patt)=f_{\Ga_n(\Patt)}.
$$

\prob{prob:Shallit_general_pattern-asymp}
(Crude asymptotics in the generalized Shallit's problem.)
Let $m_n(\Patt)$ be defined as in \bref{prob:Shallit_general_pattern-graph}.
Let $\rho_\Patt$ be the unique root of the equation
\beq{gen-Shallit-coef}
 \rho=\sum_{n\in \Patt} (n+1)\rho^{-n-1}.
\eeq
in $(1,\infty)$.
Then
\beq{crude-genShallit}
 \lim_{n\to\infty}\frac{m_n(\Patt)}{n}=\lambda_\Patt,
\eeq
where
\beq{genShallit-lambda}
\lambda_\Patt=\rho_\Patt+\sum_{n\in \Patt}\rho_\Patt^{-n-1}=\sum_{n\in \Patt}(n+2)\rho_\Patt^{-n-1}.
\eeq

\subprob{prob:Shallit-root-part-cases}
Determine the values of $\rho_\Patt$ and $\lambda_\Patt$ for the Shallit's pattern
$\Patt=\ZZ_{\ge 0}$ and the pattern of Problem~\bref{prob:Shallit-no-double-edges}, $\Patt=\ZZ_{\ge 1}$.

%%%%%%%%%%%%%%%%%%%%%%%%%%%%%%%%%%%

\hardsubprob{prob:Shallit_general_pattern_refined}
Prove or disprove: for any pattern $\Patt$ there exists the limit
$$
 C_\Patt=\lim_{n\to\infty} (m_n(\Patt)-n \rho_\Patt).
$$   

\iffalse
[\sol{prob:Shallit_general_pattern_refined}  (Unsuccessful)
Suppose $\Patt$ is finite. (The general case follows from this
particular.) More precisely, let $T\subset[0:M-1]$.
We will introduce a dynamical system on the space $\RR^M$
and investigate its limit behaviour for a one-parameter
family of initial conditions.

Let $Y=[y_1,\dots,y_M]$ be the vector in the phase space.
]
\fi

\hardprob{prob:Shallit-Kn-random-directions}
(Shallit-type problems for random configurations.)
Consider the (non-oriented) complete graph $K_{n+1}$
with vertices labeled $0,1,\dots,n$.
Define a random oriented graph $\Ga_n$ by assigning 
directions to the edges of $K_{n+1}$ as follows.
First, the ``contour'' of $\Ga_n$ will have a deterministic orientation: it involves the arcs $(0\to 1), (1\to 2),\dots, (n\to 1)$. Every diagonal of $K_{n+1}$ will be assigned either of the two possible directions independently with equal probability $1/2$. Thus for any $i$ and $j$ such that $i-j\mod (n+1)\notin \{0,\pm 1\}$ 
$$
 \Prob((i\to j)\in\AA(\Ga_n))=
 \Prob((j\to i)\in\AA(\Ga_n))=
 \frac{1}{2}.
$$
Find the asymptotics of the expectation
$
 \prE(f_{\Ga_n})
$
as $n\to\infty$.

%%%%%%%%%%%%%%%%%%%%%%%%%%%%%%%%%%%%%%%%%%%%%%%%%%%%%%%%%%%%%%%%%%%
%\section{Summary}
%
%1. Critical point equations (necessary condition for extremum).
%
%For strongly connected graphs, extremum --- at $\xx>0$.
%
%Existence, uniqueness of solution.
%
%2. Monotonicity: (a) of min value; (b) of minimizing vector $\xx$ (componentwise)
%with respect to: (i) order $\tt>\tt'$ componentwise; (ii) dilatation-order (all orders --- partial)
%
%Counterexamples.
%
%Function $f(\xtt)$ is strongly monotone if $\xtt\leq \xtt'$ (componentwise) implies $f(\xtt)\leq f(\xtt')$.
% Weakly monotone (dilatation-monotone), if $f(\xtt)\leq f(c\xtt)$ for $c\geq 1$.
%
%Value of minimum: strongly monotone (obvious, by comparing strengths of constraints).
%
%Constant height tree: strong monotonicity of every weight $x_e(\xtt)$.
%
%Any tree: weak monotonicity of every weight.
%
%General strongly connected arc-weighted graph: components are not even weakly monotone. ("Great example" with denominators 43). 
%
%Example showing that for general tree $x_e(\xtt)$ may be not strictly-monotone (two-arc tree)
%
%
%Algorithm of computation of min for tree-sums.
%
%Extremal problems.
%

\clearpage

\addtocounter{secaux}{1}
\section%{Part II. Solutions and comments}
{Solutions and comments}
\label{part:solutions}
%\addcontentsline{toc}{section}{Part II. Solutions and comments}
%\addtocontents{toc}{\protect\contentsline{section}{\protect{}Part II. Solutions and comments}{}{}}

\subsection*{Section~\ref{sec:general}}
\addcontentsline{toc}{subsection}{Section~1}

\sol{prob:consfeasible}\
By the linear independence of the rows of the matrix $A$, its row rank  equals $r$, hence its column rank is also $r$ and the range of $A$ has dimension $r$.

\sol{prob:noncompact-example}\
This is an example with two variables and one constraint:  $x_1 x_2^{-1}=1$. Clearly, $\inf (x_1+x_2)=\lim\limits_{x_1=x_2\to 0}(x_1+x_2)=0$ and there is no minimizer.

\sol{prob:change-weight}\
If $\xx$ is an admissible vector for the constraints $A\log\xx=\log\xtt$ and $x_i' w_i'=x_i w_i$, then the vector $\xx'$ is an admissible 
vector for the constraints $A\log\xx'=\log\xtt'$, and $\langle\ww,\xx\rangle= \langle\ww',\xx'\rangle$.

%\sol{prob:relaxcons1}  (Relaxation of constraints). 

%\sol{prob:monotonicity}  (Monotonicity with respect to weights). 

\solL{prob:monotonicity}{A}\  
For any admissible vector $\xx$ 
we have $\alpha\langle \ww',\xx\rangle+
(1-\alpha)\langle \ww'',\xx\rangle\geq\alpha f(\ww',\xtt',A)+(1-\alpha)f(\ww'',\xtt',A)$. Taking infimum
over $\xx$ yields the result.

\medskip\noindent{\bf Remark.} The results of the type
\bnref{prob:monotonicity}, \bnref{prob:concave-weights} are true
in any optimization problem of the form $\langle \ww,\xx\rangle\to \min$ under constraints $\xx\in\mathcal{D}$, where the set $\mathcal{D}$ of admissible vectors does not depend on $\ww$.

\sol{prob:notmonotone-wrt-products}
Let $n=2$, the constraints be $x_1=t_1$ and $x_1 x_2=t_2$. Then $\langle \ww,\xx\rangle=w_1 t_1 +w_2t_2/t_1$,
so $(d/dt_1)\langle \ww,\xx\rangle=w_1-w_2t_2/t_1^2$, which can be negative. See also \bref{prob:ex-quasibasis}.

\sol{prob:consineq}\ Relaxation of constraints cannot increase the lower bound.
To show that strict inequalities are possible, let us take an example as in \bref{prob:notmonotone-wrt-products}, where $f(\ww,\xtt',A)<f(\ww,\xtt,A)$ and $\xtt\leq\xtt'$.
There exists a vector $\xx'$ 
such that $A\log\xx'=\xtt'$
and $\langle \ww,\xx'\rangle<f(\ww,\xtt,A)$.
The vector $\xx'$ is feasible for the system of constraints 
$A\log\xx\geq\xtt$.
Hence $f_\geq(\ww,\xtt,A)<f(\ww,\xtt,A)$.

Taking any trivial example with $A_{ji}>0$,  $f(\ww,\xtt',A)<f(\ww,\xtt,A)$ and $\xtt'\leq\xtt$,
we similarly see that $f_\leq(\ww,\xtt,A)<f(\ww,\xtt,A)$
is possible.

\sol{prob:zeroweight}\ 
First, $f(\hat\ww,\hat\xtt,\hat A)\leq \inf_{w_1>0} f(\ww,\xtt,A)$ by \bref{prob:relaxcons1}. 
Let us prove the reversed inequality.
Suppose that $\hat\xx_\epsilon$ is an approximate $\epsilon$-minimizer for the truncated problem with constraint matrix $\hat A$ (possibly empty).
That is, $\langle\hat\xx_\epsilon,\hat\ww\rangle<f(\hat\ww,\hat\xtt,\hat A)+\epsilon$.
Define $\xx_\epsilon$ as the unique admissible vector for the non-truncated problem whose truncation (by crossing out the component $x_n$) is $\hat\xx_\epsilon$.     
(Simply put, recover $x_n$ from the $r$-th constraint deleted for truncation.) Then
$\langle\xx_\epsilon,\ww\rangle=\langle\hat\xx_\epsilon,\hat\ww\rangle+w_n x_n$. Hence
$\liminf_{w_n\to 0} f(\ww,\xtt,A)\leq \liminf_{w_n\to 0}\langle\xx_\epsilon,\ww\rangle\leq  f(\hat\ww,\hat\xtt,\hat A)+\epsilon$.
Since $\epsilon>0$ is arbitrary, the inequality $\liminf_{w_i\to 0} f(\ww,\xtt,A)\leq  f(\hat\ww,\hat\xtt,\hat A)$ follows. 

%Remark: The proof does not depend on existence of minimizer.

\sol{prob:logconvex}
If $\xx^{(i)}$ are admissible vectors for the systems of constraints $A\log\xx=\log\xtt^{(i)}$, then the vector
$\xx=\prod_1^k (\xx^{(i)})^{p_1}$
is an admissible vector for the system of constraints
$A\log\xx=\log\xtt$. By H\"older's inequality 
(see e.g.\ \cite[\S~18]{BeckBel_1961}), %we have
$$
 \langle\xx,\ww\rangle\leq  \prod_{i=1}^k\langle\xx^{(i)},\ww^{(i)}\rangle^{p_i}.
$$
Minimizing the right-hand side (more precisely, passing to infimum in each factor independently), we see that
$\inf_{\xx:\,A\log\xx=\xtt}\langle\xx,\ww\rangle\leq \prod_{i=1}^k f(\ww,\xtt^{(i)},A)$.

%\sol{prob:pnorm-convexity}\ 
\solL{prob:pnorm}{A}\ 
Let vectors $\yy$, $\zz$ satisfy the constraints $A\log\yy=\log\xtt$
and $A\log\zz=\log\xtt$. Put $\xx=\yy^{1-\theta}\zz^{\theta}$. Then $A\log\xx=\log\xtt$. 
By H\"older's inequality,
$$
\langle \ww,\xx^{p_\theta}\rangle^{1/p_\theta}\leq 
\langle \ww,\yy^{p_0}\rangle^{(1-\theta)/p_0}
\langle \ww,\zz^{p_1}\rangle^{\theta /p_1}.
$$
The claim follows, since the minimization over $\yy$ and $\zz$ in the right-hand side can be done independently.

%\sol{prob:psum-convexity}\
\solL{prob:pnorm}{B}\ 
{\bf Comments.} 1. Another similar statement is: $p\mapsto F(p)$ is convex if and only if $u\mapsto F(1/u)$ is convex
\cite[Ch.~3, {\bf 119}]{HLP}.

2. In the case of resolved constraints: $r=n$, $x_i=t_i$, $i=1,\dots,n$, --- the inequalities of problems A and B become
the familiar log-convexity property of $p$-norms, which is a reformulation of H\"older's inequality, and the Lyapunoff inequality \cite[Ch.~2, {\bf 18}]{HLP}, \cite[\S~V.3 ]{Mitrinovic_1973}.

\sol{prob:critpteqs}\ The system for the essential variables is obtained by the method of Lagrange's multipliers.
And if $x_i$ is a nonessential variable, then $x_i=0$ in a minimizer. The corresponding equation $w_i x_i=\sum 0\lambda_j$ is trivially satisfied.  

\sol{prob:AMGM}
In this case, the equations \eqref{critpteqs} with a single $\lambda$ variable take the form $w_i x_i=\lambda \rho_i$, and the result
readily follows.

Another way is to derive the inequality with arbitrary $w_i$ and $\rho_i$ from the standard AM-GM inequality by using \bref{prob:change-weight} and \bref{prob:dilatation}.

%\sol{prob:AMGMsmall}
\solL{prob:AMGM}{A}
Let $t'=t x_i^{-\rho_i}$. %Let the vector $\xx^{(i)}$ be obtained from $\xx$ by deletion of the $i$-th component.
Then
$$
\langle\xx,\ww\rangle\geq  \langle\hat\xx^{(i)},\hat\ww^{(i)}\rangle\ \geq f(\hat\ww^{(i)},t',\hat A^{(i)}),
$$
and the result follows by \bref{prob:AMGM}.

\sol{prob:compact-sublevel}
Since the constraints are compact, there exist positive $\rho_1,\dots,\rho_n$ and $\tau\in\RR$ such that $\sum\rho_i \log x_i=\tau$
for any admissible vector $\xx$. 
By \brefA{prob:AMGM}, for every $i=1,\dots,n$ there exists $\epsilon_i>0$
such that $x_i\geq \epsilon_i$ for $\xx\in X_c$.  Thus $X_c$ is the intersection of the closed subset in $\RR_{>0}^n$ defined by the constraints
and the compact set in $\RR_{>0}^n$ defined by the inequalities $\epsilon_i\leq x_i\leq c/w_i$ ($i=1,\dots,n$). Since $c<f(\ww,\xtt,A)$, the set $X_c$ is non-empty.

%\sol{prob:almost-compact-means-nonzero}
\solL{prob:compact-sublevel}{A}\
We may assume that the first row $A_1=(A_{1i})$ of the matrix $A$ is a nonnegative, nonzero vector.
By \bref{prob:relaxcons1} and \bref{prob:AMGM}, $f(\ww,\xtt,A)\geq f(\ww,t_1,A_1)>0$.

\smallskip\noindent
{\bf Remark.} A system of constraints satisfying the condition of this problem is not necessarily compact.
For instance, there can be two disjoint subsets of variables and independent constraints imposed on them;
compact for the first subset and non-compact for the second.

\sol{prob:existminimizer} 
We may assume that all variables are essential. (Non-essential variables equal 0 in a minimizer.)  

\smallskip
1. Suppose the system of constraints is compact. 
Pick some $c>f(\ww,\xtt,A)$. A minimizer for \eqref{inf1} is a point where the continuous function 
$\xx\mapsto\langle\xx,\ww\rangle$ attains its minimum on the compact non-empty set $X_c$ defined in \bref{prob:compact-sublevel}.

\smallskip
2. Suppose $\xx^*$ is a minimizer.
We will exhibit a linear combination of constraints with positive coefficients.
By \bref{prob:critpteqs}, the vector space spanned by the rows of the matrix $A$ contains the vector with components $w_i x^*_i>0$.  
The required linear combination of constraints is $\sum w_i  x^*_i \log x_i=\tau$ with some $\tau\in\RR$. 

\sol{prob:uniqsol} Suppose $\xx\neq\xx'$ are two minimizers with corresponding Lagrange's multipliers in system \eqref{critpteqs}.
By Lagrange's mean value theorem, $\log x_i'-\log x_i=(x_i'-x_i)/u_i$ with some $u_i$ that lies between $x_i$ and $x_i'$. 
Put $x_i'-x_i=y_i$ and $\lambda'_j-\lambda_j=\mu_j$. The system \eqref{critpteqs} becomes
$$
\ba{l}\dst
  w_i y_i =\sum_{k=1}^r \mu_k A_{ki} \quad (i=1,\dots,n),
\eqline{2}
 \sum_{i=1}^n A_{ji} \frac{y_i}{u_i}=0 \quad (j=1,\dots,r).
\ea
$$
By assumption, $\yy\neq 0$, hence some $\mu_j\neq 0$. Eliminating the $y$-variables we obtain
$$
 \sum_{i=1}^n \sum_{k=1}^r \frac{A_{ki} A_{ji}}{w_i u_i}\,\mu_k=0 \quad (j=1,\dots,r).
$$
We may assume that there are no non-essential variables (equivalently, all non-essenital $x_i=x_i'=0$ can be thrown out).
Then $w_i u_i>0$ and the formula
$$
 (\xvv,\xvv')=\sum_{i=1}^n \frac{v_i v'_i}{w_i u_i}
$$
defines a positive-definite inner product in $\RR^n$.
The $r\times r$ matrix $G$ with entries
$$
 G_{jk}=\sum_{i=1}^n  \frac{A_{ki} A_{ji}}{w_i u_i} 
$$
is the Gram matrix of the set of $r$ rows of the matrix $A$ with respect to this inner product.
Since the rows of $A$ are linearly independent, $\det G\neq 0$.  Then equation $G\vec{\mu}=0$ implies $\vec{\mu}=0$, a contradiction.

\sol{prob:nonuniqsol-polynomial}
Take $\mathbf{g}(x_1,x_2) = (x_1+x_2^2,x_1^2+x_2)$. Then $\mathbf{g}(1,0) = \mathbf{g}(0,1) = (1,1)$.

\sol{prob:symmetry}\
Note: this result critically depends on the uniqueness of minimizer.

\sol{prob:gradient_flow}
Let us first check that for any $s>0$ the vector $\xx(s)$ is admissible. 
Indeed, 
$$
 \frac{d (A\log\xx)_j}{ds}=\sum_{i=1}^n A_{ji}\frac{\dot x_i}{x_i}=\sum_{i=1}^n A_{ji} (P\xx-\xx)_i=0,
$$
since $\xx-P\xx$ is the orthogonal projection of $\xx$ onto $(\Ran A^T)^\bot=\ker A$.
 
Now,
$$
  \frac{d\langle\xx,\xone\rangle}{ds}=\langle P\xx-\xx,\,\xx\rangle\leq 0.
 $$
Hence, for any $s>0$, by \bref{prob:compact-sublevel}  $\xx(s)$ lies in the compact set $X_c$ with $c=\langle\xx_0,\xone\rangle$.
If $\xx$ is a fixed point, then $\langle P\xx,\xx\rangle=\|\xx\|^2$, that is, $\xx\in\Ran A^T$. Then $\xx$ satisfies the system
\eqref{critpteqs}, hence it is the minimizer.

\smallskip
\noindent
{\bf Remark.}
The dynamical system defined in this problem is not suitable for a practical computation of the minimizer as the manifold of admissible vectors
is not stable under the flow. 

\sol{prob:newton_method}\
(a) Eliminating the unknowns $v$ and $u$, we obtain a linear system to solve for $\{\mu_k\}_1^r$ with $r\times r$  nondegenerate matrix 
$G_{jk}=\sum_{i=1}^n A_{ji}A_{ki}(w_i x_i)^{-1}$.  % Cf.\ solution \ref{prob:uniqsol}. 
A simple check shows that $\xx+\xuu\in N_x$ and $\yy+\xvv\in M_y$.

\smallskip
(b) If $(u_i)$, $(v_i)$ are determined from the system \eqref{newton-iter} and 
$\tilde\eps_i=(y_i+v_i)-\log(x_i+u_i)$, then
$$
 \tilde\eps_i=(y_i-\log x_i)+\left(v_i-\frac{u_i}{x_i}\right)+\left(\frac{u_i}{x_i}-\log\left(1+\frac{u_i}{x_i}\right)\right) =\frac{u_i}{x_i}-\log\left(1+\frac{u_i}{x_i}\right).
$$

By Eqs.\ \eqref{critpteqs}, $\langle\xx,\ww\rangle=f(\ww,\xtt',A)$, where
$\log\xtt'=\log\xtt-\sum_i A_{ji}\eps_i$.
Hence if $\|\eps\|$ is small, then the pair $(\xx,\yy)$ is close to $(\xx_*,\log\xx_*)$, where $\xx_*$ is the minimizer. In particular, the norm of the inverse matrix $G^{-1}$ in (a) is uniformly bounded.
Therefore there exist sufficiently large $C_1$ and $C_2$ such that
$\max_i |u_i/x_i|\leq C_2\|\eps\|$ whenever $C_1\|\eps\|<1$.

Take $C=\max(C_1, 2C_2, C_2^2)$. Using the fact that $|\delta-\log(1+\delta)|<|\delta|^2$ for $|\delta|\leq 1/2$,
we get: if $C\|\eps\|<1$, then $|u_i/x_i|<1/2$ and 
$$
C\tilde\eps_i<C(u_i/x_i)^2\leq C(C_2\|\eps\|)^2\leq (C\|\eps\|)^2,
$$
as required.

%-------
\sol{prob:Jacobian_for_minimizer} 
%The stated formulas are trivial if the constraints are rigid ($r=n$). So let us assume that $r<n$.  
(a) Suppose first that $\ww$ is held constant and let $\tau_j=\log t_j$. 
Eliminating the $x$-variables in the system \eqref{critpteqs}, we get 
$$
\sum_{k=1}^r \sum_{i=1}^n (w_i x_i)^{-1} A_{ji}A_{ki}\, d\lambda_k=d\tau_j, \quad j=1,\dots,r.
$$
Let $D_{\mathbf{a}}$ denote the diagonal $n\times n$ matrix with diagonal entries $a_i$. We have $B=A D_{(\ww\xx)^{-1/2}} $ and $G=BB^T$.
(Here $(\ww\xx)^{-1/2}$ means the vector with components $(w_i x_i)^{-1/2}$.) 
Then $\vec{d\lambda}=G^{-1} \vec{d\tau}$, so $d\xx=D_{\ww}^{-1}A^T\vec{d\lambda}=D_{\ww}^{-1} A^T G^{-1}\vec{d\tau}$, which is \eqref{dxdt} in a matrix form.

\smallskip
(b) Now suppose that $\xtt$ is held constant and $\ww$ varies. 
%
%Let $\vec{\tau}'=\vec{\tau}+A\log\ww$.  Reduce the minimization problem to the one with constant weight
%vector $\ww'=\xone$ according to \bref{prob:change-weight}. 
%
Along with the problem with data $(\ww,\xtt,A)$ consider the
problem with data $(\xone,\xtt',A)$, where $\log\xtt'=\log\xtt+A\log\ww$. 
Let $\xx'$ be the minimizer for the latter. Then $x_i'=w_i x_i$ (cf.\ solution of \bref{prob:change-weight}).
%We have, by the previous,
According to (a),
$
 d\xx'=A^T G^{-1} \vec{d\tau'}.
$
But $\vec{d\tau'}=A D_{\ww}^{-1}\,d\ww$ and $d\xx'=D_{\ww} d\xx+D_{\xx}d\ww$. Therefore
$$
 d\xx=D_{\ww}^{-1} (A^T G^{-1} AD_{\ww}^{-1}-D_{\xx}) d\ww.
$$
It is a standard fact that $P=B^T G^{-1} B$ is %the matrix of 
the orthogonal projector onto $\Ran B^T$ (the row space of the matrix $B$).  
The orhtogonal projector onto its orthogonal complement $\ker B$ is $Q=I-P$. The obtained formula can be rewritten as
$$
 d\xx=-D_{(\xx/\ww)^{1/2}} Q D_{(\xx/\ww)^{1/2}} d\ww.  
$$ 
The formula \eqref{dxdw} follows.

\smallskip
(c) This result does not depend on a particular form of constraints. Clearly, if any of the parameters of the minimization problem
undergoes a small variation, the corresponding variation $\delta\xx$ of the minimizer must be orthogonal to $\ww$. 
Therefore $\delta\langle\ww,\xx\rangle=\langle\delta\ww, \xx\rangle$ and the result follows.

\sol{prob:noncommutative} From the constraints it follows that the matrix $Y$ has the form
%Put $y_{11}=y_{22}=x$, $y_{12}=y_{21}=y$. 
$Y=\mat{cc}{x & y\\y & x}$ with $2(x+y)=\log t$. The eigenvalues of $Y$ are $x\pm y$. Hence $\Tr\, e^Y=e^{x+y}+e^{x-y}=e^{(\log t)/2}+e^{x-y}$.
%Since the unknowns $x$ and $y$ are bound by just one linear constraint $2(x+y)=\log t$,
%we have 
Since $\inf (x-y)=-\infty$, we conclude that $\inf_Y \Tr\, e^Y=\sqrt{t}$ and there is no minimizer.

\sol{prob:prod-lin-constraints}
Let $x_1=u$, $x_2=v$, $x_3=w$, $x_4=uv$,
$x_5=vw$, $x_6=wu$.
We are asked to find the extrema of the weighted sum
$0x_1+0x_2+0x_3+1x_4+1x_5+1x_6$ under the product constraints
$x_1x_2x_4^{-1}=x_2x_3x_5^{-1}=x_3x_1x_6^{-1}=1$,
$x_1 x_2 x_3=p$, and the linear constraint $x_1+x_2+x_3=s$.

The system of critical point equations
is similar to 
\eqref{critpteqs} and contains one additional Lagrange's multiplier $\mu$ corresponding to the linear constraint:
$$
 \ba{l}
   -(\lambda_1+\lambda_3+\lambda_4) x_1^{-1}-\mu=0, 
\\ 
   -(\lambda_1+\lambda_2+\lambda_4) x_2^{-1}-\mu=0,
\\
   -(\lambda_2+\lambda_3+\lambda_4) x_3^{-1}-\mu=0,
\\
   1+\lambda_1 x_4^{-1}=0,
\\
   1+\lambda_2 x_5^{-1}=0,
\\
   1+\lambda_3 x_6^{-1}=0,   
 \ea
$$
and five equations of constraints.

Eliminating first $\lambda_1$, $\lambda_2$, $\lambda_3$
and then $\lambda_4$ and $\mu$, we reduce the above six  equations to one relations between the $x$-variables:
\iffalse
$$
 \mu=(x_4+x_6-\lambda_4)x_1^{-1}=
 (x_4+x_5-\lambda_4)x_2^{-1}=
 (x_5+x_6-\lambda_4)x_3^{-1}
$$
%
$$
 \lambda_4=
  \frac{(x_4+x_6)x_1^{-1}-(x_4+x_5)x_2^{-1}}{x_1^{-1}-x_2^{-1}}
=\frac{(x_4+x_6)x_1^{-1}-(x_5+x_6)x_3^{-1}}{x_1^{-1}-x_3^{-1}}
$$
%
$$
 ((x_4+x_6)x_3-(x_5+x_6)x_1)(x_2-x_1)=
  ((x_4+x_6)x_2-(x_4+x_5)x_1)(x_3-x_1)
$$
%
\fi
$$
 (x_1-x_2)x_4+(x_2-x_3)x_5+(x_3-x_1)x_6=0.
$$
Returning to the variables $u$, $v$ and $w$ and simplifying, we find
$$
 (u-v)(v-w)(w-u)=0.
$$

Consider the case $v=w$. Expressing $u$ in two ways from the equations of costraints: 
$
 u=s-2v=p/v^2
$,
we get the cubic equation for $v$:
$$
 2v^3-sv^2+p=0.
$$
It always has one negative real root, while the existence of  positive roots depends on the sign of the discriminant: they exist whenever $s^3-27p\geq 0$. 
%(The existence condition is of course but the AM-GM inequality.) -- Not quite, since AM-GM deals with positive values.

To each of the two positive roots $v$ there corresponds the extremal value
$$
f_{\mathrm{extr}}=2u v+v^2=2(s-2v)v+v^2=
 (2s-3v)v.
$$
One of these values is the maximum and the other minimum.

In general (unless $s^3=27p$), the maximum and minimum are not equal and are each attained at three disctinct critical points. There is no uniqueness,
unlike in the case of pure product constraints.
% (\bref{prob:uniqsol}). 
The extremizers are less symmetric than the problem, in contrast with \bref{prob:symmetry}.

%%-----------------------
\subsection*{Section~\ref{sec:graphs-basics}}
\addcontentsline{toc}{subsection}{Section~2}

\sol{prob:ineq-circuits-arcs}\ 
Let $|\AA|=m$. The set of all sequences of arcs of length $k$ without repetitions has cardinality ${m\choose k} k!=m!/(m-k)!$. Since $\sum_{k=1}^m 1/(m-k)!<e$,
the required estimate follows. 

%\sol{prob:ineq-circuits-arcs-sharp}\
\solL{prob:ineq-circuits-arcs}{A}\ 
Observing that to every \emph{cyclic}\ sequence of length $k$ without repetitions there correspond $k$
usual sequences of length $k$, we get   
%\beq{sharp-Eul-est}
$$
|\CC|\leq m!\left(\frac{1}{1(m-1)!}+\frac{1}{2(m-2)!}+\dots+\frac{1}{m\,0!}\right)=(m-1)!(e+o(1)).
$$
%\eeq
Sharpness of the constant $e$ can be seen by looking at the %non-simple 
graph with one node and $m$ loops
or, more generally, at a topologically equivalent to it union of $m$ circulant one-cycle graphs having one common node.
%(Fig.\ \ref{fig:flower-graph}).
(Note: the number of different supports of circuits
%sets that constitute a circuit 
in such a graph is $2^m-1$, cf.\ %\bref{prob:count-circuits-vs-supports}
\brefA{prob:count-circuits-cycles-small-n}.)

\sol{prob:countcycles-Kn}\ 
The number of closed paths 
in $\Knoloops{n}$ of length $k\geq 2$ with a marked origin and no repeated nodes is 
$
 {n-1\choose k-1} (k-1)!=(n-1)!/(n-k)!
$.
There are $n$ ways to choose the marked origin and to every cycle of length $k$ there correspond $k$ closed paths.   
Hence the stated formula for $|\CC_0(\Knoloops{n})|$. In $\Kwloops{n}$, the $n$ cycles of length $1$ (loops) are added. 

\sol{prob:countEul-circuits-Kn}\ 
Let $\gamma$ be an Eulerian circuit in $\Knoloops{n}$. It passes through every node $n$ times. 
Adding the $(n+1)$-th node to $\Knoloops{n}$,
we can insert the cycle $(i\to n+1\to i)$ in $\gamma$ in $n$
ways for each $i=1,\dots,n$. Hence $k_{n+1}\geq k_n n^n$.

The logarithmic asymptotical estimate follows, for example,
by the Euler-Maclaurin formula. 
%(In fact, $\log\prod_{d=1}^{n} d^d=\frac{1}{2}n^2\ln n-\frac{1}{4}n^2+\frac{1}{2}n\ln n+\frac{1}{12}\log n+C+o(1)$.)

\sol{prob:countcircuits-Kn}\
{\bf Comment.} Let $G_n$ be the graph with $\VV=\{(i,j), 1\leq i,j\leq n\}$, and arcs of the form $(i,j)\to (j,k)$.
In other words, the nodes of $G_n$ correspond to the arcs of $\Knoloops{n}$ and the arcs of $G_n$ correspond to 
%the two-arc paths between nodes $i$ and $k$
paths of length 2
in $\hat K_n$. Then $|\CC(\Knoloops{n})|=|\CC_0(G_n)|$. Hence, according to 
%\bref{prob:ineq-circuits-nodes},
\brefA{prob:countcycles-Kn},
$|\CC(\Knoloops{n})|<(e+o(1))(n^2-1)!$.

Using the lower estimate from \bref{prob:countEul-circuits-Kn}, we get the double-sided %asymptotical 
inequality
$$
 \frac{1}{2}-o(1)\leq\frac{\log|\CC(\Knoloops{n})|}{n^2\log n}\leq 2+o(1). 
$$
The author thinks that the left bound is the actual logarithmic asymptotics. 

%In another direction, one can try to find a recurrent formula for $|\CC(\Knoloops{n})|$.

\sol{prob:bases-circuits-to-cycles}
The linear span of the columns of the arc-cycle incidence matrix for the system of cycles $\cup_{L\in\LL}[L]$ contains all the columns of the arc-circuit incidence matrix for the system $\LL$.

%\sol{prob:basis-fungraph}{B}\
\solL{prob:relevant-fungraph}{B}\
{\bf Comment.} 
In the standard terminology of dynamical systems theory, the set $\VV_0$ is the omega-limit set of the adjacency map $\phi$.
The map $\phi:\,\VV\to\VV$ can be extended %(e.g.\ by a piecewise-linear interpolation)
to a continuous self-map of the graph $\Ga$ (with natural topology of a union of segments). 
The paper \cite{HricMalek_2006} lists possible structures of omega-limit sets of continuous self-maps of graphs in general. 
 
\sol{prob:cycle_basis_linalg}
By definition of a basis of curcuits, columns of the incidence matrix for $\LL$ are linearly independent and their span $V$ coincides with subspace in $\RR^\AA$ spanned by all vectors $\xx_L=(x_a, a\in L)$, $L\in \CC(\Ga)$.
We have to prove that $V=U$. It is obvious that $V\subset U$; equivalently, $V\bot\Ran d$, hence, $V^\bot\supset\Ran d$. 
It remains to show that $V^\bot\subset \Ran d$. Let $\xx\in V^\bot$. Let us construct a $0$-cochain $\yy$ such that $d\yy=\xx$.
Pick some node $v_*\in\VV(\Ga)$ and let $y_{v_*}=0$. For any other node $v$ define $y_v=\sum_{a\in\pi} x_a$, where $\pi$ is some path from
$v_*$ to $v$. Since the sum of $x_a$ over any circuit is zero, the value $y_v$ is well-defined (does not depend on the chosen path) and it is obvious that $d\yy=\xx$.

%\sol{prob:basis_size_inv} 
\solL{prob:cycle_basis_linalg}{A}\
Hint: The general case reduces to the case of a strongly connected graph due to %\bref{prob:basis-reduction-to-strong}.
\brefB{prob:relevant_arcs}.
 
%\sol{prob:Euler_formula} 
\solL{prob:cycle_basis_linalg}{B}\
It suffices to prove the formula in the case $r=1$. Then, using notation from \bref{prob:cycle_basis_linalg}, we have
$|\LL|=\dim U=|\AA|-\dim\Ran d=|\AA|-|\VV|+\dim\ker d$. It is easy to see that $\ker d$ is one-dimensional (spanned by the vector 
$\xx=\xone$). %, i.e.\ $x_a=1$ $\forall a\in\AA$.

%\sol{prob:basis-of-cycles-exercise}
\solL{prob:cycle_basis_linalg}{C}\
% Answer: 
The three basis cycles are $(1\to 2\to 4\to 1)$, $(1\to 3\to 4\to 1)$ and $(1\to 2\to 3\to 4\to 1)$.   
(Consistent with Euler's formula, which predicts that $|\LL|=6-4+1=3$.)

%\sol{prob:cycle_basis-counterex}
\solL{prob:cycle_basis_linalg}{D}\
Let $\VV(\Ga)=\{1,2,3\}$, $\AA(\Ga)=\{(1\to 2),(1\to 3), (2\to 3)\}$. Then $\LL=\emptyset$, while $U$ is a one-dimensional subspace spanned by the vector
$\xx$ with components $x_{12}=x_{23}=1$, $x_{13}=-1$.

%%-----------------------
\subsection*{Section~\ref{sec:arcsums}}
\addcontentsline{toc}{subsection}{Section~3}

\sol{prob:mean-min-fungraph}\ 
Due to 
\brefB{prob:delete-unused-arcs},
%\bref{prob:min-for-fungraph}, 
$ \prE f_{\Ga}$ equals the average size of the set $|\VV(\lim\Ga)|$ for a random functional graph with $n$ nodes.
The nodes of $\VV(\lim\Ga)$ are called the \emph{cyclic points}\ of the adjacency map $\phi$. 

Let $\phi:[n]\to[n]$, where $[n]=\{1,2,\dots,n\}$. Let $\mu_\phi:[n]\to\{0,1\}$ be the indicator function of the set of cyclic points of the map $\phi$.
The number of cyclic points of $\phi$ is $|\lim\phi|=\sum_{i=1}^n \mu_\phi(i)$. 
The expectation $\prE \mu_\phi(i)$ over all $n^n$ maps $\phi$ is independent of $i$.    
Therefore $\prE|\lim\phi|=n\prE \mu_\phi(1)=n\,\Prob(1\;\mbox{\rm is cyclic})$.
Now,
$$
 \Prob(1\;\mbox{\rm is cyclic})=
\sum_{k=1}^{n}\frac{1}{n}\prod_{j=0}^{k-1}\frac{n-j}{n},
$$
where the $k$-th summand is the probability that the point $1$ is cyclic and $k$ is its least period. Hence
$$
\ba{rl}
\prE|\lim\phi|&\dst
=\sum_{k=1}^{n} \frac{n!}{(n-k)!\,n^{k-1}}
%=n\sum_{k=1}^n{n-1\choose k-1}\frac{(k-1)!}{n^{k-1}}
%\\[3ex]&\dst 
%=n\sum_{m=0}^{n-1}{n-1\choose m} \,\int_0^\infty e^{-nx} x^m\,dx=
%\\[3ex]&\dst 
=n\int_0^\infty e^{-nx} (x+1)^{n-1}\,dx.
\ea
$$
Since $\ln(e^{-x}(x+1))=-x^2/2+O(x^3)$ near the point of maximum $x=0$, the integral in the last line 
is asymptotic to $\int_0^\infty e^{-nx^2/2}\,dx=\sqrt{\pi/2n}$,
so $\prE|\lim\phi|\sim \sqrt{\pi n/2}$.

\medskip\noindent
{\bf Comment. } 
Combinatorial problems concerning statistics of random maps, were first studied in the 1950s
and a formula describing the distribution of the number of cyclic points has been known since then, see 
\cite{Harris_1960}. It seems however that the asymptotics of $\prE |\lim\phi|$ (although it is but a rather simple consequence) appeared in print somewhat later, see \cite[Th.~3.4]{Moon_1970} and a more general result in \S~3.1 of Kolchin's comprehensive book \cite{Kolchin_1986}. 

An elegant short proof of the asymptotical result, which 
bypasses the explicit distribution formula and employs, without naming it, the formalism of the theory of combinatorial species (and can serve, by the way, as a mini-introduction to it!) is found in \cite{FlOl_1990}(line 2 in Eq.\ (17) and Theorem 2(ii)).

The given proof was shown to me by G.V.\ Kalachev (private communication, 2016).

%\sol{prob:non-anti-automorphic_minimizer}
\solL{prob:isomorph}{C}\
Consider the basis of cycles $\LL=\{\gamma_1,\gamma_2\}$, where $\gamma_1=\{a,a'\}$, $\gamma_2=\{a,b,c\}$.
%The critical point equations are
%$$
% x_{a}=\lambda_1+\lambda_2,\quad
% x_{a'}=\lambda_1,\quad x_b=x_c=\lambda_2,
% \quad x_{a}x_{a'}=x_{a}x_b x_c=1. 
%$$
%Hence 
%$$
% \lambda_2=x_{a1}^{-1/2}, \quad \lambda_1=x_{a1}^{-1},   
%$$
%and we get the equation for $x_{a}$
From the system \eqref{critpteqs} we find that $x_a$ satisfies the equation
$$
 x_{a}=x_{a}^{-1}+x_{a}^{-1/2}.
$$
In particular, $x_a>1$ and $x_{a'}=1/x_a<1$. (Numerically, $x_a\approx 1.4902$ and $f_{\Ga}\approx 3.7996$.)

There is no contradiction with
\brefB{prob:isomorph}: 
%\bref{prob:automorphic_minimizer}: 
the map in (a) is not an anti-automorphism,
since it does not reverse the arcs $b$ and $c$, while in (b)  the anti-automorphism maps the arcs $a$ and $a'$ to themselves.

\sol{prob:ex-quasibasis}
a) The critical point system \eqref{critpteqs} consists of the constraints $x_b x_c x_d=1$, $x_a x_c x_e=1$ and the equations
$w_a x_a=w_e x_e=\lambda_2$, $w_b x_b=w_d x_d=\lambda_1$, $w_c x_c=\lambda_1+\lambda_2$. 
Eliminating $x_c$, $x_d$, $x_e$ and $\lambda_{1,2}$,
%By elimination of
%$$
% x_d=\frac{w_b x_b}{w_d},\quad x_e=\frac{w_a x_a}{w_e}, \quad x_c=\frac{w_a x_a+w_b x_b}{w_c},
%$$
we get the system for $x_a$, $x_b$
$$
w_b x_b^2(w_a x_a+w_b x_b)=w_c w_d,
\qquad
w_a x_a^2(w_a x_a+w_b x_b)=w_c w_e.
$$
It follows that
$$
 \frac{x_b}{x_a}=\sqrt{\frac{ w_a w_d}{w_b w_e}}.
$$
Substituting and solving for $x_a$ we find
$$
x_a=\left(\frac{w_e}{w_a}\right)^{1/2}\left(\frac{w_c}{\sqrt{w_aw_e}+\sqrt{w_bw_d}}\right)^{1/3}.
$$
Then $\lambda_2=x_aw_a$, $\lambda_1=w_b x_b$ and $f_{\Ga,\hat\LL}(\ww,\xone)=3(\lambda_1+\lambda_2)$. The answer is
$$
 f_{\Ga,\hat\LL}(\ww,\xone)=3 w_c^{1/3}\left(\sqrt{w_aw_e}+\sqrt{w_bw_d}\right)^{2/3}.
$$
In the case of a symmetric weight the answer is simplified:
$$
 \hat m(\ww)=f_{\Ga,\hat\LL}(\ww,\xone)=3\left(4 w_a w_d w_c\right)^{1/3},
$$
and the components of the minimizer are%
\begin{NoHyper}%
\footnote{If one wants to calculate the minimizer $\xx$ using \bref{prob:Jacobian_for_minimizer}(c), differentiation must be carried out before 
symmetrization.}:
 $x_a=x_b=(2^{-1} w_a^{-2} w_c w_d)^{1/3}$, $x_d=x_e=(2^{-1} w_a w_c w_d^{-2})^{1/3}$,
$x_c=(4 w_a w_d w_c^{-2})^{1/3}$.  
%[The corresponding value $p_{\gamma_3}=2^{1/4}$.] 
\end{NoHyper}

\smallskip
b) Now the system \eqref{critpteqs} contains, in addition to the constraints $x_b x_c x_d=1$, $x_a x_c x_e=1$, $x_a x_b x_c=t$,
the equations $w_ax_a=\lambda_2+\lambda_3$, $w_bx_b=\lambda_1+\lambda_3$, $w_cx_c=\lambda_1+\lambda_2+\lambda_3$,
$w_dx_d=\lambda_1$, $w_ex_e=\lambda_2$.

The transposition of nodes $1\leftrightarrow 3$ is an anti-automorphism of $\Ga$. In the case of symmetric weight,
it leaves the problem data invariant, hence the minimizer possesses the symmetry $x_a=x_b$, $x_d=x_e$ (hence $\lambda_1=\lambda_2$).
Eliminating variables 
$$
 x_c=\frac{t}{x_a^2}, \quad x_d=\frac{1}{x_a x_c}=\frac{x_a}{t}
$$ 
and expressing $\lambda_1+\lambda_3$ in two ways: as $w_1 x_1$ and as $w_c x_c-w_d x_d$, we obtain an equation for $x_a$
$$
 w_a x_a = \frac{w_c t}{x_a^2}-\frac{w_d x_a}{t}.
$$
It follows that
$$
 x_a=\left(\frac{w_c t^2}{w_a t+w_d}\right)^{1/3}
$$
and we find
$$
m(\ww,t)=3\left(\frac{w_c (w_a t+w_d)^2}{t}\right)^{1/3}.    % \quad % 3 w_c^{1/3} t^{-1/3} \left(w_a t+w_d\right)^{2/3}.
$$

c) Substituting the values $x_a$, $x_b$, $x_c$ of the minimizer found in (a) we get
$$
 p_{\gamma_3}(\xx)=x_a x_b x_c=w_d/w_a
$$
and it is easy to see that $m(\ww,w_d/w_a)$ coinsides with $\hat m(\ww)$.

Differentiating $\log m(\ww,t)$ we get
$$
\frac{d \log m(\ww,t)}{dt}=-\frac{1}{3t}+\frac{2}{3}\,\frac{w_a}{w_a t+w_d}=\frac{w_a t-w_d}{3(w_a t+d)}.
$$
Therefore $\min_t m(\ww,t)$ is indeed attained at $t=w_d/w_a$.

\bigskip\noindent{\bf Remark. } Problem (b) can be solved in a closed form even without symmetry assumption.
The calculation is rather long and we only give the final result for $\ww=\xone$ and arbitrary values of $t_1$, $t_2$, $t_3$:
$$
 f_{\Ga,\LL}(\xone,\xtt)=3\left(\frac{(t_1+t_3)(t_2+t_3)}{t_3}\right)^{1/3}.
$$

\iffalse
$$
\lambda_1=t_1\left(\frac{t_2+t_3}{t_3(t_1+t_3)^2}\right)^{1/3},
\qquad 
\lambda_2=t_2\left(\frac{t_1+t_3}{t_3(t_2+t_3)^2}\right)^{1/3},
$$
%
$$
\lambda_1+\lambda_3=\left(\frac{t_3^2(t_1+t_3)}{(t_2+t_3)^2}\right)^{1/3},
\qquad 
\lambda_2+\lambda_3=\left(\frac{t_3^2(t_2+t_3)}{(t_1+t_3)^2}\right)^{1/3},
$$
%
$$
\lambda_1+\lambda_2+\lambda_3=\left(\frac{(t_1+t_3)(t_2+t_3)}{t_3}\right)^{1/3}.
$$
\fi

\subsection*{Section~\ref{sec:graphs-variants}}
\addcontentsline{toc}{subsection}{Section~4}

\sol{prob:minquotient-mingraph}\
If $\yy$ is an admissible vector for \eqref{minquotient}, then the vector $\xx$ with $x_a=y_{\alpha(a)}/y_{\beta(a)}$ is an admissible vector
for the corresponding homogeneous problem \eqref{mingraph}, so $f_{\Ga}(\ww)\leq f^{\div}_{\Ga}(\ww)$. Conversely, if $\xx>0$ is an admissible vector
for the homogeneous problem \eqref{mingraph}, then there exists a node-weight vector $\yy$ with $x_a=y_{\alpha(a)}/y_{\beta(a)}$ $\forall a\in\AA$.
(That is, the equation $d\log \yy=\log\xx$ has a solution;
the cosntruction is the same as in the solution of 
\bref{prob:cycle_basis_linalg}.)
 
Taking a sequence $\xx^{(i)}>0$ of admissible vectors converging to the minimizer and a corrsponding sequence $\yy^{(i)}$ we
get $f^{\div}_{\Ga}(\ww)\leq\lim S^{\div}(\yy^{(i)}|\ww)=\lim\langle\ww,\xx^{(i)}\rangle=f_{\Ga}(\ww)$.

\medskip\noindent{\bf Comment.} The problem \eqref{minquotient} can be considered as the homogeneous problem \eqref{mingraph} with an additional ``differential'' constraint: $\log\xx\in \Ran d$, where $d:\,\RR^\VV\to\RR^\AA$ is the operator introduced in \bref{prob:cycle_basis_linalg}. That is, $\log\xx$ is required to be a {\em 1-coboundary}. This problem says that this condition is automatically satisfied. The reason, in  essense, is a topological fact --- Poincar\'{e}'s lemma for directed graphs.

%\sol{prob:minquotient-nominimizer}
\solL{prob:minquotient-mingraph}{A}\
Let $\xx^*$ be the minimizer  for the homogeneous problem \eqref{mingraph}. If $\Ga$ is strongly connected,  then all $x^*_a>0$, hence
there exists a corresponding minimizer $\yy$ for problem \eqref{minquotient}. If $\Ga$ is not strongly connected then $x^*_a=0$ precisely
for those arcs $a$ that are not relevant \bref{prob:relevant_arcs}, \bref{prob:delete-unused-arcs}. A corresponding minimizer $\yy$ would than 
have $y_{\beta(a)}>0$ and $y_{\alpha(a)}=0$. Therefore the node $\alpha(a)$ must have indegree zero; this is equivalent to the stated condition.  

%%%%%%%%%%%%%%%%
\sol{prob:Nesbitt}
Given $\yy=(y_1, y_2, y_3)$, define $\xuu=(u_1,u_2,u_3)$ as the solution of the system of 3 linear equations $y_i=u_j+u_k$ with $\{i,j,k\}=\{1,2,3\}$.
Explicitly, $u_i=\frac{1}{2}(y_j+y_k-y_i)$. In the $y$-variables the Nesbitt inequality becomes
$$
 \frac{1}{2}\left(\sum_{i\neq j}\frac{y_i}{y_j}-3\right)\geq \frac{3}{2}.
$$
It has the required form $S_\Ga^\div(\yy)\geq f^\div_\Ga$ if  $\Ga$ is the complete directed graph $\Knoloops{3}$ on $3$ nodes without loops
(previously seen in \bref{prob:count-circuits-cycles-small-n}). Here $f^\div_\Ga=6$ and any minimizer is of the form $\yy=c\xone$, $c>0$.

%\sol{prob:Shapiro}
\solL{prob:Nesbitt}{A}\
Given $\yy=(y_1, \dots, y_n)$, define $\xuu=(u_1,\dots,u_n)$ as the solution of the system of $n$ linear equations $y_i=u_{i+1}+u_{i+2}$ 
(with $u_{n+1}=u_1$, $u_{n+2}=u_2$). For $n$ odd, the matrix of the system is nondegenerate.  
Explicitly, $u_1=\frac{1}{2}\sum_{i=1}^{n} (-1)^{i-1} y_n$ and other values $u_i$ are obtained by a cyclic permutation of indices. 
Written in the $y$-variables, the Shapiro inequality reads
$$
 \frac{1}{2}\left(\sum_{i\neq j}\eps_{ij}\frac{y_i}{y_j}-n\right)\geq C_n,
$$
where $\eps_{ij}=(-1)^{i-j}$ if $i<j$ and $\eps_{ij}=(-1)^{i-j+1}$ if $i>j$.
The corresponding graph $\Ga$ is $\Knoloops{n}$ defined in \bref{prob:count-circuits-cycles-small-n}. The required 
weights are obvious from the above formula.

\medskip\noindent
{\bf Comment.} Nesbitt's inequality \cite{Nesbitt_1903} has been ever popular in students' math competition training. In 1954  H.~Shapiro \cite{Shapiro_1954} 
proposed  a generalization of it and conjectured that $C_n=n/2$ (in which case $\yy=\xone$ would be a minimizer). Shapiro's conjecture 
turned out to be correct only for finitely many $n$: even $n\leq 12$ and odd $n\leq 23$; in fact $\lim_{n\to\infty}(2/n)C_n\approx 0.98913$ as shown by Drinfeld \cite{Drinfeld_1971}.  
The theory developed here does not work in the presence of sign-changing
weights. Whenever $C_n<n/2$ in Shapiro's problem, a minimizer does not exist (although the condition stated in 
\brefA{prob:minquotient-mingraph}
%\bref{prob:minquotient-nominimizer} 
is fulfilled)
and approximate minimizers have certain components arbitrarily small. Compared to our minimization problem
\eqref{inf1} and its instances \eqref{mingraph}, \eqref{minquotient}, Shapiro's problem requires totally different approach and is overall much harder.
Further information and references can be found in the surveys \cite{Clausing_1992} and \cite{Fink_1998}.   

Interestingly, the reduction of Shapiro's problem to the variables $y$ as described in this solution is not even new: it is found in \cite{Troesch_1985}!

%P.H.~Diananda suggested an yet more general cyclic inequality (where, as before, one should identify $u_{n+j}$ with $u_j$),
%\beq{Diananda}
% \sum_{i=1}^n \frac{u_i}{u_{i+1}+\dots+u_{i+k}}\geq C_{n,k}.
%\eeq
%In \cite{Sadov_2016} we give details and show that $\inf_{2\leq k\leq n}C_{n,k}$ lies between $\ln 2$ and a constant $\approx
%0.9305$. The inequality \eqref{Diananda}, again, can be formally written as $S_\Ga^\div(\yy|\ww)\geq f^\div_\Ga(\ww)$ provided $k$ does not divide $n$. 
%Here, again, $\Ga=\Knoloops{n}$ and $\ww$ is a sign-changing weight, which can be explititly found as in the above solution. 

Cyclic sums similar to Shapiro's with $k>2$ terms in the denominators were studied by P.H.~Diananda. A result similar to Drinfeld's (although without the claim of optimality of the constant) is obtained in \cite{Sadov_2016}. Diananda's sums can, too, be formally written as $S_\Ga^\div(\yy|\ww)\geq f^\div_\Ga(\ww)$ provided $k$ does not divide $n$. 
Here, again, $\Ga=\Knoloops{n}$ and $\ww$ is a certain sign-changing weight.

%%%%%%%%%%%%%%%%
% Harmonic sums
\subsubsection*{Subsection~\ref{ssec:harmonic}}
\addcontentsline{toc}{subsubsection}{Subsection~4.2}

\sol{prob:harmonic-graph}% 
Put $x_v=e^{h(v)}$ ($v\in\VV$). The matrix of constraints
contains $r=|\VV_1|+1$ rows. The entries in the first $|\VV_1|$ rows are
$A_{vv'}$, $v\in\VV_1$, $v'\in\VV$, of which the following are nonzero:
$A_{vv}=d^+(v)$, $A_{vv'}=-1$ ($v\in\VV^+(v)$).

One more row (we will use `$0$' as the corresponding row index) describes the boundary constraint:
$A_{0v}=1$ for $v\in\VV_0$; $A_{0v}=0$ for $v\in\VV_1$.

We obtain: $F^\Delta(\tau)=f(\xone,\xtt,A)$, cf.\ Eq.\ \eqref{inf1}, where $\xtt=[1,\dots,1,t]^T$, $t=e^\tau$.

The vector $\xx=e^{\tau/|\VV_0|}\xone$, corresponding to the constant harmonic function $h(v)\equiv \tau/|\VV_0|$, is admissible, whence the upper estimate \eqref{harmonic-min} for the minimum.

\sol{prob:harmonic-graph-lb}% 
(a) Let us prove that the vector $\xx=e^{\tau/|\VV_0|}\xone$ is not just admissible (as pointed out in the previous solution) but is the minimizer. 
%By uniqueness, it is the minimizer and $f[\phi]=me^{t/n}$.

According to the equations of extremizer \eqref{critpteqs}, we have to show that there exists a vector of Lagrange's multipliers $\Lambda=[ (\lambda_v)_{v\in\VV_1}; \lambda_0]$ such that
$A^T\Lambda=\xone$, or explicitly: 
$$
\ba{ll}
 d(v)\lambda_{v}-\sum_{v'\in\VV(v)\cap\VV_1}\lambda_{v'}=1
\qquad &(v\in\VV_1),\\[2ex]
 \lambda_0+\sum_{v'\in\VV(v)\cap\VV_1} \lambda_{v'}=1
\qquad &(v\in\VV_0).
\ea
$$

Put $\lambda_v=\lambda_1$ for all $v\in\VV_0$
By the regularity assumption, every equation from the first group becomes 
$
 q\lambda_1=1
$
and every equation from the second group becomes $\lambda_0+p\lambda_1=1$. Just take
%The values sought for are 
$\lambda_1=1/q$ and $\lambda_0=1-p/q$.

\iffalse
Clearly, the equations of constraints are satisfied by this solution. It remains to show that there exists a set of Lagrange multipliers to satisfy the linear (first)
part of the system \eqref{critpteqs}, which is to say that
the vector $\xx=\xone$ lies in $\Ran A^T$. We will prove the equivalent fact: that $\xx=\xone$ is orthogonal to $\Ker A$. In the matrix $(A_{ji})$, the index $j$ (counting the constraints) takes the values $v\in\VV_1$ corresponding to the harmonic equation at the inner nodes and $0$ (just to distinguish) corresponding to the
condition imposed on the boundary values.

The index $i$ (counting the variables) takes the values $v\in \VV$.

The nonzero entries of the matrix $A$ are as follows.
 
If $j=v\in\VV_1$, then $A_{vv}=d^+(v)$, $A_{vv'}=-1$ for $v'\in\VV^+(v)$.

If $j=0$, then $A_{0v}=1$ for $v\in \VV_0$.

The condition $\xuu\in\Ker A$ has the form
%$$
%\ba{l}
%  \sum_{v'\in\VV} A_{vv'}u_{v'}=0, \qquad v\in\VV_1;
%\\
%  \sum_{v\in\VV_0} A_{0v}u_{v}=0,
%\ea
%$$
%that is,
$$
\ba{l}
  \sum_{v'\in\VV^+(v)} (u_v-u_{v'})=0, \qquad v\in\VV_1;
\\
  \sum_{v\in\VV_0} u_{v}=0.
\ea
$$

\fi 
%%%%%%%%%%%%%%%%%%%%%%%%%%%%%%%%%%%%%%%%%%

\medskip
(b) To avoid calculations, we refer to the symmetry argument \bref{prob:symmetry}. The problem is invariant under the involution on the index set $i\leftrightarrow n-i$. Therefore the minimizer must have $x_0=x_n$. A harmonic function constant on the boundary is constant.
Hence $x_i=x_0$ for all $i$. The result follows.

%%%%%%%%%%%%%%%%%%%%%%%%%%%%%%%%%%%%%%%%%%
\iffalse
Explicit calculations:
$$
 \mat{cccc}{1 &-1 & 0 & 0 \\
            0 &2 &-1 & 0\\
            0 &-1 &2 &-1\\
            0 &0 &-1 & 2\\
            1 &0 &0  &-1 }
$$

$$
 2x-y=1,
\quad
 -x+2y-z=1,
\quad
 -y+2z=1
$$

$$
 z+x=3, \quad y=2x-1=2z-1
$$

$$
 x=z=3/2, y=2.
$$

$$
 \frac32 e_1+2e_2+\frac32 e_3=\mat{c}{-3/2\\ 1\\1\\1\\-3/2}
$$

$$
\frac52 e_0+ \frac32 e_1+2e_2+\frac32 e_3=\xone_5
$$

If $v'\in\VV_1$, then
$$
 \sum_{v\in\VV_1} A_{vv'}+A_{0v'}=A_{v'v'}+\sum_{v\in\VV^-(v')\cap\VV_1} A_{vv'}+0=??
$$
\fi
%%%%%%%%%%%%%%%%%%%%%%%%%%%%%%%%

\sol{prob:harmonic-example} The matrix $A$ is
$$
 A=\mat{rrrr}{3 & -1 & -1 & -1\\
              -1 & 2 & -1 & 0\\
               0 & 0 & 1 & 1}.
$$
The first group of the necessary conditions for extremum \eqref{critpteqs}  has the form
$\xx\in\Ran A^T$, that is, $\xx\bot\ker A$, which is the linear equation
$x_1+3x_2+5x_3-5x_4=0$.

Including the equations of constraints $x_1^3=x_2x_3x_4$,
$x_2^2=x_1x_3$, $x_3x_4=t=e^\tau$, we get a system of 4 equations with 4 unknowns, which must have a unique positive solution according to the theory.

Solving the equations of constraints for 
$x_2,x_3,x_4$, we get
$$
 \log x_2=3\log x_1-\tau,
\quad\log x_3=5\log x_1-2\tau,
\quad \log x_4=-5\log x_1+3\tau.
$$

The linear relation between $x_{1\div 4}$ yields the algebraic equation for $x_1$ %(with $t=e^\tau$)
$$
 x_1+3t^{-1} x_1^3+5t^{-2} x_1^5-5t^3 x_1^{-5}=0.
$$
Substituting $x_1=\xi t^{1/2}$, we obtain the equation for $\xi$ free of $t$: 
$5\xi^{10}+3\xi^8+\xi^6-5=0$, which can be solved numerically ($\xi=0.9355637\dots$).

The final result is
$$
 F^\Delta(\tau)=Ce^{\tau/2},
$$
where $C=F^\Delta(0)\approx 3.86638136$.
(The ratio $C/4\approx 0.96659534$.)

\sol{prob:harmonic-const-ratio}
Let $|\VV|=n$, $|\VV_1|=m$, so $|\VV_0|=n-m$.
The scheme of calculation follows the sample case \bref{prob:harmonic-example}. 

The null space of the matrix $A$ is one-dimensional; let it be spanned by a vector $\xaa=[a_1,\dots,a_n]$. 

Let $\xee$ be the $(m+1)$-dimensional vector $\xee=[0,\dots,0,1]^T$.
The critical point system \eqref{critpteqs}
reduced to the $x$-variables has the form
$$
 \ba{l}
   A\log\xx=\tau \xee,\\
   (\xx,\xaa)=0.
\ea
$$

Since the constant function $h\equiv 1$ is harmonic, we have
$$
 A\xone=(n-m)\xee.
$$
Therefore a solution of the equation
$A\log\xx=\tau\xee$ has the form
$\log\xx=p\tau\xone+y\xaa$ with $p=(n-m)^{-1}$ and
the coefficient $y$ to be determined.

Writing componentwise, we get
$$
 x_i=e^{p\tau+ya_i}.
$$
Denote
$$
 F(\tau,y)=\sum_{i=1}^n e^{p\tau+ya_i}.
$$
and let $G(\tau,y)=\partial_y F(\tau,y)$.
The condition $(\xx,\xaa)=0$ takes the form
$G(\tau,y)=0$, which is an equation to be solved for $y$. Let $y_*=y_*(\tau)$ be the (unique) positive root and $\xx^*$ be the corresponding vector. 
Note that 
$$F^\Delta(\tau)=\sum x_i^*=F(\tau,y_*(\tau)).$$
   
We have
$$
 \frac{d}{d\tau}\left(F^\Delta(\tau)e^{-p\tau}\right)=
\frac{\partial}{\partial\tau}\left(e^{-p\tau} F(\tau,y_*)\right)
+e^{-p\tau}G(\tau,y_*)\frac{\partial y_*}{\partial\tau}.
$$
Both summands in the right-hand side equal zero and we obtain the result as claimed.

%%%%%%%%%%%%%%%%
\subsection*{Section~\ref{sec:trees}}
\addcontentsline{toc}{subsection}{Section~5}
% Tree-sums

% WRONG!
%\sol{prob:tree-unitproduct}
%If $\xtt=\xone$, then the system of constraints can be written as $\prod_{\pi(\lambda',\lambda'')}=1$ for all pairs of leaves $\lambda'\neq\lambda''$, where $\pi(\lambda',\lambda'')$ is the path connecting those leaves. (In such a form the number of constraints is generally excessive but one can select a minimal independent subset.)

\sol{prob:treemin-ineq}
If, for some $\lambda_i$, it occurs that $P_{\pi(\lambda_i)}(\yy)=t'_{\lambda_i}>t_{\lambda_i}$, then $y_{\lambda_i}$ can be replaced by the smaller value $y_{\lambda_i} t_{\lambda_i}/t'_{\lambda_i}$; other maximal chain products are not affected.

\sol{prob:tree-quotient}
We construct $F$  as a forest whose nodes correspond one-to-one to the summands in \eqref{treesum-quotient}.
Nodes of $F$ will have labels of two kinds: (i) nodes labeled $\hat v$, where $v$ is a node (possibly a leaf) of $T$, $v\neq\rho$;
(ii) nodes labeled $\lambda'$, where $\lambda\in \Lambda(T)$.  
The nodes of the first, resp., second kind will correspond to the terms of the  first, resp., second sum in \eqref{treesum-quotient}.
Every non-leaf node of $T$ except the root gives rise to a unique node of $F$, while every leaf of $T$ gives rise to
exactly two nodes of $F$, labeled $\hat\lambda$ and $\lambda'$.  
There is no node in $F$ corresponding to $\rho$.

The roots of $F$ are those nodes $\hat v$, for which
$\alpha(v)=\rho$ in $T$.
Thus the set of roots of $F$ is in one-to-one correspondence with the set $\alpha^{-1}(\rho)$.

The parent function $\alpha_F$ for $F$ is defined by the rules:
\\[0.5ex]\hspace*{1em}
(i) $\alpha_F(\hat v)=(\alpha(v))\hat{}$  (note that $\alpha(v)\neq \rho$).
\\[0.5ex]\hspace*{1em} 
(ii) $\alpha_F(\hat\lambda)=\lambda'$ for $\lambda\in\Lambda(T)$. 
The nodes labeled $\lambda'$ are the leaves of $F$.

\smallskip
The vector of node values $\hat\yy$ for $F$ is defined as follows: 
\\[0.5ex]\hspace*{1em}
(i) $\hat y_{v}=y_{\alpha(v)}/y_{v}$ (including the case $v\in\Lambda$);
\\[0.5ex]\hspace*{1em} 
(ii) $\hat y_{\lambda'}=y_{\lambda}/x_{\rho}$.

\smallskip
By construction we have $Q_T(\xtt|\yy)=\sum_{v\in\VV(F)}{\hat y_v}$. Also, the maximal chain products in $F$ are 
$P_{\pi(\lambda')}(\hat\yy)=t_\lambda$. The identity $\inf_{\yy>0} Q_T(\xtt|\yy)=m_F(\xtt)$ follows at once.

%%%%%

\sol{prob:tree-to-general}
Let $\sim$  be the equivalence relation on $\VV(T)$ that collapses all leafs into a single node $\lambda_*$. Define the graph $\Ga_0=T/\sim$. 
The graph $\Ga$ is obtained by adding to $\Ga_0$ the arc $(\lambda_*\to \rho)$. Thus $\VV(\Ga)=\VV(\Ga_0)=(\VV(T)\setminus\Lambda(T))\cup\{\lambda_*\}$.

There is a one-to-one correspondence $e_*:\VV(T)\to \AA(\Ga)$ between nodes of $T$ and arcs of $\Ga$. 

(i) If $v\neq\rho$ and $v\notin\Lambda(T)$, then $e_*(v)=e(v)$, the unique arc ending in $v$.

(ii) If $\lambda\in\Lambda(T)$, then $e_*(\lambda)=(\alpha(\lambda)\to\lambda_*)$.

(iii) $e_*(\rho)=(\lambda_*\to\rho)$. (This arc does not correspond to any arc in $T$.)  

\smallskip
The graph $\Ga$ is strongly connected, hence by Euler's formula 
\brefB{prob:cycle_basis_linalg},
%\bref{prob:Euler_formula}, 
the cardinality of any basic of circuits for $\Ga$
is $|\AA(\Ga)|-|\VV(\Ga)|+1=|\VV(T)|-(|\VV(T)|-|\Lambda(T)|+1)+1=|\Lambda(T)|$. 
The $|\Lambda(T)|$ cycles $\gamma(\lambda)=e_*(\pi(\lambda))\cup e_*(\rho)$
are obviously independent,
hence the system $\LL$ formed by them is a basis of cycles in $\Ga$.

Let $\xx$ be a vector labeled by the arcs of $\Ga$, 
$\yy$ be a vector labeled by the nodes of $T$, and
$x_{e_*(v)}=y_v$. Clearly, $p_{\gamma(\lambda)}(\xx)=P_{\pi(\lambda)}(\yy)$, where the left-hand side pertains to $\Ga$ and the right-hand side to $T$,
so the constraints $P_{\pi(\lambda)}(\yy)=t_\lambda$ for problem \eqref{mintree} translate into  the constraints $p_{\gamma(\lambda)}(\xx)=t_\lambda$ 
for problem \eqref{mingraph}.
Also $\sum_{v\in\VV(T)} y_v=\sum_{v\in\VV(T)} x_{e_*(v)}=\sum_{a\in\AA(\Ga)} x_a=\langle\xx,\xone\rangle$. 
The formula $m_T(\xtt)=f_{\Ga}(\xone,\xtt)$ follows.

(We have silently assumed that $\rho\notin\Lambda(T)$. In the excepted case, $T$
is the one-node tree. The corresponding graph $\Ga$ is then  one node, one loop graph.)

\sol{prob:tree-minimizer}\
Referring to the reduction \bref{prob:tree-to-general}
(or arguing in some different way), one can recast the %constrained 
minimization problem \eqref{mintree}
in the unconstrained form
$$
 \sum_{v\in\VV\setminus\Lambda} y_v+\sum_{\lambda\in\Lambda} \frac{t_\lambda}{P_{[\rho,\alpha(\lambda)]}}\,\to\,\min,
$$
which, by noticing that $(\partial/\partial y_v) \log P_{[\rho,\alpha(\lambda)]}=1/y_v$ if $\lambda\in\Lambda_v$ and $0$ otherwise,
immediately leads to the necessary conditions of extremum \eqref{tree-minimizer}.%
\begin{NoHyper}%
\footnote{There is no need in the Lagrange multipliers and no notational collision occurs in the context of trees,
where we use the symbol $\lambda$ to denote leaves of a tree.}
\end{NoHyper}

\sol{prob:tree-recurrrence}
Put $y_\rho=1/k$. Then the constraints in \eqref{mintree} take the form $P_{[v,\lambda]}=kt_\lambda$ with $v\in\alpha^{-1}(\rho)$, $\lambda\in\Lambda_v$.
The minimization problems for the trees $T_v$, $v\in\alpha^{-1}(\rho)$, are independent, hence the recurrence relations.

(This result is similar to \bref{prob:dynprog}.)

%%%
\sol{prob:const-height-tree-radicals}
(a) If $\yy(c)$ is the minimizer corresponding to the product vector $c\xtt$ ($c>0$), then obviously $\yy(c_1)=(c_1/c_2)^{1/\ell}\yy(c_2)$;
hence $m_T(\xtt)$ is homogeneous of order $1/\ell$.

\smallskip
(b) We have $\min_{k>0} (k^{-1}+ak^{1/h})=\ell(\frac{a}{h})^{h/\ell}$ and $\min$ is attained at $k^*=(\frac{a}{h})^{-h/\ell}$.
Put $a=\sum_{v\in\alpha^{-1}(\rho)} m_{T_v}(\xtt_v)$. The expressibility of $m_T(\xtt)$
in radicals follows by induction: the recurrence \eqref{rec-trees} provides the induction step.
 
\smallskip
(c) Continuing the above argument, note that $k^*$ is expressible in radicals,
and the induction step for the components of the minimizer holds due to the recurrence relations \eqref{rec-trees-minimizer}. 

\smallskip\noindent{\bf Remark.}
Alternatively, (c) can be done by using just the recurrence \eqref{rec-trees} for $m_T$ and referring to \bref{prob:Jacobian_for_minimizer}(c).
For this approach to work we have to generalize problem \eqref{mintree} in an obvious way, introducing node weights $w_v$;
this does not cause any change in formulations and proofs of (a) and (b).   

%%%%% Monotonicity

\sol{prob:dilatation-monotone-trees}\ 
Proof by induction on the height of $T$. The induction base %$h=0$ 
--- the root-only tree --- is trivial.
Now suppose $\VV(T)\neq\{\rho\}$, consider the trees $T_v$, $v\in\alpha^{-1}(\rho)$, re-label them using numerical superscripts as $T^{(i)}$,
$1\leq i\leq |\alpha^{-1}(\rho)|$.
and denote the corresponding product vector -- to -- minimizer maps by $M_i$ (as an abbreviation of $M_{T^{(i)}}$).  
By the induction hypothesis, the maps $M_i$ are dilatation-monotone.

Let $\mu(\xtt)=y^{(T)}_\rho$ be the component of the vector
function $M_T(\cdot)$ corresponding to the root. That is, $\mu(\xtt)$ is the value of the root of the minimizer.  Let $\mu_i(\cdot)$ are the similarly defined functions for the trees $T^{(i)}$.

The critical point equations \eqref{tree-minimizer} imply
$$
y^{(T)}_\rho=\sum_{\lambda\in\Lambda(T)} y^{(T)}_\lambda=\sum_{v\in\alpha^{-1}(\rho)}\sum_{\lambda\in\Lambda(T_v)} y^{(T)}_\lambda =\sum_{v\in\alpha^{-1}(\rho)} y^{(T)}_v.
$$
Hence
$$
 \mu(\xtt)=\sum_i \mu_i\left(\frac{\xtt}{\mu(\xtt)}\right).
$$  

Fix the vector $\xtt$ and consider the scalar function $g(s)=\sum_{i} \mu_i(s\xtt)$. It is monotone (increasing). % in $s$. 
Given $r>0$, the positive real number $\xi=\mu(r\xtt)$ satisfies the equation $\xi=g(r/\xi)$. It follows that $\xi$ as a function of $r$
is increasing and $r/\xi$ is also increasing. Now, the root component $y^{(T)}_\rho$ of the minimizer $M_T(r\xtt)$ is $\xi$, while other components are found among components of the minimizers $M_i((r/\xi)\xtt^{(i)})$; all are increasing in $r$.

\noindent
\begin{minipage}[c]{0.7\textwidth} 
%\sol{prob:no-strong-monotone-trees}
\solL{prob:dilatation-monotone-trees}{A}\
Consider the tree $T$ of constant height 3 with 4 nodes and 2 leaves as shown. The maximal chains are
$[0,2]=\{0,1,2\}$, $[0,3]=\{0,1,3\}$. 
%(In the nomenclature of Table~\ref{tab:trees}),  $T=[[,]]$.) 
Let the maximal chain products be
$P_{[0,2]}=a$, $P_{[0,3]}=b$.  It is easy to find the minimizer explicitly; the values of the leaves are
\end{minipage}
\hspace{2em}\raisebox{-3mm}{
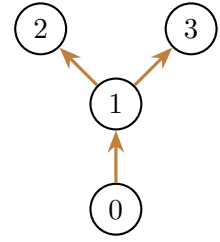
\begin{minipage}[c]{0.2\textwidth} 
%\vspace*{1ex}
\begin{tikzpicture}
\begin{scope}[every node/.style={circle,thick,draw}]
    \node (R) at (1,0) {0};
    \node (A) at (1,1.4) {1};
    \node (B) at (0,2.4) {2};
    \node (C) at (2,2.4) {3};
\end{scope}

\begin{scope}[>={Stealth[brown]},
              every node/.style={fill=white,circle},
              every edge/.style={draw=brown,very thick}]
    \path [->] (R) edge (A);
    \path [->] (A) edge (B);
    \path [->] (A) edge (C);
\end{scope}\end{tikzpicture}
\end{minipage}
}% end \raisebox
$$
  y_2=\frac{\sqrt{a}}{\left(\sqrt{a}+\sqrt{b}\right)^{1/3}},
  \qquad
  y_3= 
  \frac{\sqrt{b}}{\left(\sqrt{a}+\sqrt{b}\right)^{1/3}}.
$$
Clearly, the map $(a,b)\mapsto(y_2,y_3)$ is not strongly monotone.

%\sol{prob:no-dilatation-monotone-graph}
\solL{prob:dilatation-monotone-trees}{B}\
It seems that some effort is needed to produce a counterexample of low complexity. 
Here is one.   

Let $\Ga$ be a graph with three distinguished nodes $\xi,\eta,\zeta$ and the paths between them: 
$$
  \begin{array}{l}
  \pi_1:\;\; \xi\to\eta,\quad\mbox{\rm with $p$ arcs},\\
  \pi'_1:\;\; \eta\to\xi,\quad\mbox{\rm with $p'$ arcs},\\  
  \pi_2:\;\; \eta\to\zeta,\quad\mbox{\rm with $q$ arcs},\\ 
  \pi'_2:\;\; \zeta\to\eta,\quad\mbox{\rm with $q'$ arcs},\\ 
  \pi_3:\;\; \zeta\to\xi,\quad\mbox{\rm with $r$ arcs}.
  \end{array}
$$
The parameters $p,p',q,q',r$ will be chosen later.
The basis of cycles $\LL$ consists of the cycles $L_1=\pi_1\cup\pi'_1$,  $L_2=\pi_2\cup\pi'_2$, $L_3=\pi_1\cup\pi_2\cup\pi_3$.

By \eqref{critpteqs} we express the components of the minimizer as linear combinations of Lagrange's multipliers:
$$
  \begin{array}{l}
  a\in\pi_1\;\Rightarrow \; x_a=\lambda_1+\lambda_3,\\
  a\in\pi'_1\;\Rightarrow \; x_a=\lambda_1,\\  
  a\in\pi_2\;\Rightarrow \; x_a=\lambda_2+\lambda_3,\\ 
  a\in\pi'_2\;\Rightarrow \; x_a=\lambda_2,\\ 
  a\in\pi_3\;\Rightarrow \; x_a=\lambda_3.
  \end{array}
$$
The product constraints yield the system of 3 equations for $\lambda_1,\lambda_2,\lambda_3$,
$$
\begin{array}{l}
 (\lambda_1+\lambda_3)^{p} \lambda_1^{p'}=t_1,
 \\[1ex]
 (\lambda_2+\lambda_3)^{q} \lambda_2^{q'}=t_2,
 \\[1ex]
(\lambda_1+\lambda_3)^{p}(\lambda_2+\lambda_3)^{q} \lambda_3^{r} =t_3.
\end{array}
$$
Let $\lambda_i(s)$ be the values corresponding to the product vector $s\xtt$, where $s$ is the dilatation parameter. 
The result of logarithmic differentiation at $s=1$ in the matrix form reads
$$
\mat{ccc}{a+v_1 & 0 & a\\[0.5ex] 0 & b+v_2 & b\\[0.5ex] a & b & a+b+v_3}\,\mat{c}{d\lambda_1/ds\\[0.5ex]
d\lambda_2/ds\\[0.5ex]
d\lambda_3/ds}
=\mat{c}{1\\[0.5ex]1\\[0.5ex]1},
$$
where
$$
a=\frac{p}{\lambda_1+\lambda_3}, \quad b=\frac{q}{\lambda_2+\lambda_3}, \quad
v_1=\frac{p'}{\lambda_1}, \quad
v_2=\frac{q'}{\lambda_2}, \quad
v_3=\frac{r}{\lambda_3}.  
$$
It follows, in particular, that for $a\in\pi_3$
$$
 \frac{dx_a}{ds}=\frac{d\lambda_3}{ds}
=\frac{v_1 v_2-ab}{D},
$$
where 
$$
 D=ab(v_1+v_2+v_3)+a v_2(v_1+v_3)+bv_1(v_2+v_3)+v_1 v_2 v_3
$$
is the determinant of the matrix.

The parameters $\lambda_i$, $i=1,2,3$, can be assigned any positive values and the parameters $p$, $q$, $p'$, $q'$, $r$ can be assigned any positive integer values.
This can be done so as to make $d\lambda_3/ds<0$.

%. Considering their values as fixed as well as the values of $p'$ and $q'$, and making
%$p$, $q$ sufficiently large, we get $v_1 v_2-ab<0$.
%If now $r$ is large enough so as to make 
%$D>0$, we get $d\lambda_3/ds<0$.

\smallskip
For a concrete numerical example, put $\lambda_1=\lambda_2=\lambda_3=1$ and
take $p=p'=q'=1$, $q=5$, and $r=2$. 
The corresponding cycle product  vector is
$\xtt=(2,32,64)$, and  
$\left.d\lambda_3/ds\right|_{s=1}=-1/64$.

\medskip\noindent{\bf Remark.} A scrupulous reader can verify formula \eqref{dxdt} (Problem \bref{prob:Jacobian_for_minimizer}a) in this example.

%%%%%%%%%%%%%%%
\sol{prob:tree-linearization}\
The map $M$ can be written as a composite map 
$(y'_\lambda)_{\lambda\in\Lambda}\mapsto(y_v)_{v\in\VV\setminus\Lambda} \mapsto(y_\lambda)_{\lambda\in\Lambda}$ 
given by the right-hand sides of the critical point equations \eqref{tree-minimizer}.
Linearizing at $\yy^*$ we get
$$
\begin{array}{l}
\displaystyle
\delta y_v=\sum_{\lambda\in\Lambda_v} \delta y'_\lambda,\quad v\in\VV\setminus\Lambda,
\\[2ex]
\displaystyle
\frac{\delta y_\lambda}{y^*_\lambda}=-\sum_{v\in[\rho,\alpha^{-1}(\lambda)]}\frac{\delta y_v}{y^*_v}, \quad \lambda\in\Lambda.
\end{array}
$$
Eliminating the variables $\delta y_v$ we obtain the vector equation $\delta \hat\yy=B\delta \hat \yy'$ with matrix $B$ as stated in problem.

%\sol{prob:tree-linearization-spectrum}
\solL{prob:tree-linearization}{A}\
Put $\xi_\lambda=\delta y_\lambda/\sqrt{y^*_\lambda}$. Rewriting the linearized iteration in variables $\xi_\lambda$, we see that
the matrix $B$ is similar to the symmetric matrix $\tilde B=-\tilde Q \tilde Q^T$,
where
$$
 \tilde Q_{\lambda,v}=Q_{\lambda,v}\sqrt{\frac{y^*_\lambda}{y^*_v}}.
$$ 
Clearly, the matrix $-\tilde B$ is positive-definite and has nonnegative entries. Its largest eigenvalue equals its norm as an operator in the 
space $\RR^{|\Lambda|}$ with standard Euclidean norm.
%
\iffalse
For any set of vectors $\{f_i\}$ the norm of their Gram matrix $(G_{ij})=(\langle f_i,f_j\rangle)$
is bounded by $\sum\|f_i\|^2$. Applying this to our case we get 
$$
\|B\|\leq \sum_{\lambda\in\Lambda} \sum_{v\in\VV\setminus\Lambda} Q^2_{\lambda,v}\frac{y^*_\lambda}{y^*_v}
$$
Note: $Q^2_{\lambda,v}=Q_{\lambda,v}$. %, so squaring does not change it)
For a fixed $v$, by \eqref{tree-minimizer}, $\sum_\lambda Q_{\lambda,v} y^*_\lambda=y^*_v$. Hence
$$
\|B\|\leq \sum_{v\in\VV\setminus\Lambda} 1 =|\VV|-|\Lambda|.
$$ 
\fi
%
To estimate $\|\tilde B\|$ let us take a vector $\xuu=(u_\lambda)$, $\lambda\in\Lambda$ such that $\sum u^2_\lambda=1$ and estimate
$|\langle \tilde B\xuu,\xuu\rangle|=\|\tilde Q\xuu\|^2$. We have
$$
(\tilde Q\xuu)_v %=\sum_{\lambda} Q_{\lambda,v}\sqrt{\frac{y^*_\lambda}{y^*_v}}u_\lambda 
=\sum_{\lambda\in\Lambda_v} \sqrt{\frac{y^*_\lambda}{y^*_v}}u_\lambda.
$$
By the Cauchy-Schwarz inequality, 
$$
 |(\tilde Q\xuu)_v|^2\leq \sum_{\lambda\in\Lambda_v} \frac{y^*_\lambda}{y^*_v} \, \sum_{\lambda\in\Lambda_v} u_\lambda^2
 =\sum_{\lambda\in\Lambda_v} u_\lambda^2. 
$$
(The last equality holds due to \eqref{tree-minimizer}.)
Therefore
$$
\|\tilde Q\xuu\|^2\leq \sum_{v\in\VV\setminus\Lambda} \sum_{\lambda\in\Lambda_v} u_\lambda^2
=\sum_{\lambda} c_\lambda u_\lambda^2, 
$$
where
$$
 c_\lambda=\sum_{v\in[\rho,\alpha^{-1}(\lambda)]}1=|\pi(\lambda)|-1\leq h,
$$
and the inequality $\|B\|=\|\tilde B\|\leq h$ follows.

%\sol{prob:tree-numerical-iteration}
\solL{prob:tree-linearization}{B}\
It suffices to show that the spectrum of linearization of the iteration operator at $\hat\yy^*$ lies in $(-1,1)$.
The linearization is $B_\tau=\tau B+(1-\tau)I$ and its spectrum, according to %\bref{prob:tree-linearization-spectrum}, 
\brefA{prob:tree-linearization},
lies
in $(-\tau h+1-\tau, 1-\tau)$. The condition stated in the problem is equivalent to $-\tau h+1-\tau>-1$. 

\sol{prob:trees-analytical-numerical}
In Table~\ref{tab:treevalues}, analytical and numerical answers are presented for all 17 trees with at most 5 nodes and the smallest tree $T_{6,*}$ with three distinct leave values of the minimizer.

The trees, with root on top, are depicted below.

\bigskip
\noindent
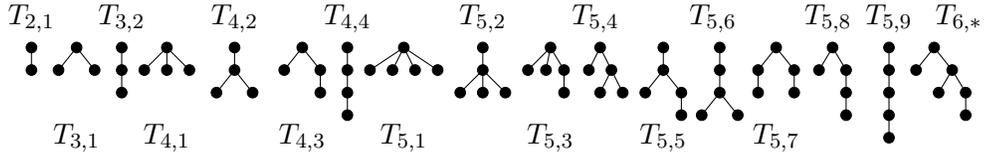
\hspace*{-0.2em}
\begin{tikzpicture}[x=0.3cm,y=0.3cm]
%--T2,1---------------------------------------
\draw [line width=0.4pt] (0,10)-- (0,9);
\draw [fill=black] (0,10) circle (2pt);
\draw [fill=black] (0,9) circle (2pt);
%--T3,1---------------------------------------
\draw [line width=0.4pt] (2,10)-- (1.2,9);
\draw [line width=0.4pt] (2,10)-- (2.8,9);
\draw [fill=black] (2,10) circle (2pt);
\draw [fill=black] (1.2,9) circle (2pt);
\draw [fill=black] (2.8,9) circle (2pt);
%--T3,2---------------------------------------
\draw [line width=0.4pt] (4,10)-- (4,9);
\draw [line width=0.4pt] (4,9)-- (4,8);
\draw [fill=black] (4,10) circle (2pt);
\draw [fill=black] (4,9) circle (2pt);
\draw [fill=black] (4,8) circle (2pt);
%--T4,1---------------------------------------
\draw [line width=0.4pt] (6,10)-- (5,9);
\draw [line width=0.4pt] (6,10)-- (6,9);
\draw [line width=0.4pt] (6,10)-- (7,9);
\draw [fill=black] (6,10) circle (2pt);
\draw [fill=black] (5,9) circle (2pt);
\draw [fill=black] (6,9) circle (2pt);
\draw [fill=black] (7,9) circle (2pt);
%--T4,2---------------------------------------
\draw [line width=0.4pt] (9,10)-- (9,9);
\draw [line width=0.4pt] (9,9)-- (8.2,8);
\draw [line width=0.4pt] (9,9)-- (9.8,8);
\draw [fill=black] (9,10) circle (2pt);
\draw [fill=black] (9,9) circle (2pt);
\draw [fill=black] (8.2,8) circle (2pt);
\draw [fill=black] (9.8,8) circle (2pt);
%--T4,3---------------------------------------
\draw [line width=0.4pt] (12,10)-- (11.2,9);
\draw [line width=0.4pt] (12,10)-- (12.8,9);
\draw [line width=0.4pt] (12.8,9)-- (12.8,8);
\draw [fill=black] (12,10) circle (2pt);
\draw [fill=black] (11.2,9) circle (2pt);
\draw [fill=black] (12.8,9) circle (2pt);
\draw [fill=black] (12.8,8) circle (2pt);
%--T4,4---------------------------------------
\draw [line width=0.4pt] (14,10)-- (14,9);
\draw [line width=0.4pt] (14,9)-- (14,8);
\draw [line width=0.4pt] (14,8)-- (14,7);
\draw [fill=black] (14,10) circle (2pt);
\draw [fill=black] (14,9) circle (2pt);
\draw [fill=black] (14,8) circle (2pt);
\draw [fill=black] (14,7) circle (2pt);
%--------------------------------------------
%--T5,1---------------------------------------
\draw [line width=0.4pt] (16.5,10)-- (15,9);
\draw [line width=0.4pt] (16.5,10)-- (16,9);
\draw [line width=0.4pt] (16.5,10)-- (17,9);
\draw [line width=0.4pt] (16.5,10)-- (18,9);
\draw [fill=black] (16.5,10) circle (2pt);
\draw [fill=black] (15,9) circle (2pt);
\draw [fill=black] (16,9) circle (2pt);
\draw [fill=black] (17,9) circle (2pt);
\draw [fill=black] (18,9) circle (2pt);
%--T5,2---------------------------------------
\draw [line width=0.4pt] (20,10)-- (20,9);
\draw [line width=0.4pt] (20,9)-- (19,8);
\draw [line width=0.4pt] (20,9)-- (20,8);
\draw [line width=0.4pt] (20,9)-- (21,8);
\draw [fill=black] (20,10) circle (2pt);
\draw [fill=black] (20,9) circle (2pt);
\draw [fill=black] (19,8) circle (2pt);
\draw [fill=black] (20,8) circle (2pt);
\draw [fill=black] (21,8) circle (2pt);
%--T5,3---------------------------------------
\draw [line width=0.4pt] (23,10)-- (23.6,9);
\draw [line width=0.4pt] (23.6,9)-- (23.6,8);
\draw [line width=0.4pt] (23,10)-- (22,9);
\draw [line width=0.4pt] (23,10)-- (22.8,9);
\draw [fill=black] (23,10) circle (2pt);
\draw [fill=black] (23.6,9) circle (2pt);
\draw [fill=black] (23.6,8) circle (2pt);
\draw [fill=black] (22,9) circle (2pt);
\draw [fill=black] (22.8,9) circle (2pt);
%--T5,4---------------------------------------
\draw [line width=0.4pt] (25.2,10)-- (24.7,9);
\draw [line width=0.4pt] (25.2,10)-- (25.7,9);
\draw [line width=0.4pt] (25.7,9)-- (25.2,8);
\draw [line width=0.4pt] (25.7,9)-- (26.2,8);
\draw [fill=black] (25.2,10) circle (2pt);
\draw [fill=black] (24.7,9) circle (2pt);
\draw [fill=black] (25.7,9) circle (2pt);
\draw [fill=black] (25.2,8) circle (2pt);
\draw [fill=black] (26.2,8) circle (2pt);
%--T5,5---------------------------------------
\draw [line width=0.4pt] (28,10)-- (28,9);
\draw [line width=0.4pt] (28,9)-- (27.2,8);
\draw [line width=0.4pt] (28,9)-- (28.8,8);
\draw [line width=0.4pt] (28.8,8)-- (28.8,7);
\draw [fill=black] (28,10) circle (2pt);
\draw [fill=black] (28,9) circle (2pt);
\draw [fill=black] (27.2,8) circle (2pt);
\draw [fill=black] (28.8,8) circle (2pt);
\draw [fill=black] (28.8,7) circle (2pt);
%--T5,6---------------------------------------
\draw [line width=0.4pt] (30.5,10)-- (30.5,9);
\draw [line width=0.4pt] (30.5,9)-- (30.5,8);
\draw [line width=0.4pt] (30.5,8)-- (29.7,7);
\draw [line width=0.4pt] (30.5,8)-- (31.3,7);
\draw [fill=black] (30.5,10) circle (2pt);
\draw [fill=black] (30.5,9) circle (2pt);
\draw [fill=black] (30.5,8) circle (2pt);
\draw [fill=black] (29.7,7) circle (2pt);
\draw [fill=black] (31.3,7) circle (2pt);
%--T5,7---------------------------------------
\draw [line width=0.4pt] (33,10)-- (32.2,9);
\draw [line width=0.4pt] (33,10)-- (33.8,9);
\draw [line width=0.4pt] (32.2,9)-- (32.2,8);
\draw [line width=0.4pt] (33.8,9)-- (33.8,8);
\draw [fill=black] (33,10) circle (2pt);
\draw [fill=black] (32.2,9) circle (2pt);
\draw [fill=black] (33.8,9) circle (2pt);
\draw [fill=black] (32.2,8) circle (2pt);
\draw [fill=black] (33.8,8) circle (2pt);
%--T5,8---------------------------------------
\draw [line width=0.4pt] (35.5,10)-- (34.9,9);
\draw [line width=0.4pt] (35.5,10)-- (36.1,9);
\draw [line width=0.4pt] (36.1,9)-- (36.1,8);
\draw [line width=0.4pt] (36.1,8)-- (36.1,7);
\draw [fill=black] (35.5,10) circle (2pt);
\draw [fill=black] (34.9,9) circle (2pt);
\draw [fill=black] (36.1,9) circle (2pt);
\draw [fill=black] (36.1,8) circle (2pt);
\draw [fill=black] (36.1,7) circle (2pt);
%--T5,9---------------------------------------
\draw [line width=0.4pt] (38,10)-- (38,9);
\draw [line width=0.4pt] (38,9)-- (38,8);
\draw [line width=0.4pt] (38,8)-- (38,7);
\draw [line width=0.4pt] (38,7)-- (38,6);
\draw [fill=black] (38,10) circle (2pt);
\draw [fill=black] (38,9) circle (2pt);
\draw [fill=black] (38,8) circle (2pt);
\draw [fill=black] (38,7) circle (2pt);
\draw [fill=black] (38,6) circle (2pt);
%--T6,*---------------------------------------
\draw [line width=0.4pt] (40,10)-- (39.2,9);
\draw [line width=0.4pt] (40,10)-- (40.8,9);
\draw [line width=0.4pt] (40.8,9)-- (40.2,8);
\draw [line width=0.4pt] (40.8,9)-- (41.4,8);
\draw [line width=0.4pt] (41.4,8)-- (41.4,7);
\draw [fill=black] (40,10) circle (2pt);
\draw [fill=black] (39.2,9) circle (2pt);
\draw [fill=black] (40.8,9) circle (2pt);
\draw [fill=black] (40.2,8) circle (2pt);
\draw [fill=black] (41.4,8) circle (2pt);
\draw [fill=black] (41.4,7) circle (2pt);
%---------------------------------------------
%\begin{scriptsize}
\draw (0,11.3) node {$T_{2,1}$};
\draw (2,6) node {$T_{3,1}$};
\draw (4,11.3) node {$T_{3,2}$};
\draw (6,6) node {$T_{4,1}$};
\draw (9,11.3) node {$T_{4,2}$};
\draw (12,6) node {$T_{4,3}$};
\draw (14,11.3) node {$T_{4,4}$};
\draw (16.5,6) node {$T_{5,1}$};
\draw (20,11.3) node {$T_{5,2}$};
\draw (23,6) node {$T_{5,3}$};
\draw (25,11.3) node {$T_{5,4}$};
\draw (28,6) node {$T_{5,5}$};
\draw (30.2,11.3) node {$T_{5,6}$};
\draw (33,6) node {$T_{5,7}$};
\draw (35.3,11.3) node {$T_{5,8}$};
\draw (38,11.3) node {$T_{5,9}$};
\draw (41.1,11.3) node {$T_{6,*}$};
%\end{scriptsize}
\end{tikzpicture}

\bigskip
In the table every tree is represented by a parenthetical expression
with leaf values corresponding to the minimizer.

\smallskip
Recursive description of coding trees by parenthetical expressions: 

\smallskip
(i) $[w]$ is the one-vertex tree with vertex value $w$.  

\smallskip
(ii) If $T_1,\dots, T_m$ are parenthetical expressions 
for $m$ disjoint rooted trees, then the concatenated string $T_1\dots T_m$ denotes the corresponding forest.
If there are $k$ copies of the same tree $T_i$,
we use the abbreviation $T_i^k$ instead of repeating
the string $T_i$ $k$ times.

\smallskip
(iii) $F$ is the parenthetical expression for a forest, then
$[F]$ denoted the tree where the forest's member trees become the branches attached to the new root.

\begin{table}[h]
\begin{tabular}{c|c|c|c|c|c} 
\# & $\min$ & $\ell$ & $h$ & $w$ & tree \\ \hline
$T_{1,1}$ & $1.0$ & $1$ & $0$ & $1$ & $[1.0]$\\ 
%\end{tabular}
%
%\bigskip
\hline
%\medskip
%\begin{tabular}{c|c|c|c|c|c} 
%\# & $\min$ & $\lf$ & $\hh$ & $\wh$ & tree \\ \hline
$T_{2,1}$ & $2.0$ & $1$ & $1$ & $1$ & $[[1.0]]$\\ 
%\end{tabular}
%
%\bigskip
\hline
%\medskip
%\begin{tabular}{c|c|c|c|c|c} 
% & $\min$ & $\lf$ & $\hh$ & $\wh$ & tree \\ \hline
$T_{3,1}$ & $2.828427$ & $2$ & $1$ & $2$ & $[[0.707107]^2]$\\ 
$T_{3,2}$ & $3.0$ & $1$ & $2$ & $1$ & $[[[1.0]]]$\\ 
%\end{tabular}
%
%\bigskip
\hline
%\medskip
%\begin{tabular}{c|c|c|c|c|c} 
% & $\min$ & $\lf$ & $\hh$ & $\wh$ & tree \\ \hline
$T_{4,1}$ & $3.464102$ & $3$ & $1$ & $3$ & $[[0.577350]^3]$\\ 
$T_{4,2}$ & $3.779763$ & $2$ & $2$ & $2$ & $[[[0.629961]^2]]$\\ 
$T_{4,3}$ & $3.799605$ & $2$ & $2$ & $2$ & $[[0.671044][[0.819173]]]$\\ 
$T_{4,4}$ & $4.0$ & $1$ & $3$ & $1$ & $[[[[1.0]]]]$\\ 
%\end{tabular}
%
%\bigskip
\hline
%\medskip
%\begin{tabular}{c|c|c|c|c|c} 
% & $\min$ & $\lf$ & $\hh$ & $\wh$ & tree \\ \hline
$T_{5,1}$ & $4.0$ & $4$ & $1$ & $4$ & $[[0.5]^4]$\\ 
$T_{5,2}$ & $4.326749$ & $3$ & $2$ & $3$ & $[[[0.480750]^3]]$\\ 
$T_{5,3}$ & $4.401338$ & $3$ & $2$ & $3$ & $[[0.546097]^2 [[0.738984]]]$\\ 
$T_{5,4}$ & $4.457410$ & $3$ & $2$ & $2$ & $[[0.593905][[0.544933]^2]]$\\ 
$T_{5,5}$ & $4.729032$ & $2$ & $3$ & $2$ & $[[[0.569841][[0.754877]]]]$\\ 
$T_{5,6}$ & $4.756828$ & $2$ & $3$ & $2$ & $[[[[0.594604]^2]]]$\\ 
$T_{5,7}$ & $4.762203$ & $2$ & $2$ & $2$ & $[[[0.793701]]^2]$\\ 
$T_{5,8}$ & $4.787079$ & $2$ & $3$ & $2$ & $[[0.655866][[[0.868837]]]]$\\ 
$T_{5,9}$ & $5.0$ & $1$ & $4$ & $1$ & $[[[[[1.0]]]]]$\\ 
%\end{tabular}
%
%\bigskip
\hline
$T_{6,*}$ & $5.377468$ & $3$ & $3$ & $2$ & $[[0.571120][[0.484072][[0.695753]]]]$
\end{tabular}
\caption{Minima $m_T$ and minimizers for small trees;
here: $\ell$ --- number of leaves, 
$h$ --- height,
$w$ --- width (maximum cardinality of a level). 
}
\label{tab:treevalues}
\end{table}

%%%%%%%%%%%%%%%%%%%%%%%%%%%%%%%%
\medskip

Algebraic values for the leaf components of the minimizers are given below. They are listed in order of their numerical counterparts in Table~\ref{tab:treevalues}, but without inner parentheses. The symbol $x$ denotes the unique positive root of the corresponding irreducible polynomial $P(\cdot)$.

$T_{1,1}$, $T_{2,1}$, $T_{3,2}$, $T_{4,4}$, $T_{5,9}$: $[1]$.

$T_{k,1}$: $[(k-1)^{-1/2}]$ ($k=3,4,5$).

$T_{k,2}$: $[(k-2)^{-2/3}]$ ($k=4,5$).

$T_{k,3}$: $[x^2,x]$, $P(x)=(k-3)x^4+x^3-1$, ($k=4,5$).

$T_{5,4}$: $[2x^2, x]$, $P(x)=4x^4+4x^3-1$.

$T_{5,5}$: $[x^2,x]$, $P(x)=x^3+x^2-1$.

$T_{5,6}$: $[2^{-3/4}]$; 
$T_{5,7}$: $[2^{-1/3}]$.

$T_{5,8}$: $[x^3,x]$, $P(x)=x^6+x^4-1$.

$T_{6,*}$: $[x^3(x+1),x^2,x]$, $P(x)=x^4(x+1)^2(x^2+1)-1$.

\sol{prob:snowflake}
Consider a generalized snowflake with degree $p$ of the central vertex and the degrees $q+1$ of the middle ring vertices. In the standard snowflake, $p=6$, $q=3$.

Denote the center by $C$, the middle ring vertices by $B_i$ and the leaves by $A_{ij}$ ($i=1,\dots,p$, $j=1,\dots,q$). The minimizing vector will be denoted $\yy$. 

\smallskip
(a) Root at the center. Let $x=y(A_{ij})$ (the same for all, by symmetry). Then by the equations satisfied by the minimizer, $y(B_i)=qx$, $y(C)=pqx$, and $pq^2 x^3=1$.
Hence $x=(pq^2)^{-1/3}$ and 
$m_T=pq y(A_{ij})+py(B_i)+y(C)=3pqx=3(p^2q)^{1/3}$.

\smallskip
(b) Root at $B_1$. There are two types of leaves. Let
$y(A_{ij})=x$ for $i\neq 1$ and $y(A_{1j})=\xi$.
We have the additive relations
$y(B_i)=qx$ ($i\neq 1$), $y(C)=(p-1)qx$, $y(B_1)=(p-1)qx+q\xi$. The maximal chain
product equations $\xi=x\cdot qx\cdot (p-1)qx$, $\xi\cdot((p-1)qx+q\xi)=1$ follow. Therefore
$x$ is the positive root of the 6-th degree equation $q^3(p-1)^2 x^4(1+q^2 x^2)=1$. 
The result is
$m_T=3(p-1)qx+2q\xi=(p-1)qx(4+2q^2x^2)$.  

\smallskip
(c) Root at $A_{1,1}$. Again, there are two type of leaves. Let $y(A_{ij})=x$ for $i\neq 1$ and $y(A_{1j})=\xi$ ($j\neq 1$). We have
$y(B_i)=qx$ ($i\neq 1$), $y(C)=(p-1)qx$, $y(B_1)=y(A_{11})=(p-1)qx+(q-1)\xi$. The product equations are $\xi=x\cdot qx\cdot (p-1)qx$ and
$\xi ((p-1)qx+(q-1)\xi)^2=1$. 
Hence $x$ satisfies the 9-th degree equation
$(p-1)^3 q^4 x^5 (1+q(q-1)x^2)^2=1$.
The result is $m_T=5(p-1)qx+2(q-1)\xi=(p-1)qx(5+3q(q-1)x^2)$.

\smallskip
For $p=6$, $q=3$ the numerical answers are:
$m_T\approx 14.28661$ in the case (a), $m_T\approx 12.68813$ in the case (b), and
$m_T\approx 12.18572$ in the case (c).

%%%%%%%%%%%%%%%%
\subsection*{Section~\ref{sec:extremal_problems}}
\addcontentsline{toc}{subsection}{Section~6}
% Extremal problems

%\subsection*{Section~\ref{ssec:trees-extremal}}
%\addcontentsline{toc}{subsubsection}{Section~6.1}
% Extremal problems for node-weighted rooted trees

%\sol{prob:min-palm-tree-eqneib}
\solL{prob:min-palm-tree}{A}\
Raising both parts of the equation
$
(n-\ell+1)\ell^{1/(n-\ell+1)}=(n-\ell)(\ell+1)^{1/(n-\ell)}
$
to the power $(n-\ell+1)(n-\ell)$, we obtain two factorizations of the same integer. Since $\ell$ and $\ell+1$ are coprime as well as $n-\ell$ and $n-\ell+1$,
we conclude that $(n-\ell)^{n-\ell+1}=\ell$
and $(n-\ell+1)^{n-\ell}=\ell+1$. That is, $k^{k+1}=(k+1)^k+1$, where $k=n-\ell$. Rewriting it as
$k=(1+1/k)^k+k^{-k}$, we see that $k=2$ is the only solution. It corresponds to the equality 
$m_{\Pl(10,8)}=m_{\Pl(10,9)}=6$.

\sol{prob:mintree-to-max} 
The node value vector $\yy=\xone$ satisfies the homogeneous constraints for any tree and yeilds $\sum y_v=|\VV(T)|$.
Therefore always $m_T\leq|\VV(T)|$. It is easy to see that $\yy=\xone$ satisfies the conditions of extremum \eqref{tree-minimizer}
if and only if $T$ is a linear tree. 

\sol{prob:mintree-to-max-ell}\
Let $T\in\mathcal{T}(n,\ell)$, $\ell\geq x$.
Put $y_\lambda=1/x$ for $\lambda\in\Lambda$,
$y_v=1$ for $\rho\neq v\in\VV\setminus\Lambda$
and $y_\rho=x$. The so defined node value vector $\yy$ 
satisfies the homogeneous constraints.  Therefore  $m_T\leq
\sum_{v\in\VV} y_v= \ell/x+(n-\ell-1)+x$.

%\sol{prob:mintree-to-max-sqrk}\
\solL{prob:mintree-to-max-ell}{A}\
It is easy to see that in the case $\ell=x^2$ with integer $x$ the vector $\yy$ introduced in the previous solution satisfies the conditions of extremum \eqref{tree-minimizer}
if and only if the tree $T$ consists of 
$x$ branches of the type $\Pl(m,x)$ (possibly with different values of $m$) attached to the common root.
% Amalgamated direct sum of posets 

\medskip\noindent
{\bf Remark.} Another case where the upper bound of Problem \bref{prob:mintree-to-max-ell} is sharp is 
the tree of height 1:
$n-\ell=1$, cf.\ $T_{k,1}$ in solution of \bref{prob:trees-analytical-numerical}.

\iffalse
Introduce also trees that look like ``palm bushes''.
Let $n=1+n_1+\dots+n_r$ 
and $\ell=\ell_1+\dots+\ell_r$.
The tree $\Pl(n_1,\dots,n_r; \ell_1,\dots,\ell_r)\in\mathcal{T}(n,\ell)$
is defined as the union of the trees $\Pl(n_j+1,\ell_j)$,
$j=1,\dots,r$, that share the common root.
%The trees $\Pl(5,3;2,2)$ and $\Pl(4,4;2,2)$ are depicted. 
\fi
%%%%

\noindent
\begin{minipage}{0.6\textwidth}
%\sol{prob:mintree-to-max-not-unique}
\solL{prob:mintree-to-max-ell}{B}\
Two non-isomorhic rooted trees maximizing $m_T$
in $\mathcal{T}(7,4)$ are depicted on the right. 
However they are isomorphic as unrooted trees. 
\end{minipage}
\hspace{1.5em}
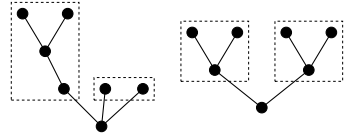
\begin{minipage}{0.18\textwidth}
\begin{tikzpicture}[x=0.5cm,y=0.5cm]
%--G(7,3)---------------------------------------
\draw [line width=0.4pt] (2.5,0)--(1.5,1)--(1,2);
\draw [line width=0.4pt] (0.4,3)--(1,2)--(1.6,3);
\draw [line width=0.4pt] (2.6,1)--(2.5,0)--(3.6,1);
\draw [fill=black] (2.5,0) circle (2pt);
\draw [fill=black] (1.5,1) circle (2pt);
\draw [fill=black] (1,2) circle (2pt);
\draw [fill=black] (0.4,3) circle (2pt);
\draw [fill=black] (1.6,3) circle (2pt);
\draw [fill=black] (2.6,1) circle (2pt);
\draw [fill=black] (3.6,1) circle (2pt);
\draw [line width=0.1pt, dash pattern=on 1pt off 1pt] (2.3,0.7)--(3.9,0.7)--(3.9,1.3)--(2.3,1.3)--cycle;
\draw [line width=0.1pt, dash pattern=on 1pt off 1pt] (0.1,0.7)--(1.9,0.7)--(1.9,3.3)--(0.1,3.3)--cycle;
\end{tikzpicture}
\end{minipage}
\hspace{-0.5em}
\begin{minipage}{0.2\textwidth}
\begin{tikzpicture}[x=0.5cm,y=0.5cm]
%--G(7,3)---------------------------------------
\draw [line width=0.4pt] (0,2)--(0.6,1)--(1.2,2);
\draw [line width=0.4pt] (2.5,2)--(3.1,1)--(3.7,2);
\draw [line width=0.4pt] (0.6,1)--(1.85,0)--(3.1,1);
\draw [fill=black] (0.6,1) circle (2pt);
\draw [fill=black] (0,2) circle (2pt);
\draw [fill=black] (1.2,2) circle (2pt);
\draw [fill=black] (2.5,2) circle (2pt);
\draw [fill=black] (3.1,1) circle (2pt);
\draw [fill=black] (3.7,2) circle (2pt);
\draw [fill=black] (1.85,0) circle (2pt);
\draw [line width=0.1pt, dash pattern=on 1pt off 1pt] (-0.3,0.7)--(1.5,0.7)--(1.5,2.3)--(-0.3,2.3)--cycle;
\draw [line width=0.1pt, dash pattern=on 1pt off 1pt] (2.2,0.7)--(4,0.7)--(4,2.3)--(2.2,2.3)--cycle;
\end{tikzpicture}
\end{minipage}

\medskip
A minimal example with non-isomorphic unrooted trees seems to be a pair of trees in $\mathcal{T}(13,9)$. 
Both trees have three branches $\Pl(n_1^{(i)},3)$, $\Pl(n_2^{(i)},3)$
and $\Pl(n_3^{(i)},3)$ ($i=1,2$) growing from the common root. The parameters are $n_1^{(1)}=n_2^{(1)}=n_3^{(1)}=5$
for the first tree and $n_1^{(1)}=4$, $n_2^{(1)}=5$, $n_3^{(1)}=6$ for the second tree.
As  unrooted trees, the first one has an automorphism of order 3 acting nontrivially on its non-leaves, while the second one does not. Therefore they are not isomorphic.

\sol{prob:mintree-to-max-lbnd}
Consider any tree $T\in\mathcal{T}(n,\ell)$ consisting
of $\ell$ linear branches attached to the root.
The vector $\yy$ with components $y_v=1/\sqrt{\ell}$
for any $v\neq\rho$ and $y_\rho=\sqrt{\ell}$ is the minimizer according to Eqs.\ \eqref{tree-minimizer}.
Therefore $m_T=(n-1)/\sqrt{\ell}+\sqrt{\ell}=(n-\ell+1)/\sqrt{\ell}$. This expression as a function of $\ell$ is decreasing for $\ell\leq n-1$.

\sol{prob:mintree-to-max-lbnd2}\
We will exhibit a tree $T\in\mathcal{T}(n,\ell)$ and
a vector $\xtt\leq\xone$ such that $m_T(\xtt)=2k+(n-\ell-1)$. By %\bref{prob:treemin-monotone-in-t} 
\brefA{prob:treemin-ineq}
it will follow that $m_{T}\leq 2k+(n-\ell-1)$.  

Consider a tree $T$ consisting of $a$ trees of the type $\Pl(*,k)$, $b$ trees of the type $\Pl(*,k+1)$, and $c$ trees of the type $\Pl(*,k+2)$ sharing the common root. 
The parameters $a$, $b$ and $c$ will be specified later.
Here $*$ stands for arbitrary integers, generally different, admissible by the construction.

Define the vector of node values $\yy$ as follows.
For the leaves and inner vertices (excluding the root)
of the subtrees $\Pl(*,r)$ we set $y_\lambda=1/r$,
$y_{v\notin\Lambda}=1$; here $r\in\{k,k+1,k+2\}$.
Finally, set $y_\rho=k$.

The additive conditions of extremum in \eqref{tree-minimizer} are satisfied if $a+b+c=k$.
The product conditions are satisfied if we assign the values of the maximal chain products as $t_\lambda=k/r$ for $\lambda\in\Lambda(\Pl(*,r))$.
Hence $\yy$ is the minimizer for $T$ with product constraints given by the vector $\xtt$. 
Clearly, $\xtt\leq\xone$. 

The number of leaves in $T$ is $ka+(k+1)b+(k+2)c=k^2+(b+2c)$, so we choose the parameters so as to make $\ell-k^2=b+2c$. The required number of nodes $n$ can be achieved by the choice of the heights of the palm-subtrees.

We have
$m_T(\xtt)=\sum_{\lambda\in\Lambda} y_\lambda+y_\rho+
\sum_{\rho\neq v\notin\Lambda} y_v=
2k+(n-\ell-1)$. 

Note that the equality $m_T=2k+(n-\ell-1)$ in this construction is possible only if $\xtt=\xone$,
that is, $b=c=0$, which means $\ell=k^2$.

The proof is complete.

%\sol{prob:mintree-to-max-doublebnd}
\solL{prob:mintree-to-max-lbnd2}{A}\
The double inequality is the combination of the lower bound just proved and the upper bound from
\bref{prob:mintree-to-max-ell}(i).

%\sol{prob:max-min-monotone} 
\solL{prob:mintree-to-max-lbnd2}{B}\
Fix $n$ and let $g(p)=\max_{T\in\mathcal{T}(n,k^2+p)} m_T-n-(k^2+p)-1$. The problem is to prove that
$g(p+1)<g(p)+1$ for $0\leq p\leq 2k$. 

The proof in the case $0\leq p<k$ is a simple consequence of the previous double-sided estimate.
Indeed, the lower bound implies that
$g(p)+1\geq 2k+1$ and the upper bound implies that $g(p+1)< 2\sqrt{k^2+p+1}\leq 2k+1$ if $p<k$.

For $k\leq p\leq 2k$ the same argument is sufficient to deduce a weaker monotonicity result:
$g(p+2)<g(p)+2$. 
However the author did not find a proof of the initial, stronger inequality.

\iffalse
Let $T\in\mathcal{T}(n,\ell)$ be a tree that maximizes $m_T$ in its class. Assuming that $\ell\geq 2$, we will exhibit a tree $T'\in\mathcal{T}(n,\ell-1)$ and a vector $\xtt\in\RR^{\ell-1}$ such that 
$\xone\geq\xtt$ and 
$m_{T'}(\xtt)>m_T$. Then by \bref{prob:treemin-monotone-in-t} it follows that
$m_{T'}=m_{T'}(\xone)>m_T$. 

Let $\lambda_*\in\Lambda(T)$ and $v_*=\alpha(\lambda_*)$.

Consider first the case where $d^+(v_*)>1$. 
Define the tree $T'$ as follows. 
Let $v_{**}=\alpha(v_*)$. 

1. Detach the branch $T_{v_*}$ from the node $v_{**}$.

2. Remove the leaf $\lambda_*$ from $T_{v_*}$.

3. Introduce a new node $v'$ and make it the parent of $v_*$. Call the so modified branch $T'_{v'}$.

4. Finally, attach $T'_{v'}$ to the node $v_{**}$ to get the tree $T'$. 

Supppose the vector $\yy=(y_v)_{v\in\VV(T)}$ is the minimizer for $T$. 
,,,
\fi

\sol{prob:mintree-to-min-local}
Suppose that $T$ has at least two branching nodes (i.e.\ nodes with outdegree $>1$) and a vector of values $\yy$ indexed by $\VV(T)$ with maximal path products equal to $t$ is given. We will construct a tree $\tilde T$ whose nodes will have the same labels and such that:

\smallskip
 (i) The number of branching nodes in $\tilde T$
is less than the number of branching nodes in $T$.

\smallskip
(ii) Maximal path products in $\tilde T$ corresponding to the vector of values $\yy$ mapped onto $T$ are  $\le t$.  
Then by 
\brefA{prob:treemin-ineq}
%\bref{prob:treemin-monotone-in-t} 
it follows that
$m_{\tilde T}(\xtt)\leq m_T(\xtt)$.

\smallskip
Let $v_1$ and $v_2$ be two branching nodes in $T$. There are two possibilities: (a) $v_1$ and $v_2$ are comparable (that is, there exists a maximal chain through both $v_1$ and $v_2$; (b) $v_1$ and $v_2$ are incomparable.

\smallskip
Case (a). Without loss of generality,  We may assume that $v_1< v_2$ and there are no branching nodes between $v_1$ and $v_2$. Denote $I=(v_1,v_2]$. Consider two subcases:

\smallskip
(a1) $P_{I}\ge 1$. The tree $\tilde T$ is constructed in two steps:
  
  1. Remove all nodes $v\in I$. (In other words re-plant the forest $\cup_{v\in T_{\alpha^{-1}(v_2)}} T_v$ at the root $v_1$ and remove the old trunk $I$.) 
 
  2. Attach the removed chain $I$ to any leaf of the former tree $T_{v_2}$.
    
\smallskip  
(a2)   $P_I< 1$. Again, the tree $\tilde T$ is constructed in two steps:

1. Cut off all branches outgoing from $v_1$ except the one containing $v_2$. (In other words, remove all subtrees $T_v$ with $\alpha(v)=v_1$ and $v\not\le v_2$.)

2. Re-plant the removed trees at the root $v_2$. 

\smallskip
  In both cases it is easy to see that the constructed tree $\tilde T$ satisfies our requirements.

\smallskip
Case (b). Since there are two incomparable non-leaves, the root is a branch node, and the situation is reduced to the case (a). 

%%%%%%%%%  
\sol{prob:mintree-to-min-special}\
By \bref{prob:mintree-to-min-local}, we may restrict attention to low-branching trees. The case of a linear tree is trivial. Suppose there is a branching point $v_*$ and two linear branches, $[v_1, \lambda_1]$ and $[v_2,\lambda_2]$, are attached to it.
Let $r_1$, $r_2$ be their respective lengths.
Let $\tau$ be the common value of the products $P_{[v_1,\lambda_1]}=P_{[v_2,\lambda_2]}$.  
The minimizing values of the nodes in the branch $[v_i,\lambda_i]$, $i=1,2$, 
are the constants $y^{(i)}=\tau^{1/r_i}$ (by the AM-GM inequality).
Hence the sum of the values in the two branches is
$$
 r_1 \tau^{1/r_1}+r_2\tau^{1/r_2}.
$$ 
If $|r_1-r_2|>1$, then it is possible to find integers $r_1'$ and $r_2'$ with $r_1'+r_2'=r_1+r_2$ and $|r_1'-r_2'|<|r_1-r_2|$.
Consider the function
$g(u,x)=xe^{u/x}$. We have
$$
 \frac{d^2 g(u,x)}{dx^2}=\frac{u^2}{x^3} e^{u/x}\ge 0.
$$ 
Thus, for any $u\in\RR$ the function $g(u,\cdot)$ is concave (linear, if $u=0$).
The concavity of $g\left(\ln\tau,\cdot\right)$ implies the inequality
$$
 r'_1 \tau^{1/r'_1}+r'_2\tau^{1/r'_2}<  r_1 \tau^{1/r_1}+r_2\tau^{1/r_2}.
$$ 
The existence of a minimizing tree of almost constant height follows.

\sol{prob:mintree-to-min-local-lb}
We may assume, by \bref{prob:mintree-to-min-special}, that $T$ is a low branching tree  of almost-constant height.
The case $\ell=1$ (the linear tree) is trivial, so suppose that $\ell\ge 2$, therefore there is a branching node; call it $v_*$.
Let $q$ be the number of nodes in the ``trunk'' $[\rho,v_*]$ (where $\rho$ is the root). The number of nodes in the chains
outgoing from $v_*$ (not including $v_*$) is therefore $n-q$ and the average length of those chains is $(n-q)/\ell$.
Let $\tau$ be the common value of the chain products for these chains. 

Let $\yy$ be the minimizer for the tree $T$.
As in the solution \bref{prob:mintree-to-min-special}, by  concavity of the function $g(\ln\tau,\cdot)$ we have:
$$
 \sum_{v>v_*} y_v\ge (n-q)\tau^{\frac{\ell}{n-q}}.
$$
Let $s$ be the common value of the nodes in $[\rho,v_*]$, i.e. $y_v=s$ for $v\in[\rho,v_*]$.
Then 
$t=\tau s^{q}$, so $s=(t/\tau)^{1/q}$. 
Therefore
$$
m_T(t\xone)=\sum_{v\in\VV(T)} y_v \ge qs+ (n-q)\tau^{\frac{\ell}{n-q}}=q\left(\frac{t}{\tau}\right)^{1/q}+(n-q)\tau^{\frac{\ell}{n-q}}.
$$
Putting $\alpha=1/q$, $\beta=\ell/(n-q)$, $A=qt^{1/q}$, $B=n-q$ in  a prepared formula
$$
 \min_\tau (A\tau^{-\alpha}+B\tau^\beta)=\left( \frac{(\alpha+\beta)^{\alpha+\beta}}{\alpha^\alpha \beta^\beta}\,A^\beta B^\alpha
 \right)^{\frac{1}{\alpha+\beta}},
$$
we find
$$
 m_T(t\xone) \ge \left(n+(\ell-1)q\right) \left(t\ell^{-q}\right)^{\frac{\ell}{n+(\ell-1)q}}.
$$
Put $n+(\ell-1)q=\xi$ and $\ell-1=r$. Then $q=(\xi-n)/r$.
In this notation, the inequality becomes
$$
 \ln m_T(t\xone)\ge \ln\xi+\frac{\ell}{\xi}\ln t +\frac{(n-\xi)\ell}{r\xi}\ln\ell =-\ell R(\ell) +\ln\xi +\frac{C}{\xi},
$$                
where 
$
C=\ell\left(\ln t+nR(\ell)\right).
$
The condition for the case (a) (formula \eqref{mintree-local-lb-small-t})
takes the form $C\leq n+\ell-1$ and the case (b)
takes the form $C\geq 0$.

The function $\xi\mapsto \ln\xi+C/\xi$ is increasing for all $\xi>0$ if $C<0$. If $C>0$, then it has minimum at $\xi=C$ and
$\min (\ln\xi+C\xi)=1+\ln C$.
Note that always $\xi\ge n+\ell-1$, since $q\ge 1$.

\medskip
Case (a). Since $\xi\geq n+\ell-1\geq C$, we have 
$$
\ln m_T(t\xone)\ge -\ell R(\ell) +\ln(n+\ell-1)+\frac{C}{n+\ell-1}
=\frac{\ell \,\ln(t/\ell)}{n+\ell-1}+\ln(n+\ell-1).
$$
This is the inequality \eqref{mintree-local-lb-small-t}. 
        
\medskip
Case (b). We have the inequality
$
\ln m_T(t\xone)\ge -\ell R(\ell) +1+\ln C$,
which is equivalent to \eqref{mintree-local-lb-large-t}.

\sol{prob:mintree-to-min-special-ge}\
Let $\yy$ be the minimizer  
for the tree $T$ with product data $t\xone$.
Let $\lambda_0\in\Lambda(T)$ be the leaf with minimum value.
%$
% y_{\lambda_0}=\min_{\lambda\in\Lambda(T)} y_{\lambda}
%$. 
%By the critical point equations \eqref{tree-minimizer}, 
By Eqs.\ \eqref{tree-minimizer},
$y_{\lambda_0}\le y_v$ $\forall v\in\VV(T)$.

Put $r=\left|\pi(\lambda_0)\right|-1$; it will be the length of the trunk in a palm tree to be constructed.
Let $k=n-r$. We will show that the inequality \eqref{palm-better-than-any} holds.
 
Let us set up a bijection $\phi:\VV(\Pl(n,k))\to\VV(T)$ as follows. We identify 
one leaf of $\Pl(n,k)$ with $\lambda_0$ and the path $\pi(\lambda_0)$ in $\Pl(n,k)$ with path $\pi(\lambda_0)$ in $T$.
(Both paths have the same cardinality by the choice of $k$).  
For the rest, the bijection $\phi$ between 
$\VV(\pi(\lambda_0))$ and $\VV(T)\setminus\pi(\lambda_0)$
is defined arbitrarily.

%we arbitrarily 
%put the leafs $\lambda\in %\Lambda(\Pl(n,k))\setminus\{\lambda_0\}$ in one-to-one correspondence with nodes $v\in \VV(T)\setminus\pi(\lambda_0)$.  

%The so defined bijection can be seen as labelling of nodes of $\Pl(n,k)$ by nodes of $T$.
Consequently, we have the vector $\yy$ mapped onto $\VV(\Pl(n,k)$. To distinguish the notation for components of $\yy$ and for path products in $T$ and $\Pl(n,k)$, the latter will be written  with tilde on top.   
By construction, $\sum_{v\in\VV(\Pl(n,k))}\tilde y_v=\sum_{v\in\VV(T)} y_v$.
% is the same for $\Pl(n,k)$ and $T$. 
It remains to check that 
all maximal chain products in $\Pl(n,k))$ do not exceed $t$. 
(We rely on the equivalence of minimization problems \eqref{mintree} and \eqref{mintree-ineq}.)

%To distinguish the notation for path products in $T$ and $\Pl(n,k)$, the latter will be written  with tilde on top.
By construction, $\tilde P_{\pi(\lambda_0)}=t$.  For any
other  $\lambda\in \Lambda(\Pl(n,k))$ we have the inequality, due to the choice of $\lambda_0$, 
$$
\tilde P_{\pi(\lambda)}=\tilde y_\lambda\frac{\tilde P_{\lambda_0}}{\tilde y_{\lambda_0}}=\frac{ty_{\phi(\lambda)}}{y_{\lambda_0}}\ge t.
$$

\medskip
{\em Uniqueness of the minimizing tree}. Let us explore  the case of equality $m_{\Pl(n,k)}(t\xone_k)=m_T(t\xone_\ell)$ in the above construction.
A necessary condition, in view of the last inequality,
is $y_{\phi(\lambda)}=y_{\lambda_0}$ for all $\lambda\in\Lambda(\Pl(n,k))\setminus\{\lambda_0\}$.
But this means that $y_v=y_{\lambda_0}$ for all $v\in\VV(T)\setminus\pi(\lambda_0)$. Therefore
$\VV(T)\setminus\pi(\lambda_0)\subset\Lambda(T)$. Hence
$T$ is a palm tree with trunk $[\rho,\alpha(\lambda_0)]$,
hence $T$ is isomorphic to $\Pl(n,k)$.

%\sol{prob:mintree-local-counterex} 
\solL{prob:mintree-to-min-special-ge}{A}\
Let $T_\ell$ be the unique tree of constant height 2 with $n=2\ell+1$ nodes and $\ell$ leaves.
It is easy to find the explicit solution of the system \eqref{tree-minimizer} with all $t_\lambda=1$:
$$
 y_v=\left\{\ba{ll}\ell^{-1/3}, & v\neq \rho,\\[0.5ex]
                           \ell^{2/3}, &v=\rho
                       \ea
                       \right.
$$
Hence $m_{T_\ell}=\ell^{-1/3}\cdot 2\ell+\ell^{2/3}=3\ell^{2/3}$.

On the other hand, by \eqref{min-palm-tree-homo}, $m_{\Pl(n,\ell)}>n>\ell/2$.
Hence $\Pl(n,\ell)$ certainly is not a minimizer for the homogeneous problem in $\mathcal{T}(n,\ell)$ if $\ell\ge 6^3$, i.e. if $n>432$. 

\medskip
\noindent
{\bf Remark.} We did not attempt to construct a counterexample with smallest possible $n$. Note that at least for $n\le 8$, as Table 1 shows,
the minimizing tree for the homogeneous problem in $\mathcal{T}(n,\ell)$ is always a palm tree.

\sol{prob:mintree-to-min-global-homo}\  
This follows from \bref{prob:mintree-to-min-special-ge} and formula~\eqref{min-palm-tree-homo}. 

%\sol{prob:mintree-to-min-global-explicit}
\solL{prob:mintree-to-min-global-homo}{A}\  
{\em Lower bound. }
Put $n+1=r$, $x=r-\ell$, and denote 
$$
 f(x,r)=\ln m_{\Pl(r+1,r-x)}=\ln x+\frac{\ln(r-x)}{x}.
$$
We will target the inequality $\exp f(x,r)>e\ln(r-\ln r)$, which is slightly stronger than \eqref{log-lower-bound-mT-trees}, which states that 
$\exp f(x,r)>e\ln(r-1-\ln (r-1))$

Let us temporarily forget that in our setting $x$ must be an integer and consider the function $x\mapsto f(x,r)$ with real $x\in (1,r-1)$. 
It has the unique critical point $x=\xi(r)$ determined by the equation
\beq{xi-critpt}
\xi=\frac{\ln (r-\xi)}{1-1/(r-\xi)}.
\eeq
Let $\xi=r\eta$, $0<\eta<1$. The above equation can be rewritten in the form suitable for iteration,
$$
 \eta=\frac{\ln r}{r}+\frac{g(\eta)}{r},
$$
where 
$$
 g(\eta)=\ln(1-\eta)+\frac{\eta}{1-\eta}=\frac{\eta^2}{2}+O\left(\eta^3\right).
$$
It follows that 
$$
 \eta=\frac{\ln r}{r}+\frac{(\ln r)^2+o(1)}{2r^3}
$$
and therefore
$$
\xi=\ln r+\frac{(\ln r)^2+o(1)}{2r^2}.
$$
Calculation (helped by expression of $\xi^{-1}\ln(r-\xi)$ via \eqref{xi-critpt}) then yields 
$$
 \min_x f(x,r)=f(\xi,r)=\ln\ln r+1-\frac{1}{r}-\frac{\ln r +o(1)}{2r^2}.
$$
On the other hand, 
$$
\ln (r-\ln r)=(\ln r)\left(1-\frac{1}{r}-\frac{\ln r +o(1)}{2r^2}\right),
$$
hence
$$
 \ln\ln(r-\ln r)=\ln\ln r-\frac{1}{r} -\frac{\ln r+1+o(1)}{2r^2}.
$$
We see that $\min_x f(x,r)>\ln\ln(r-\ln r)$ for sufficiently large $r$.
For small $r$ the inequality can be verified numerically (and we are only interested in integer values of $r$).

To round off the proof, one might establish analytically a lower value of ``sufficiently large'' $r$ by keeping track of particulars in the $o$-terms during
the derivation of the asymptotics.  
We do not go into it here, as \cite{KalachevSadov_2017} contains a compact but somewhat tricky proof valid for all $t>2$.
(Also proved there is a curious upper bound $\min_{1<x<r-1} f(x,r)<1+\ln\ln(r-\ln(r-\ln r))$.)

\medskip
{\em Asymptotic formula for $\min m_T$, $T\in\mathcal{T}(n)$. }
To determine the minimum in question we have to minimize $f(x,r)$ over {\em integer}\ values of $x$.
The optimal value of $x$ has the form $x_*=\xi(r)+O(1)$, where $\xi(r)$ is the critical point that we dealt with in the first part of the solution.
Obviously, the $O(1)$ accuracy is the best we can have if $x_*$ must be an integer.
Therefore $x_*=\ln r+O(1)$. Now, if we simply substitute this to the espression for $f(x_*,r)$, the result will be
$
f(x_*,r)=\ln\ln r+1+O\left(1/{\ln r})\right),
$ 
and  
$$
 \exp({f(x_*,r)})= e\ln r+O(1).
$$
The remainder here is worse than the one claimed in the problem. And indeed, we can do better with a little trick. Let us write $x_*=\ln r+\delta$
and work out the asymptotics of $f(x_*,r)$ with remainder in terms of $\delta$.   
Thus
$$
\ba{l}\dst
 \ln(\ln r+\delta)=\ln\ln r+\frac{\delta}{\ln r}+O\left(\frac{\delta^2}{(\ln r)^2}\right),
\\[2ex]
 \dst
  \ln (r-(\ln r+\delta))=(\ln r)\left(1+O\left(\frac{1}{r}\right)\right) ,
\\[2ex]
 \dst
 \frac{1}{\ln r+\delta}=\frac{1}{\ln r}\left(1-\frac{\delta}{\ln r}+O\left(\frac{\delta^2}{(\ln r)^2}\right)\right).
 \ea
$$
Putting everything together and only now substituting $\delta=O(1)$, we get
$$
f(\ln r+\delta,\,r)=\ln\ln r+1+O\left(\frac{1}{(\ln r)^2}\right).
$$
There follows 
$$
\exp\left(f(x_*,\,r\right)=e\ln r +O\left(\frac{1}{\ln r}\right),
$$
with remainder term as desired.

%\sol{prob:mintree-global-from-local}
\solL{prob:mintree-to-min-global-homo}{B}\  
{\bf Remark.} The crude lower bound in this problem
is obtained by evaluating the function $xe^{-x}$ at the left endpoint of the interval $[R(n-1),1]$, where it has minimum. It is wrong though to imply that the tree $\Pl(n,n-1)$ is the global minimizer in $\mathcal{T}(n)$. The more precise analysis in %\bref{prob:mintree-to-min-global-explicit} 
\brefA{prob:mintree-to-min-global-homo}
shows that the minimizing tree has not exactly $n-1$ but $n-O(n/\ln n)$ leaves. 

\sol{prob:randompalmbush}\
{\bf Comment.}\ One may want to explore the statistics of palm bushes and to compare it with statistics of general rooted trees. We quote some known results that may serve as benchmarks.  
%\iffalse

It is known \cite[Th.~3.5]{Moon_1970} that the number of all labelled, unordered rooted trees with $n$ vertices and $k\ge 2$ leaves is
$$
 R(n,k)=\frac{n!}{k!} S(n-2,n-k),
$$ 
where $S(n,k)$ denotes the Stirling numbers of the second kind \cite[\S~1.4]{Stanley_1999}. It is also known \cite[p.~6--7]{Czabarka_et_al} that for large $n$, the numbers $S(n,k)$ have approximately normal distribution with parameters
$$
\ba{l}\dst 
 \prE(S(n,\cdot))=\frac{B_{n+1}}{B_n}-1 \;\sim\frac{n}{\ln n},
 \\[2ex]
 \dst
  \mathcal{D}^2(S(n,\cdot))=\frac{B_{n+2}}{B_n}-\left(\frac{B_{n+1}}{B_n}\right)^2-1\;\sim\frac{n}{(\ln n)^2},
  \ea
$$
where $B_n$ are the Bell numbers.

Also in \cite[Corollary 4.2]{Czabarka_et_al} the asymptotic normality of $R(n,k)$ is stated. A ``typical'' labelled tree on $n$ nodes has about $n-(n/\ln n)(1+o(1))$ leaves. 

%As the author showed elsewhere, 
It is known
% since Polya (then E.M. Wright and others)
 that ``almost all'' labelled trees do not have nontrivial isomorphisms, so the distribution of unlabelled trees is roughly the same.

\iffalse
This prompts the question: how close is the mean value $\langle m_{T}\rangle_{T\in\mathcal{T}(n)}$
to the lower bound for $m_T$, $T\in\mathcal{T}(n,\ell)$ with $n-\ell\sim n/\ln n$ from \bref{prob:mintree-to-min-local-lb-homo},
which is $e\ln n(1+o(1))$? 
\fi

%\hardprob{prob:tree-average}
%Let $\mathcal{T}_{n,\ell}$ be the set of all labeled rooted trees with $n$ nodes and $\ell$ leaves, the root having fixed label $0$. Suppose that
%the leaves have fixed labels $n-\ell,\dots,n-1$ so that the leaf $\lambda_i$ has label $n-i$ and the  vector $\xtt=(t_i)_{i=1,\dots,\ell}$ is given,
%with components' indexation corresponding to the leaves' labels.     
%measure on

%\bigskip
%How about $\max m_T$ and $\min m_T$ if height of $T$ is restricted?
%\fi

\subsubsection*{Subsection~\ref{ssec:graphs-extremal}}
\addcontentsline{toc}{subsubsection}{Subsection~6.2}

\sol{prob:maxvalue-graph}
For any graph $\Ga\in\mathfrak{G}_m$, taking the vector of arc values $\xx=\xone$ we get the sum $m$, hence
$f_\Ga\le m$. And for the cycle graph $\Ga$ we have $f_\Ga=m$ by the AM-GM inequality.

\sol{prob:maxgraph}
Suppose that $f_\Ga=m$.  Then $\xx=\xone$ is a minimizer for $\Ga$. By uniqueness, it is {\em the}\ minimizer.  
From the necessary condition of extremum (first line in Eq.~\eqref{critpteqs}) we get $1=\sum_{L\in\LL} \lambda_L A_{L,a}$, $\forall a\in\AA(\Ga)$,
where $(A_{L,a})$ is the arc-circuit incidence matrix for some basis of circuits $\LL$ for $\Ga$.
%Therefore, if 
%\beq{potential_condition}
%\sum_a A_{L,a}u_a=0\quad
%\forall L\in\LL, 
%\eeq
%then $\sum_a u_a=0$.    
%    
%On the other hand, the condition \eqref{potential_condition} is equivalent to the existence of a potential function $\phi:\VV(\Ga)\to\RR$
%such that $u_a=\phi(\alpha(a))-\phi(\beta(a))$. We have
    
Consider an arbitrary function $\phi:\VV(\Ga)\to\RR$ and
let $u_a=\phi(\alpha(a))-\phi(\beta(a))$ for $a\in\AA(\Ga)$. Clearly, $\sum_{a\in L} u_a=0$ for any
circuit $L$. Therefore $\xuu\in\ker A$. Since $\xone\in\Ran A^T$, we have $(\xuu,\xone)=\sum_{a\in\AA} u_a=0$. But 
$$
 \sum_{a\in\AA(\Ga)} u_a=\sum_{a\in\AA(\Ga)} \phi(\alpha(a)) - \sum_{a\in\AA(\Ga)} \phi(\beta(a))
 =\sum_{v\in\VV(\Ga)} \phi(v) \left(d^+(v)-d^-(v)\right).
$$
The right-hand side must vanish for any $\phi$, therefore $d^+(v)=d^-(v)$ for all $v\in\VV(\Ga)$. This is a necessary and sufficient
condition for $\Ga$ to be an Eulerian graph. (Since $\Ga$ is known to be strongly connected.)

\sol{prob:mingraph-to-special}\ 
Let us use the interpretation of $f_\Ga$ as the minimum value of the sum of quotients, 
$f_\Ga=f^\div_\Ga=\min_{\yy} S_\Ga(\yy)$,
see \S\,\ref{ssec:sums-of-quotients} and %Problem~
\bref{prob:minquotient-mingraph}.

Suppose $\Ga\in\mathfrak{G}_{m,n}$ and $\yy$ is a minimizer for $S_\Ga(\cdot)$. 
We will construct a special graph $\tilde\Ga$ on the same set $\VV$ of nodes.
To construct $\tilde\Ga$ means to define its arcs.
Define the nodes $v_*\neq v^*$ so that
$$
  y_*\leq y_v\leq y^*\quad \forall v\in\VV.
$$

Step 1. Remove (if necessary) some arcs from $\AA(\Ga)$ to obtain the graph $\Ga_0$
(on the same set of nodes $\VV$) with the property:
there is exactly one path from $v_*$ to $v^*$, exactly one path from $v^*$ to $v_*$, and every node  $v\in\VV\setminus\{v_*,v^*\}$ 
lies on exactly one of these paths.

The graph $\tilde\Gamma$ will include $\Gamma_0$ as a subgraph, 
%: $\VV(\tilde\Gamma)=\VV(\Gamma_0)=\VV(\Gamma)$ and
that is,
$\AA(\tilde\Gamma)\supset\AA(\Gamma_0)$.)
Thus the graph $\tilde\Ga$ is guaranteed to be strongly connected.
 
\smallskip
Step 2. Let $r=|\AA(\Ga)\setminus\AA(\Ga_0)|$. We complete the construction of $\tilde\Ga$ by including $r$ distinct arcs from $v_*$ to $v^*$ in $\AA(\tilde \Ga)$. 

\smallskip
We have
$$
 S_\Ga(\yy)-S_{\tilde\Ga}(\yy)=\sum_{a\in\AA(\Ga)\setminus\AA(\Ga_0)} \left(\frac{y_{\alpha(a)}}{y_{\beta(a)}}-\frac{y_*}{y^*}\right)\ge 0.
$$
Since $S_\Gamma(\yy)=f_\Gamma$ and $S_{\tilde\Gamma}(\yy)\geq f_\Gamma$,
we conclude that $f_{\tilde\Gamma}\leq f_{\tilde\Gamma}$, as required.

\sol{prob:mingraph-to-tree} 
For  $n=m$ a check is straightforward. If $2\le n\le m-1$, then
the construction described in the solution of \bref{prob:tree-to-general} defines a one-to one correspondence between the set of low-branching trees
$\mathcal{T}^*(m,\ell)$ and the set of special graphs $\mathfrak{G}^*_{m,n}$.
For a graph $\Ga$ and a tree $T$ corresponding to each other, the equality $f_\Ga=m_T$ holds.

It remains to apply the results \bref{prob:mintree-to-min-local} and \bref{prob:mingraph-to-special}.

\solL{prob:mingraph-local-lb-homo}{B}\
This final result of the main theoretical thread of this work trivially follows from \brefA{prob:mintree-to-min-global-homo}. 

\subsubsection*{Subsection~\ref{ssec:teaser-minsum}}
\addcontentsline{toc}{subsubsection}{Subsection~6.3}

\sol{prob:extrval-teaser}
(a) There are various ways to construct an example. Here is one. Put
$$
\omega_i=\left\{
\ba{ll}
\{1,2,\dots,n-1\}\;&\mbox{\rm if}\;\; i=1;\\
\{n\}\;&\mbox{\rm if}\;\; 2\le i\le n,
\ea
\right.
$$
and
$$
x_i=\left\{
\ba{ll}
\eps^2\;&\mbox{\rm if}\;\; i=1;\\
\eps\;&\mbox{\rm if}\;\; 2\le i\le n-1;\\
1\;&\mbox{\rm if}\;\; i=n.\
\ea
\right.
$$ 
Assuming that $\eps<1$ we have $Y_\omega(\xx)=2+(n-2)\eps$.

For $n=2021$, taking $\eps=10^{-5}$, we get $2+(n-2)\eps=2.02019<2.021$.

\smallskip
Note that the described assignment of sets is reducible: take $\Sigma=\{n\}$.

\medskip
(b) 
Put $\omega(i) =\{i+1\}$ for $i = 1,\dots,k-1$; $\omega(k) = \{k+1,\dots,n\}$,
$\omega(i) = \{1\}$ for $i = k +1, \dots,n$.
(This assignment of sets is suggested by extremal properties of
palm trees considered in \S\,\ref{ssec:trees-extremal}.) 

Now take
$k = 6$, $x_i=3^{1-i}$ for $i=1,\dots,k$; $x_{k+1}=\dots=x_n=3^{-k}$, which yields $Y_\omega(\xx)=21-\frac{172}{729}\approx 20.764$.

\smallskip\noindent 
A similar, slightly worse example is obtained with
 $k = 7$, $x_i=(5/2)^{1-i}$ for $i=1,\dots,k$;  $x_{k+1}=\dots=x_n=(5/2)^{-k}$, 
which yields $Y_\omega(\xx)\approx 20.800$.
 
\medskip
(c) An assignment of sets $i\mapsto\omega(i)$ determines the graph $\Ga$ with $\VV(\Ga)=[1:n]$ and the set of arcs
$$
\AA(\Ga)=\{(i\to j), j\in\omega(i)\}.
$$
Irreducibility of $\omega$ is equivalent to $\Ga$ being strongly
connected.

If $(i_1,i_2,\dots, i_{r+1}=i_1)$ is any cycle of nodes in $\Ga$,
then we have
$$
\prod_{k=1}^r \frac{x_{i_k}}{\min(x_m\,|\,m\in\omega(i_{k}))}
\ge \prod_{k=1}^r \frac{x_{i_k}}{x_{i_{k+1}}}=1.
$$
Applying the inequality \eqref{log-lower-bound-graphs},
we get the estimate
$$
Y_{\omega}(\xx)\ge e\ln (n-\ln n)
$$
or, if one traces the origin of this result back to Problem
\ref{prob:mintree-to-min-global-homo}, the estimate
$$
Y_{\omega}(\xx)\ge \min_{1\leq\ell\leq n-1}(n-\ell+1)\ell^{\frac{1}{n-\ell+1}},
$$
which is less explicit but slightly more precise.

For $n=2021$ the former estimate yields $Y_{\omega}(\xx)\ge 20.679$ and the latter (with $\ell=8$ in the right-hand side)
$Y_{\omega}(\xx)\ge 20.704$.

\subsection*{Section~\ref{sec:shallit}}
\addcontentsline{toc}{subsection}{Section~7}

\sol{prob:Shallit}\
{\bf Comment.}\
The problem was published in the Problems section (later discontinued) of SIAM Reviews. A solution \cite{Shallit_sol_1995} appeared in due time. 
The numerical constant is
$$
 C=1.36945140\dots.
$$

As the first step, a reformulation of the problem in the with $O(n)$ summands was introduced; see
its interpretation in the current framework in Problem \bref{prob:Shallit-alt-graph}.

Shallit's problem thus has been referred to as solved, see e.g. \cite{Finch_2003}. However, the argument presented in \cite{Shallit_sol_1995} is not rigorous: some intermediate claims are based on numerical observations. 
A complete proof, containing justification of all steps of the argument in \cite{Shallit_sol_1995}, is given by the present author in \cite{Sadov_Shallit}. Also, a more precise asymptotic formula is proved there:
$$
\min S_n(\xx)-3n+C=O((2-\sqrt{3})^n).
$$

\sol{prob:Shallit-as-particular-case}
Let $\Ga_n$ be the graph with $\VV(\Ga_n)=[0:n]$ and the set of arcs 
$$
\AA(\Ga_n)=\{a_{i},i=1,\dots,n\}\cup \{b_{ij},1\leq i\leq j\leq n\},
$$
where
$$
a_i=(i-1\to i)
%$$
%and
%$$
,\qquad
b_{ij}=(j\to i-1).
$$
Assign the value $x_i$ to the arc $a_i$ and the value $p_{ij}=(x_i\dots x_j)^{-1}$ to the arc
$b_{ij}$. We get
$$
 S_n(\xx)=\sum_{i=1}^n x_i +\sum_{1\leq i\leq j\leq n} p_{ij}.
$$
The cycles $(a_i,a_{i+1},\dots,a_{j};b_{ij})$ form a basis of cycles in $\Gamma_n$ and all cyclic products are equal to $1$. Hence
$
  \min_{\xx>0}S_n(\xx)=f_{\Ga_n}
$.

\begin{figure*}
\begin{center}
% Defining arrowhead positions in the middle
\tikzset{->-/.style={color=blue, dash pattern= on 2pt off 2pt, decoration={
  markings,
  mark=at position .6 with {\arrow{Stealth}}},postaction={decorate}}}
\tikzset{-!>-/.style={color=blue, decoration={
  markings,
  mark=at position .56 with {\arrow{Stealth}}},postaction={decorate}}}
\tikzset{->>-/.style={color=red, decoration={
  markings,
  mark=at position .7 with {\arrow{Stealth}}},postaction={decorate}}}
\tikzset{->!-/.style={color=red, dash pattern= on 2pt off 2pt, decoration={
  markings,
  mark=at position .7 with {\arrow{Stealth}}},postaction={decorate}}}

%%%%%%%%%%%%%%
%% Graph for Problem 76
%%%%%%%%%%%%%%
\begin{tikzpicture}[x=1cm,y=1cm]
% Defining positions of nodes
\begin{scriptsize}
\node[label=left:0] (A) at (0,0) {};
\node[label=135:1] (B) at (0.5858,1.4142){};
\node[label=above:2] (C) at (2,2){};
\node[label=45:3] (D) at (3.4142,1.4142){};
\node[label=right:4] (E) at (4,0){};
\end{scriptsize}

%---------------------------------------
% forward edges and arrows
\draw[->>-] (A) to (B);
\draw[->>-] (B) to (C);
\draw[->>-] (C) to (D);
\draw[->>-] (D) to (E);

%---------------------------------------
% return arrows along the perimeter
\draw[->-] (B) to [bend right=22.5] (A);
\draw[->-] (C) to [bend right=22.5] (B);
\draw[->-] (D) to [bend right=22.5] (C);
\draw[->-] (E) to [bend right=22.5] (D);

%---------------------------------------
% inner return arrows 
\draw[-!>-] (E) to (A);
\draw[-!>-] (D) to (A);
\draw[-!>-] (C) to (A);
\draw[-!>-] (E) to (B);
\draw[-!>-] (D) to (B);
\draw[-!>-] (E) to (C);

% nodes
\filldraw  (A) circle (2pt);
\filldraw  (B) circle (2pt);
\filldraw  (C) circle (2pt);
\filldraw  (D) circle (2pt);
\filldraw  (E) circle (2pt);
\end{tikzpicture}
\hspace{3em}
%%%%%%%%%%%%%%
%% Graph for Problem 77
%%%%%%%%%%%%%%
\begin{tikzpicture}[x=1cm,y=1cm]
%%%%%%%%%

% Defining positions of nodes
\begin{scriptsize}
\node[label=left:0] (A) at (0,0) {};
\node[label={135:1}] (B) at (0.5858,1.4142){};
\node[label=above:2] (C) at (2,2){};
\node[label=45:3] (D) at (3.4142,1.4142){};
\node[label=right:4] (E) at (4,0){};
\end{scriptsize}

%---------------------------------------
% a-arrows
\draw[->>-] (A) to [bend left=15] (B);
\draw[->>-] (B) to [bend left=15] (C);
\draw[->>-] (C) to [bend left=15] (D);
\draw[->>-] (D) to [bend left=15] (E);

%---------------------------------------
% bar a-arrows
\draw[->!-] (B) to [bend left=15] (A);
\draw[->!-] (C) to [bend left=15] (B);
\draw[->!-] (D) to [bend left=15] (C);
\draw[->!-] (E) to [bend left=15] (D);

%---------------------------------------
% b-arrows 
\draw[-!>-] (A) to (C);
\draw[-!>-] (C) to (E);
\draw[->-] (D) to (B);

% nodes
\filldraw  (A) circle (2pt);
\filldraw  (B) circle (2pt);
\filldraw  (C) circle (2pt);
\filldraw  (D) circle (2pt);
\filldraw  (E) circle (2pt);
\end{tikzpicture}
\end{center}

\hspace*{1em}
\begin{minipage}[t]{0.45\textwidth}
%Left: 
The graph $\Gamma_4$ for {\bf 76}. In $\Gamma_n$ for {\bf 78},
dashed blue arrows are absent.
\end{minipage}
\hspace{3.5em}
\begin{minipage}[t]{0.4\textwidth}
%Right: 
The graph $\Gamma'_4$ for {\bf 77}. 
\end{minipage}

\end{figure*}

\sol{prob:Shallit-alt-graph}
Rewrite the definition of Shallit's sum in the form
$$
 S_n(\xx)-\sum_{i=1}^n x_i=\sum_{i\geq 1}\frac{1}{x_i}\left(1+\frac{1}{x_{i+1}}\left(1+\frac{1}{x_{i+2}}\left(\dots\left(1+\frac{1}{x_n}\right)\dots\right)\right)\right).
$$
Define the new variables $x'_i$ as the summands in the right-hand side. That is, put
$$
x'_1=\frac{1+x'_2}{x_1},\;\dots,
\;
x'_i=\frac{1+x'_{i+1}}{x_i},\;\dots,
\;
x'_n=\frac{1}{x_n}.
$$
Then 
$$
 x_i=\frac{1+x'_{i+1}}{x'_i}, \quad i=1,\dots,n-1,
\quad \mbox{and}\quad
x_n=\frac{1}{x'_n},
$$
so that
$$
 S_n(\xx)=\sum_{i=1}^n \left(x'_i+\frac{1}{x'_i}\right)+\sum_{i=1}^{n-1} \frac{x'_{i+1}}{x'_i}.
$$
The corresponding graph $\Gamma'_n$ has the set of nodes $\VV=[0:n]$ and the set of arcs
$$
\AA=\{a_{i},i=1,\dots,n\}\cup \{\bar{a}_{i},i=1,\dots,n\}
\cup\{b_{i},i=1,\dots,n-1\},
$$ 
where
$$
a_i=(i-1\to i),\qquad
\bar{a}_{i}=(i\to i-1),
$$
and
$$
b_{i}=\begin{cases}
(i-1\to i+1),\; i\equiv 1\!\!\mod 2,\\
(i+1\to i-1),\; i\equiv 0\!\!\mod 2.
\end{cases}
$$
The values are assigned to the arc as follows:
$$
 a_i\mapsto\begin{cases}x'_i,\; i=0\mod 2\\
                  1/x'_i,\; i=1\mod 2
            \end{cases},
\qquad
 \bar{a}_i\mapsto\begin{cases}1/x'_i,\; i=0\mod 2\\
                  x'_i,\; i=1\mod 2
            \end{cases},
$$
and
$$
b_i\mapsto \frac{x'_{i+1}}{x'_i}.
$$
A natural basis of cycles in $\Ga'_n$ is
$$
\ba{rccl}
\mathcal{L}&=&&\dst \left\{\{\bar{a}_{2i-1},\bar{a}_{2i}, b_{2i-1}\},\,1\leq i\leq \lfloor n/2\rfloor \right\}
\\[1ex] &&\cup&
 \left\{\{{a}_{2i},{a}_{2i+1}, b_{2i}\},\,1\leq i\leq \lfloor (n-1)/2\rfloor \right\}
\\[1ex] &&\cup&
\left\{\{a_i,\bar a_i\},\,i=1,\dots,n\right\}.
\ea
$$
Thus $|\VV|=n+1$, $|\AA|=3n-1$, $|\mathcal{L}|=2n-1$,
and Euler's formula \brefB{prob:cycle_basis_linalg}
%\bref{prob:Euler_formula} 
confirms that we
didn't miss anything.

It is trivial to check that the cycle products under the described assignment are all equal to 1, while the sum of the arc values
equals $S_n(\xx)$. We conclude that
$\min_{\xx} S_n(\xx)=f_{\Ga_n'}$, as required. 

\sol{prob:Shallit-no-double-edges}\
{\bf Comment.}\ The existence of a one-term asymptotics
$\min_{\xx}\hat S_n(\xx)=n\hat\lambda+o(n)$ follows from
\bref{prob:Shallit_general_pattern-asymp} with constant
$\hat\lambda\approx 2.48$, see
\bnref{prob:Shallit-root-part-cases}(b).

The existence of a two-term asymptotics was proved by the author by a method similar to the one used
in \cite{Sadov_Shallit} but the proof is unpublished; this is why this problem is marked with asterisk.

The constant term in the asymptotics is
$$
\hat C\approx 2.0112096.
$$

\sol{prob:Shallit_general_pattern-graph}
Let $\Ga_n(\Patt)$ be the graph with set of nodes $\VV=[0:n]$ and the set of arcs 
$$
\AA=\AA^{(1)}\cup\AA^{(2)},
$$
where
\beq{AA1}
\AA^{(1)}=\{(i-1\to i),\,i=1,\dots,n\}
\eeq
and
\beq{AA2}
\AA^{(2)}=\{(j\to i-1),\, 1\leq i\leq j\leq n,\,j-i\in \Patt\}.
\eeq
Assign the values $x_i$ to the arcs $a_i=(i-1\to i)$ from $\AA^{(1)}$ and the values $p_{ij}=(x_1\dots x_j)^{-1}$ to the arcs
$b_{ij}=(j\to i-1)$ from $\AA^{(2)}$. We get
$$
 S_n(\xx|\Patt)=\sum x_i +\sum p_{ij},
$$
and all cyclic products of arc values in $\Ga_n(\Patt)$ are equal to $1$. Thus
$$
  m_n(\Patt)=f_{\Ga_n(\Patt)}.
$$

\sol{prob:Shallit_general_pattern-asymp}\
Assign the special admissible set of values to the arcs of the graph $\Ga_n(\Patt)$, which will provide an upper bound for $f_{\Gamma_n(\Patt)}$:
\beq{genShallit-assignments-of-values}
 \ba{l}
 x_i=\rho_\Patt, \quad i=1,\dots, n,
 \\
 p_{ij}=\rho_\Patt^{i-j-1}, \quad j-i \in \Patt.
 \ea
\eeq

Let us omit the subscript $\Patt$ for brevity.
The sum of arc values \eqref{genShallit-assignments-of-values}
equals
\begin{multline*}
 n\rho+\sum_{\Patt\ni m\leq n-1} (n-m)\rho^{-m-1}
\\ 
 =(n+1)\left(\rho+\sum_{\Patt\ni m\leq n-1} \rho^{-m-1}\right)
 -\rho-\sum_{\Patt\ni m\leq n-1} (m+1)\rho^{-m-1}
\\
=(n+1)(\lambda+\tilde O(\rho^{-n}))-2\rho+o(1)
\\
=n\lambda+\lambda-2\rho+o(1).
\end{multline*}
Hence
$$
\lim_{n\to\infty} \frac{f_{\Gamma_n(\Patt)}}{n}\leq \lambda.
$$

To obtain a lower bound for $f_{\Gamma_n(\Patt)}$,
we will assign particular weights to the arcs of the graph $\Ga_n(\Patt)$ so as to make that the above defined values
minimize the weighted sum of arc values assuming all cycle products are 1.
To determine appropriate weights, let us first show that the following cycles form a basis $\LL$
in $\Ga_n(\Patt)$:
$$
 \gamma_{ij}=(i-1\to i\to\dots\to j-1\to j\to i-1)
=\{a_i,\dots,a_j,b_{ij}\} 
 , \quad j-i\in \Patt.
$$
It is enough to check that these cycles are independent. And it follows from Euler's formula
\brefB{prob:cycle_basis_linalg}.
%\bref{prob:Euler_formula}. 
Indeed, using notation \eqref{AA1}, \eqref{AA2}, we have $|\AA^{(1)}|=n$ and $|\LL|=|\AA^{(2)}|$, 
so $|\VV|-|\AA|+|\LL|=1$, as had to be expected. 

Put $w_{b_{ij}}=1$ for $b_{ij}\in\AA^{(2)}$. Then the arc values $p_{ij}$ can be identified with 
Lagrange multipliers for the cycles $\gamma_{ij}$. 
(The cycles are in one-to-one correspondence with arcs of the set $\AA^{(2)}$. )

The arc $a_i$ belongs to the cycle $\gamma_{k\ell}$ if and only
if $k\leq i\leq \ell$. 
We have the constraints $1\leq k\leq i$, $i\leq \ell\leq n$,
so the number of admissible cycles with fixed $m=\ell-k\in \Patt$
is $\min(m+1,i,n-m,n-i+1)$
provided $m\leq n-1$.

The critical point equations for the arc values $x_i$ will be satisfied if we define $w_{a_{i}}$ so as to satisfy the first part of the critical point system
\eqref{critpteqs},
$$
 x_i w_{a_i}=\sum_{j:\,j\leq n,\,j-i\in \Patt} p_{ij}, \quad i=1,\dots, n-1,
 $$ 
that is,
 $$
 \rho w_{a_i}=\sum_{m\in \Patt} c_{i,m}\rho^{-m-1}, 
 $$
 where, for the given $n$, 
 $$
  c_{i,m}=\min(m+1,\, i,\,n-m,\,n-i+1).
 $$
Using Eq.\ \ref{gen-Shallit-coef} we see that
$\rho w_{a_i}\leq \rho$.
 Hence, always $w_{a_i}\leq 1$. By \bref{prob:monotonicity},
 $$
 m_n(\Patt)= f_{\Ga_n(\Patt)}(\xone)\geq f_{\Ga_n(\Patt)}(\ww).
 $$
It is easy to see that for any given $\eps>0$ there exist
$k\in\NN$ and $n_0\in\NN$ such that
all but $k$ values $w_{a_i}$ are within the $\eps$-neighbourhood
of $1$ for any $n\geq n_0$. It follows that
$$
 f_{\Ga_n(\Patt)}(\ww)\geq (\lambda-\eps) n.
$$
Hence
$$
\lim_{n\to\infty}\frac{f_{\Ga_n(\Patt)}}{n}\geq \lambda.
$$

%%%%%%%%%%%%%%%%%%%%%%%%%%%%%%%%%%%
%\sol{prob:Shallit-root-part-cases}
\solL{prob:Shallit_general_pattern-asymp}{A}\
(a) In the Shallit's case the equation \eqref{gen-Shallit-coef}
takes the form 
$$
 \rho=\sum_{n=1}^\infty n\rho^{-n} = \rho^{-1}\left(1-\rho^{-1}\right)^{-2},
$$
whence $|1-\rho|=1$ and the unique positive root is $\rho_\Patt=2$.
In this case
$$
\lambda_\Patt=\rho_\Patt+\frac{1}{\rho_\Patt-1}=3.
$$

\smallskip
(b) In Problem~\bref{prob:Shallit-no-double-edges}
we have the equation
$$
\rho=\sum_{n=2}^\infty n\rho^{-n} =\frac{1}{\rho}\left(1-\frac{1}{\rho}\right)^{-2}-\frac{1}{\rho}.
$$
Putting $\rho+\rho^{-1}=\sigma$ and observing that
$\rho(1-\rho^{-1})^2=\sigma-2$, we obtain
$$
 \sigma=\frac{1}{\sigma-2},
$$
whence $\sigma=1+\sqrt{2}$. Solving the quadratic equation
$\rho^2-\sigma\rho+1=0$ and choosing the root $>1$ we get
$$
\rho_\Patt=\frac{1}{2}\left(1+\sqrt{2}+\sqrt{2\sqrt{2}-1}\right)
\approx 1.883203506.
$$
In this case
$$
\lambda_\Patt=\rho_\Patt+\frac{1}{\rho_\Patt-1}-\frac{1}{\rho_\Patt}.
$$
We have
$$
 (\rho_\Patt-1)(\rho_\Patt-\sigma+1)=\rho_\Patt^2-\sigma\rho_\Patt+\sigma-1=
 \sigma-2=\sqrt{2}-1.
$$
Hence $(\rho_\Patt-1)^{-1}=(\rho_\Patt-\sqrt{2})(\sqrt{2}+1)$.
Also we have $\rho_\Patt^{-1}=\sigma-\rho_\Patt$.
Thus
\begin{multline*}
\lambda_\Patt=(\rho_\Patt-\sqrt{2})(\sqrt{2}+1)+2\rho_\Patt-\sigma
=(3+\sqrt{2})\rho_\Patt-3-2\sqrt{2}
\\
=\frac{5+4\sqrt{2}}{2}-(3+2\sqrt{2})+\frac{3+\sqrt{2}}{2}
\sqrt{2\sqrt{2}-1}
\\
=\frac{1}{2}\left(-1+(3+\sqrt{2})
\sqrt{2\sqrt{2}-1}\right)
\approx
2.484435332.
\end{multline*}
%

%\sol{prob:Shallit_general_pattern_refined}\
\solL{prob:Shallit_general_pattern-asymp}{B}\
{\bf Comment.}\
It is not difficult to prove that if the two-term asymptotics exists for any {\em finite}\ pattern, then the same is true for any pattern. 

%The possibility of a negative answer cannot be excluded by ``handwaving''. The critical point equations reformulated in terms of a fixed point of a dynamical system may have ``limit cycles.''

\sol{prob:Shallit-Kn-random-directions}\
{\bf Comment.}\
The purpose of this problem is not so much the specific inquiry as an invitation to explore statistical properties of similar optimization problems.

\end{document}